\documentclass[11pt,lettersize,reqno,fullpage]{article}
\usepackage{amsmath,amssymb,amsfonts,mathrsfs,verbatim,enumitem,pstricks,amsthm,float,graphicx}
\usepackage[all]{xy}
\usepackage{centernot, xr, accents}
\usepackage{hyperref}

\theoremstyle{plain}
\newtheorem{thm}{Theorem}[section]
\newtheorem{cor}[thm]{Corollary}
\newtheorem{lmm}[thm]{Lemma}
\newtheorem{prp}[thm]{Proposition}
\newtheorem{que}[thm]{Question}
\newtheorem{eg}[thm]{Example}
\newtheorem{defn}[thm]{Definition}

\newtheorem{claim}[thm]{Claim}
 
\newtheorem{rem}[thm]{Remark} 

\newtheorem{notn}[thm]{Notation}
\newtheorem*{mthm}{MAIN THEOREM}
\newtheorem*{ack}{Acknowledgement}

\theoremstyle{definition}
\newtheorem{step}{Step}

\numberwithin{equation}{section}

\def \hf{\hspace*{0.5cm}}                      
\def\bge{\begin{equation}}                
\def\ede{\end{equation}}                
\def\bgd{\begin{displaymath}}         
\def\edd{\end{displaymath}}            
\def\bgee{\begin{equation*}}           
\def\edee{\end{equation*}}           

\def \ni{\noindent}

\def\lra{\longrightarrow}
\def\lrab{\dashrightarrow}
\def\lan{\langle}
\def\ran{\rangle}

\def\BA{\begin{eqnarray}}
\def\EA{\end{eqnarray}}
\def\BAA{\begin{eqnarray*}}
\def\EAA{\end{eqnarray*}}
\def\Bal{\begin{align*}}
\def\Eal{\end{align*}}

\def \C{\mathbb{C}}

\def \P{\mathbb{P}}
\def \A{\mathcal{A}}
\def \E{\mathcal{E}}
\def \Num{\mathcal{N}}
\def \PP{\mathcal{P}}
\def \f{\frac}
\def \ff{\tilde{f}}

\def \p{\tilde{p}}

\def \lp{l_{\p}}

\def \w{\tilde{w}}

\def \B{\mathcal{B}}

\def \UL{\mathbb{L}} 
\def \DL{\mathcal{L}}
\def \UV{\mathbb{V}}
\def \DV{\mathcal{V}}

\def \ds{\psi}
\def \us{\Psi}

\def \ct{\mathcal{C}}

\def \y{y}
\def \x{x}
\def \a{a}
\def \lm{\lambda}
\def \gD{\gamma_{\mathcal{D}}}
\def \gP{\gamma_{_{\mathbb{P}^2}}}
\def \G{\tilde{\gamma}}

\def \25node{A_5} 
\def \62node{A_6}

\def \D{\mathcal{D}}
\def \DD{\mathcal{D}}

\def \X{\mathfrak{X}}

\def \pf{\noindent \textbf{Proof:  }}
\def \ov{\overline}

\def \td{\PP D}

\def \XC{\mathcal{X}}

\def \AA{\hat{\A}}
\def \DD{\hat{\D}}
 
\def \XX{ \hat{\XC}}

\def \J{\mathcal{J}}
\def \U{\mathcal{U}}

\def \hati{\hat}
\def \hatii{\bar}

\def \W{\mathbb{W}}
\def \WL{\mathcal{W}}
\def \Q{\mathcal{Q}}

\def \N{\nabla}

\def \cu{\gamma}

\def \T{\mathrm{P}}
\def \Z{\mathrm{A}} 
\def \Y{\mathrm{B}}

\def \kq{\kappa}
\def \etq{\eta}

\def \Bq{\mathrm{B}}

\def \Hpl{\mathrm{H}} 
\def \Ln{\mathrm{L}}

\def \du{\sqcup}

\def \hxt2{\hat{x}_{t_2}} 
\def \hyt2{\hat{y}_{t_2}}

\def \FF{\mathrm{F}}
\def \Fx{\mathrm{F}_{x_{t_2}}}
\def \Fy{\mathrm{F}_{y_{t_2}}}

\def \xt{x_{t_2}} 
\def \yt{y_{t_2}}

\def \R{\mathcal{R}}

\def \mp{\PP \D_7^{s}}
\def \mq{\PP \D_8^{s}}
\def \mr{\PP \D_6^{\vee s}}

\title{{\LARGE Enumeration of curves with two singular points}}
\author{Somnath Basu and Ritwik Mukherjee }
\date{}
\usepackage{hyperref}
\hypersetup{
	colorlinks,
	citecolor=red,
	filecolor=black,
	linkcolor=blue,
	urlcolor=black
}

\begin{document}

\maketitle

\begin{abstract}
In this paper we obtain an explicit formula for the number of curves in $\P^2$, of degree $d$, 
passing through $(d(d+3)/2 -(k+1))$ generic points 
and having one node and one codimension $k$ singularity, where $k$ is at most $6$. 
Our main tool is a classical fact from 
differential topology: the number of zeros of a generic smooth section of a vector bundle $V$ over $M$, counted with a sign, is 
the Euler class of $V$ evaluated on the fundamental class of $M$. 
\end{abstract}

\tableofcontents

\section{Introduction}
\label{introduction}
\hf\hf Enumeration of singular curves in $\P^2$ (complex projective space) is a classical 
problem in algebraic geometry. In our first paper \cite{BM13}, we used purely topological 
methods to answer the following question: 
\begin{que}
How many degree $d$-curves are there in $\P^2$, passing through $(d(d+3)/2 -k)$ generic points and having  one singularity 
of codimension $k$, where $k$ is at most $7$?  
\end{que} 
In this paper, we extend the methods applied in \cite{BM13} to enumerate curves with two singular points. More precisely, we obtain an 
explicit answer for the following question: 
\begin{que}
How many degree $d$-curves are there in $\P^2$, passing through $(d(d+3)/2 -(k+1))$ generic points, having one node and one singularity 
of codimension $k$, where $k$ is at most $6$?  
\end{que} 
\hf\hf Let us denote the space of curves of degree $d$ in $\P^2$ by $\D$. It follows that $\D \cong \P^{\delta_d}$, where $\delta_d = d(d+3)/2$. Let $\gP\lra \P^2$ be the tautological line bundle. 
A homogeneous polynomial $f$, of degree $d$ and in $3$ variables, induces a holomorphic section of the line bundle $\gP^{*d} \lra \P^2$. 
If $f$ is non-zero, then we will denote its \textit{equivalence class} in $\D$ by $\ff$. Similarly, if $p$ is a non-zero vector in $\C^3$, 
we will denote its equivalence class in $\P^2$ by $\p$ \footnote{In this paper we will use the symbol $\tilde{A}$ to denote the equivalence class of 
$A$ instead of the standard $[A]$. This will make some of the calculations in section \ref{closure_of_spaces} easier to read.}.     
\begin{defn}
\label{singularity_defn}
Let $\ff \in \D$ and $\p \in \P^2$. A point $\p \in f^{-1}(0)$ \textsf{is of singularity type} $\A_k$,
$\D_k$, $\E_6$, $\E_7$, $\E_8$ or $\XC_8$ if there exists a coordinate system
$(x,y) :(\U,\p) \lra (\C^2,0)$ such that $f^{-1}(0) \cap \U$ is given by 
\begin{align*}
\A_k: y^2 + x^{k+1}   &=0  \qquad k \geq 0, \qquad \D_k: y^2 x + x^{k-1} =0  \qquad k \geq 4, \\
\E_6: y^3+x^4 &=0,  \qquad \E_7: y^3+ y x^3=0, 
\qquad \E_8: y^3 + x^5=0, \\
\XC_8: x^4 + y^4 &=0.   
\end{align*}
\end{defn}
In more common terminology, $\p$ is a {\it smooth} point of $f^{-1}(0)$ if 
it is a singularity of type $\A_0$; a {\it simple node} if its singularity type is $\A_1$; a {\it cusp} if its type is $\A_2$; a {\it tacnode} 
if its type is $\A_3$; a {\it triple point} if its type is $\D_4$; and a {\it quadruple point} if its type is $\XC_8$. \\
\hf\hf We have several results (cf. Theorem \ref{algoa1a1}\,-\,\ref{algope6a1}, section \ref{algorithm_for_numbers})  
which can be summarized collectively as our main result. Although \eqref{algoa1a1}-\eqref{algope6a1} may appear as equalities, 
the content of each of these equations is a theorem.
\begin{mthm}
\label{main_result}
Let $\X_k$ be a singularity of type $\A_k$, $\D_k$ or $\E_k$. Denote $\Num(\A_1\X_k,n)$ to be the number of 
degree $d$ curves in $\P^2$ that pass through $\delta_d - (k+1+n)$ generic points having  
one $\A_1$ node and  one singularity of type $\X_k$ at the intersection of $n$ generic lines. \\
\textup{(i)} There is a formula for $\Num(\A_1\X_k,n)$ if $k+1 \leq 7$, provided  $d \geq \ct_{\A_1\X_k}$ where 
\bgd
\ct_{\A_1\A_k} = k+3, ~~\ct_{\A_1\D_k} = k+1, ~~\ct_{\A_1\E_6} = 6. 
\edd
\textup{(ii)} There is an algorithm to explicitly compute these numbers.     
\end{mthm}

\begin{rem}
Note that $\Num(\A_1\X_k,n)$ is zero if $n>2$, since three or more generic lines do not intersect at a common point. Moreover, $\Num(\A_1\X_k,2)$ is the the number of curves, of degree $d$, that pass through $\delta_d-(k+3)$ generic points having one $\A_1$-node and one singularity of type $\X_k$ lying at a given fixed point (since the intersection of two generic lines is a point).
\end{rem}

In \cite{BM13}, we obtained an explicit formula for 
$\Num(\X_k, n)$, the number of degree $d$ 
curves in $\P^2$ that pass through $\delta_d - (k+n)$ generic points having  
one singularity of type $\X_k$ at the intersection of $n$ generic lines.
We extend the methods applied in \cite{BM13} to obtain an explicit 
formula for $\Num(\A_1 \X_k, n)$. \\[0.1cm]

\hf \hf The numbers $\Num(\A_1\X_k,0)$ till 
$k+1\leq 7$ have also been computed by Maxim Kazarian \cite{Kaz}   
using different methods. Our results for $n=0$ 
agree with his. The bound $d \geq \ct_{\A_1\X_k}$ is imposed to ensure that the 
relevant bundle sections are transverse.$\footnote{However, this bound is not the optimal bound.}$ 
The formulas for $\Num(\A_1\A_1,n)$ and $\Num(\A_1\A_2,n)$ also appear in \cite{Z1}. 
We extend the methods applied by the author to obtain the remaining formulas. 




\section{Overview} 
\hf\hf Our main tool will be the following well known fact from topology (cf. \cite{BoTu}, Proposition 12.8).
\begin{thm} 
\label{Main_Theorem} 
Let $V\lra X$ be a vector bundle over a manifold $X$. Then the following are true: 
\hspace*{0.5cm}\textup{(1)} A generic smooth section $s: X\lra V$ is transverse to the zero set. \\
\hspace*{0.5cm}\textup{(2)} Furthermore, if $V$ and $X$ are oriented with $X$ compact then the zero set of such a section defines an integer homology class  in $X$, 
whose Poincar\'{e} dual is the Euler class of $V$. 
In particular, if the rank of $V$ is same as the dimension of $X$, 
then the signed cardinality of $s^{-1}(0)$ is the Euler class of $V$, evaluated on the fundamental class of 
$X$, i.e., 
\bgd
|\pm s^{-1}(0)| = \lan e(V), [X] \ran. 
\edd
\end{thm}
\begin{rem}
Let $X$ be a compact, complex manifold, $V$ a holomorphic vector bundle and $s$ a holomorphic section that is transverse to the zero set. If the rank of $V$ is same as the dimension of $X$, 
then the signed cardinality of $s^{-1}(0)$ is same as its actual cardinality (provided $X$ and $V$ have their natural orientations). 
\end{rem}
However, for our purposes, the requirement that $X$ is a smooth manifold is too strong. We will typically be dealing with spaces that are 
smooth but have non-smooth closure. The following result is a stronger  version of Theorem \ref{Main_Theorem}, that applies to singular spaces, provided 
the set of singular points is of real codimension two or more.   
\begin{thm} 
\label{Main_Theorem_pseudo_cycle} 
Let $M \subset \P^{N}$ be a smooth, compact algebraic variety and $X \subset M$ a smooth subvariety, not necessarily closed. Let $V \lra M$ be an oriented vector bundle, such that the rank of $V$ is same as the dimension of $X$. Then the following are true: \\
\hspace*{0.5cm}\textup{(1)} The closure of $X$ is an algebraic variety and defines a homology class.\\
\hspace*{0.5cm}\textup{(2)} The zero set of a generic smooth section $s: M \lra V$ intersects $X$ transversely and does not intersect $\ov{X}-X$ anywhere. \\
\hspace*{0.5cm}\textup{(3)} The number of zeros of such a section inside $X$, counted with signs, 
is the Euler class of $V$ evaluated on the homology class $[\ov{X}]$, i.e., 
\bgd
|\pm s^{-1}(0) \cap \ov{X}| = |\pm s^{-1}(0) \cap X| = \big\langle e(V), ~[\ov{X}] \big\rangle.
\edd
\end{thm}

\begin{rem}
All the subsequent statements we make are true provided $d$ is sufficiently large. The precise bound on $d$ is given in section \ref{bundle_sections}. Although the results of this paper are an extension of \cite{BM13}, our aim has been to keep this paper self-contained. Ideally, a reader not familar with \cite{BM13} should have no difficulty following this paper.
\end{rem}

\hf\hf We will now explain our strategy to compute  $\Num(\A_1\X_k, n)$. 
The strategy is very similar to that of computing $\Num(\X_k, n)$, which 
was the content of \cite{BM13}. Let $\mathrm{X}_1$ and  $\mathrm{X}_2$ be two subsets of $\D \times \P^2$. Then we define 
\begin{align*}
\mathrm{X}_1 \circ \mathrm{X}_2 &:= \{ (\ff, \p_1, \p_2) \in  \D \times \P^2 \times \P^2: (\ff, \p_1) \in \mathrm{X}_1, ~~(\ff, \p_2) \in \mathrm{X}_2, ~~\p_1 \neq \p_2 \}. 
\end{align*}
Next, given a subset $\mathrm{X}$ of $\D \times \P^2$ we define 
\begin{align*}
\Delta \mathrm{X} &:= \{ (\ff, \p, \p) \in  \D \times \P^2 \times \P^2: (\ff, \p) \in \mathrm{X} \}. 
\end{align*}
Similarly, let $\mathrm{X}_1$ and  $\mathrm{X}_2$ be two subsets of $\D \times \P^2$ and $\D \times \P T\P^2$ respectively.
Then we define 
\begin{align*}
\mathrm{X}_1 \circ  \mathrm{X}_2 &:= \{ (\ff, \p_1, l_{\p_2}) \in  \D \times \P^2 \times \P T\P^2: (\ff, \p_1) \in \mathrm{X}_1, ~~(\ff, l_{\p_2}) \in \mathrm{X}_2, ~~\p_1 \neq \p_2 \}. 
\end{align*}
Finally, given a subset $\mathrm{X}$ of $\D \times \P T\P^2$ we define 
\begin{align*}
\Delta \mathrm{X} &:= \{ (\ff, \p, l_{\p}) \in  \D \times \P^2 \times \P T \P^2: (\ff, l_{\p}) \in \mathrm{X} \}. 
\end{align*}
\ni The following result is clear from the definition of closure.
\begin{lmm}
\label{closure_obvious}
We have the following equality of sets 
\begin{align*}
 \ov{\mathrm{X}_1 \circ \mathrm{X}_2} & = \ov{\ov{\mathrm{X}}_1 \circ \mathrm{X}_2} = \ov{\mathrm{X}_1 \circ \ov{\mathrm{X}}_2} = \ov{\ov{\mathrm{X}}_1 \circ \ov{\mathrm{X}}_2}. 
\end{align*}
\end{lmm}
\hf\hf Given a singularity $\X_k$, we also denote by $\X_k$, the \textit{space} of 
curves of degree $d$ with a marked point $\p$ such that the curve has a singularity of type $\X_k$ at $\p$. 
Similarly, $\A_1\circ \X_k$ is the the \textit{space} of degree $d$ curves with two distinct marked points $\p_1$ and $\p_2$ such that the curve has a node at $\p_1$ and a singularity of type $\X_k$ at $\p_2$. Note that except when $\X_k=\A_1$, the space $\A_1\circ\X_k$ {\it is} the fibre product $\A_1\times_\D \X_k$.\\
\hf\hf Let $\p_1, \p_2, \ldots , \p_{\delta_d -(k+1+n)}$ be 
$\delta_d-(k+1+n)$ generic points in $\P^2$ 
and $\Ln_1, \Ln_2, \ldots, \Ln_n $ be $n$ generic lines 
in $\P^2$. Define the following sets
\begin{align}
\label{hyperplane} 
\Hpl_i& := \{ \ff \in \D: f(p_i)=0 \}, \qquad \Hpl_i^* := \{ \ff \in \D: f(p_i)=0,  \nabla f|_{p_i} \neq 0 \},  \nonumber \\
\hat{\Hpl}_i& := \Hpl_i \times \P^2 \times \P^2, \qquad \hat{\Hpl}_i^* := \Hpl_i^* \times \P^2 \times \P^2 \qquad \textnormal{and}  \qquad \hat{\Ln}_i := \D \times\P^2 \times \Ln_i.  
\end{align}
By definition, our desired number $\Num(\A_1 \X_k,n)$ is the cardinality of 
the set
\bge
\label{number_Xk_defn}
\Num(\A_1\X_k,n) := |\A_1 \circ \X_k \cap \hat{\Hpl}_1 \cap \ldots \cap \hat{\Hpl}_{\delta_d-(n+1+k)} \cap 
\hat{\Ln}_1 \cap \ldots \cap \hat{\Ln}_n|.
\ede
Let us now clarify an important point to avoid confusion: as per our notation, 
the number $\Num(\A_1 \A_1,0)$ 
is the number of degree $d$ curves through $\delta_d-2$ generic points having 
two \textit{ordered} nodes. To find the the corresponding number of curves where the 
nodes are \textit{unordered}, we have to divide by $2$.\\
\hf\hf We will now describe the various steps involved to obtain an explicit formula for $\Num(\A_1\X_k,n)$. 

\begin{step} 
Our first observation is that if $d$ is sufficiently large then $\A_1 \circ \X_k$ is a smooth algebraic variety and its closure defines a homology class.
\begin{lmm}{\bf (cf. section \ref{bundle_sections})}
The space $\A_1 \circ \X_k$ is a smooth subvariety of $\D \times \P^2 \times \P^2$ of dimension $\delta_d-k$. 
\end{lmm}
\end{step}

\begin{step}
Next we observe that if the points and lines are chosen generically, 
then the corresponding hyperplanes and lines defined in \eqref{hyperplane} will intersect our space 
$\A_1 \circ \X_k$ transversely. Moreover, they would not intersect any extra points in the closure. 
\begin{lmm}
\label{gpl2}
Let $\p_1, \p_2, \ldots , \p_{\delta_d -(k+1+n)}$ be 
$\delta_d-(k+1+n)$ generic points in $\P^2$ 
and $\Ln_1, \Ln_2, \ldots, \Ln_n $ be $n$ generic lines 
in $\P^2$. Let $\hat{\Hpl}_i$, $\hat{\Hpl}_i^*$  and $\hat{\Ln}_i$ be as defined in 
\eqref{hyperplane}.  
Then 
\bgd
\ov{\A_1 \circ \X}_k \cap \hat{\Hpl}_1 \cap \ldots \cap \hat{\Hpl}_{\delta_d -(k+n+1)} 
\cap \hat{\Ln}_1 \cap \ldots \cap \hat{\Ln}_n = 
\A_1 \circ \X_k \cap \hat{\Hpl}_1^* \cap \ldots \cap \hat{\Hpl}_{\delta_d -(k+n+1)}^* 
\cap \hat{\Ln}_1 \cap \ldots \cap \hat{\Ln}_n  
\edd
and every intersection is transverse. 
\end{lmm}
We omit the details of the proof; it follows from an application of the families transversality theorem and Bertini's theorem. 
The details of this proof  can be found in \cite{BM_Detail}. 
 
\begin{notn}
\label{tau_bundle_defn}
Let $\gD\lra \D$ and $\gP\lra\P^2$ denote the tautological line bundles. If $c_1(V)$ denotes the first Chern class of a vector bundle then we set 
\bgd
\y : = c_1(\gD^*) \in H^{2}(\D; \mathbb{Z}), \qquad 
\a := c_1(\gP^*)  \in H^{2}(\P^2; \mathbb{Z}). 
\edd
\end{notn}
  
\hf\hf As a consequence of Lemma \ref{gpl2}  
we obtain the following fact: 
\begin{lmm}
\label{gpl}
The number $\Num(\A_1\X_k,n)$ is given by~  
$$ \Num(\A_1\X_k, n) =
\big\langle (\pi_{\D}^*\y)^{\delta_d-(n+k+1)} (\pi_2^*\a)^n, ~[\ov{\A_1 \circ \X}_k] \big\rangle$$
 where ~$\pi_{\D}, \pi_1, \pi_2: \D \times \P_1^2 \times \P^2_2 
 \lra \D, \P^2_1, \P^2_2 $
 are the projection maps.  
\end{lmm}
\end{step}
\pf This follows from  Theorem \ref{Main_Theorem_pseudo_cycle} and Lemma \ref{gpl2}. \qed \\

\ni As explained in \cite{BM13}, 
the space $\X_k$ is not easy to describe directly and hence 
computing $\Num(\X_k,n)$ \textit{directly} is not easy. As a result we define another space $\PP \X_k \subset \D \times \P T\P^2.$ This is the space of curves $\ff$, of degree $d$, with a marked point $\p \in \P^2$ and a marked direction $\lp \in \P T_{\p}\P^2$, such that the curve $f$ has a singularity of type $\X_k$ at $\p$ and certain directional derivatives \textit{vanish along $\lp$}, and certain other derivatives \textit{do not vanish}. To take a simple example, $\PP \A_2$ is the space of curves $\ff$ with a marked point $\p$ and a marked direction $\lp$ such that $f$ has an $\A_2$-node at $\p$ and the Hessian is degenerate along $\lp$, but the third derivative along $\lp$ is non-zero. It turns out that this space is much easier to describe. We have defined $\PP \X_k$ in section \ref{definition_of_px}. Similarly, instead of dealing with the space $\A_1 \circ \X_k$, we deal with the space $\A_1 \circ \PP \X_k$. 

\begin{step}
Next we observe that since $\A_1 \circ \PP \X_k$ is described locally as the vanishing of certain sections that are transverse to the zero set, these are \text{smooth} algebraic varieties.

\begin{lmm}{\bf (cf. section \ref{bundle_sections})}
\label{pr_sp_pseudo}
The space 
$\A_1 \circ \PP \X_k$ is a smooth subvariety 
of $\D \times \P^2 \times \P T\P^2$ of dimension $\delta_d-(k+1)$. 
\end{lmm}

\begin{notn}
\label{tau_bundle_pv_defn}
Let $\G \lra \P T\P^2$ be the tautological line bundle. The first Chern class of the dual will be denoted by $\lm := c_1(\G^*)\in H^{2}(\P T\P^2; \mathbb{Z})$.   
\end{notn}
\ni Lemma \ref{pr_sp_pseudo} now motivates the following definition: 
 \begin{defn}
\label{up_number_defn}
 We define the number 
$\Num(\A_1 \circ \PP \X_k,n,m)$ as 
\bge
\label{num_proj}
\Num(\A_1\PP \X_k, n,m) := \big\langle \pi_{\D}^*\y^{\delta_d-(k+n+m+1)} \pi_2^*\a^n \pi_2^*\lm^m, ~[\ov{\A_1 \circ \PP \X}_k] \big\rangle. 
\ede
where ~$\pi_{\D}, \pi_1, \pi_2: \D \times \P_1^2 \times \P T\P^2_2 \lra \D, \P^2_1, \P T\P^2_2 $ are the projection maps. 
\end{defn}

\ni The next Lemma relates the numbers $\Num(\A_1\PP\X_k,n,0)$ and 
$\Num(\A_1\X_k,n)$. 
\begin{lmm}
\label{up_to_down}
The projection map $ \pi: \A_1 \circ \PP \X_k \lra \A_1 \circ \X_k $ is one to one if $\X_k = \A_k, \D_k, \E_6, \E_7$ or $\E_8$ except for $\X_k = \D_4$ when it is three to one. In particular, 
\bge
\label{up_down_equation}
\Num(\A_1\X_k, n)  = \Num(\A_1 \PP \X_k, n,0) 
\qquad \textnormal{if} ~~\X_k \neq \D_4   \qquad \textnormal{and} 
\qquad \Num(\A_1 \D_4, n) = 
\frac{\Num(\A_1 \PP \D_4,n,0)}{3}.
\ede
\end{lmm}
\pf This is identical to the proof of the corresponding lemma in \cite{BM13}.\qed
\end{step}
\ni To summarize, the \textit{definition} of $\Num(\A_1\X_k,n)$ is \eqref{number_Xk_defn}. 
Lemma \ref{gpl} equates this number to a \textit{topological} computation. 
We then introduce another number $\Num(\A_1\PP\X_k,n,m)$ in Definition \ref{up_number_defn} and relate it to 
$\Num(\A_1\X_k,n) $ in Lemma \ref{up_to_down}. 
In other words, we do not compute $\Num(\A_1\X_k,n)$ \textit{directly}; 
we compute it \textit{indirectly} by first computing $\Num(\A_1\PP \X_k, n,m)$ and then using Lemma \ref{up_to_down}. \\
\hf\hf We now give a brief idea of how to compute these numbers. 
Suppose we want to compute $\Num(\A_1\PP \X_k, n, m)$. 
We first find some singularity $\X_l$ for which 
$\Num(\A_1\PP \X_l, n,m)$ has been calculated 
and which contains $\X_k$ in its closure, i.e., we want $\PP \X_k$ to be a subset of $\ov{\PP \X}_{l}$. 
Usually, $l=k-1$ but it is not necessary. Our next task is to 
describe the closure of $ \PP \X_l$ and   $\A_1 \circ \PP \X_{l}$ explicitly as 
\begin{align}
\ov{\PP \X}_{l}&=  \PP \X_l \du \ov{\PP \X}_k \cup \B_1  \qquad \textnormal{and} \label{stratification_general_one_point}  \\
\ov{\ov{\A}_1 \circ \PP \X}_{l} & = \ov{\A}_1 \circ  \PP \X_l \du \ov{\A}_1 \circ (\ov{\PP \X}_l- \PP \X_l) \du \big( \Delta \B_2 \big)  \nonumber \\ 
& = \ov{\A}_1 \circ  \PP \X_l \du \ov{\A}_1 \circ (\ov{\PP \X}_k \cup \B_1) \du \big( \Delta \B_2 \big), \qquad 
\textnormal{where} \label{stratification_general} \\
\Delta\B_2 & := \{ (\ff, \p_1, l_{\p_2}) \in \ov{\ov{\A}_1 \circ \PP \X}_l: \p_1 = \p_2 \}. \nonumber 
\end{align}
Note that $\ov{\A_1 \circ \PP \X_{l}} = \ov{\ov{\A}_1 \circ \PP \X_l}$.
The main content of \cite{BM13} was to express $\ov{\PP \X}_l$ as in \eqref{stratification_general_one_point}. 
The main content of this paper is to 
compute $\Delta \B_2$, i.e    
expressing $\ov{\A_1 \circ \PP \X}_l$ as in \eqref{stratification_general}.
Concretely, computing $\Delta \B_2$ means figuring out what happens to a $\X_k$ singularity, when it collides with an $\A_1$-node. 
As a simple example, when two nodes collide, we get an 

\begin{figure}[h!]
\vspace*{0.2cm}
\begin{center}\includegraphics[scale = 0.5]{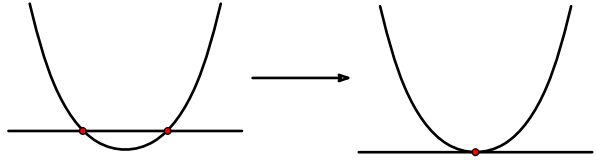}\vspace*{-0.2cm}\end{center}
\caption{Two nodes colliding into a tacnode}
\end{figure}

$\A_3$-node (which is basically the content of Lemma \ref{cl_two_pt}, statement \ref{a1a1_up_cl}). 
By definition \ref{up_number_defn} and Theorem \ref{Main_Theorem_pseudo_cycle} 
\bgd
\Num(\A_1\PP \X_k, n, m) := \big\langle e(\W_{n,m,k}^{1}), ~[\ov{\A_1 \circ \PP \X}_k] \big\rangle = \big|\pm \Q^{-1}(0) \cap \A_1 \circ \PP \X_k\big|, 
\edd
where
\bge
\label{generic_Q}  
\Q :\D \times \P^2 \times \P T\P^2  \lra \W_{n,m,k}^{\delta}:= \bigg({\textstyle \bigoplus}_{i=1}^{\delta_d -(n+m+k+\delta)}\pi_{\D}^*\gD^*\bigg)\oplus\bigg({\textstyle \bigoplus}_{i=1}^{n} 
\pi_1^*\gP^*\bigg)\oplus\bigg({\textstyle \bigoplus}_{i=1}^{m}\pi_2^*\G^*  \bigg)
\ede
is a generic smooth section. 
In \cite{BM13}, we constructed a section $\us_{\PP \X_k}$ 
of an appropriate vector bundle
\bgd
\UV_{\PP \X_k}\lra \ov{\PP \X}_{l} = \PP \X_{l} \cup \ov{\PP\X}_k \cup \B_1
\edd
with the following properties:  the section ~$ \us_{\PP \X_k}: \ov{\PP \X}_k \lra  \UV_{\PP \X_k}$
does not vanish on $\PP \X_{l}$ and it vanishes \textit{transversely} on $\PP \X_k$. 
With a similar reasoning one can show that the induced section 
$$ \pi_2^*\us_{\PP \X_k}: \ov{\A_1 \circ \PP \X}_l \lra  \pi_2^*\UV_{\PP \X_k} $$
vanishes \textit{transversely} on $\A_1 \circ \PP \X_k$.\footnote{However the bound on $d$ for which transversality 
is achieved increases.} 
Here $\pi_2$ is the following projection map 
$$ \pi_2: \D\times \P_1^2 \times \P T\P^2_2 \lra \D \times \P T \P^2_2.$$
Since $\us_{\PP \X_k}$ does not vanish on  
$\PP \X_{l}$, the section $ \pi_2^*\us_{\PP \X_k}$ does not 
vanish on $\A_1 \circ \PP \X_l$. 
Therefore, 
\begin{align}
\Big\langle e(\pi_2^*\UV_{\PP \X_{k}} \oplus \W_{n,m,k}^1), ~\big[\ov{\ov{\A}_1 \circ \PP \X}_{l}\big] \Big\rangle 
& = \Num(\A_1 \PP \X_k,n,m) + \mathcal{C}_{\A_1 \circ \B_1}(\pi_2^*\us_{\PP \X_k} \oplus \Q) \nonumber \\  
& \,\,+ \mathcal{C}_{\Delta \B_2}(\pi_2^*\us_{\PP \X_k} \oplus \Q) \label{Euler_equal_number_plus_bdry}
\end{align}
where $\mathcal{C}_{\A_1 \circ \B_1}(\pi_2^*\us_{\PP \X_k} \oplus \Q)$ and  $\mathcal{C}_{\Delta \B_2}(\pi_2^*\us_{\PP \X_k} \oplus \Q)$ 
are the contributions of the section
$\pi_2^*\us_{\PP \X_k} \oplus \Q$ to the Euler class from the points of $\A_1 \circ \B_1$ and $\Delta \B_2$ 
respectively.  
The number $\mathcal{C}_{\A_1 \circ \B_1}(\pi_2^*\us_{\PP \X_k} \oplus \Q)$ was computed 
in \cite{BM13}. The main content of this paper is to compute $\mathcal{C}_{\Delta \B_2}(\pi_2^*\us_{\PP \X_k} \oplus \Q)$.   
Once we have computed these numbers, we observe that the left hand side of \eqref{Euler_equal_number_plus_bdry} 
is computable via splitting principle and the 
fact that $\Num(\A_1\PP \X_l, n,m)$ is known. Therefore, we get a recursive formula for $\Num(\A_1 \PP \X_k,n,m)$ in terms of $\Num(\A_1 \PP \X_l,n^{\prime},m^{\prime})$ and $\Num(\PP \X_{k+1}, n, m)$. The main result of  
\cite{BM13} was to find an explicit formula for $\Num(\PP \X_{k+1}, n, m)$. Using this and iterations, 
we get an explicit formula for $\Num(\A_1 \PP \X_k,n,m)$. 
Finally, using 
Lemma \ref{up_down_equation}, we get our desired numbers $\Num(\A_1\X_k, n)$.
\begin{eg}
Suppose we wish to compute $\Num(\A_1\A_5, n)$. 
This can be deduced from the knowledge of $\Num(\A_1\PP \A_5, n,m)$. 
The obvious singularities which have $\A_5$-nodes in its closure are $\A_4$-nodes. 
In order to analyze the space $\ov{\ov{\A}_1 \circ \PP \A}_4$, 
we infer (cf. Lemma \ref{cl_two_pt}, statement \ref{a1_pa5_cl}) that
\bgd
\ov{\ov{\A}_1 \circ  \PP \A}_4 = \ov{\A}_1 \circ  \PP \A_4 \du \ov{\A}_1 \circ (\ov{\PP \A}_4- \PP \A_4) \du 
\Big( \Delta \ov{\PP \A}_6 \cup \Delta \ov{\mp} \cup \Delta \ov{\PP \E}_6   \Big). 
\edd
By \cite{BM13} (cf. Lemma \ref{cl}, statement \ref{A3cl} we conclude that 
\bgd
\ov{\ov{\A}_1 \circ  \PP \A}_4 = \ov{\A}_1 \circ  \PP \A_4 \du \ov{\A}_1 \circ (\ov{\PP \A}_5 \cup \ov{\PP \D}_5 ) \du 
\Big( \Delta \ov{\PP \A}_6 \cup \Delta \ov{\mp} \cup \Delta \ov{\PP \E}_6   \Big). 
\edd
The corresponding line bundle $\UL_{\PP \A_5} \lra \ov{\PP \A}_4$ with a 
section $\us_{\PP \A_5}$ that does not vanish on $\PP \A_4$ and vanishes transversely on 
$\PP \A_5$ is defined in section \ref{summary_vector_bundle_definitions}. 
In section \ref{bundle_sections}, we indicate that 
\[ \pi_2^*\us_{\PP \A_5}: \ov{\ov{\A}_1 \circ \PP \A}_4 \lra \pi_2^* \UL_{\PP \A_5} \] 
vanishes 
on $\A_1 \circ \PP \A_5$ transversely. 
Let $\Q$ be a generic section of the vector  bundle 
\bgd
\W_{n,m,5}^1 \lra \D \times\P^2 \times \P T\P^2. 
\edd
By \cite{BM13} 
$\pi_2^*\us_{\PP \A_5} \oplus \Q $ vanishes on all points of $\A_1 \circ \PP \D_5$ with a multiplicity 
of $2$. By Corollary \ref{a1_pak_mult_is_2_Hess_neq_0} and \ref{a1_pa4_mult_is_5_around_pe6},
$\pi_2^*\us_{\PP \A_5} \oplus \Q $ vanishes on all the points of $\Delta \PP \A_6$ and $\Delta \PP \E_6$ 
with a multiplicity of $2$ and $5$ respectively. 
Furthermore, we also show that $\Delta \ov{\mp}$ is contained inside $\Delta \ov{\PP \D}_7$. 
Since the dimension of $\Delta \PP \D_7$ is one less than the rank of $\pi_2^* \UL_{\PP \A_5} \oplus \W_{n,m,5}^1$ 
and 
$\Q$ is generic, $\pi_2^*\us_{\PP \A_5} \oplus \Q $ does 
not vanish on $\Delta \ov{\PP \D}_7$. Hence, it does not vanish on  $\Delta \ov{\mp}$. 
Therefore, we conclude that
\begin{align}
\big\langle e(\pi_2^*\UL_{\PP \A_5} \oplus \W_{n,m,5}^1), ~~[\ov{\ov{\A}_1 \circ \PP \A}_4] \big\rangle =  & \,\,\Num(\A_1\PP\A_5,n,m)\nonumber\\
& +2 \Num(\PP \D_5, n, m)+ 2 \Num(\PP \A_6, n, m)+ 5\Num(\PP \E_6, n, m).  \label{sample_computation}
\end{align}
This gives us a recursive formula for $\Num(\A_1 \PP \A_5, n, m)$ 
in terms of $\Num(\A_1 \PP \A_4, n^{\prime}, m^{\prime})$, $\Num(\PP \A_6, n, m)$, 
$\Num(\PP \D_5, n, m)$ and $\Num(\PP \E_6, n, m)$, 
which is \eqref{algopa4a1} in our algorithm. 
\end{eg}

\begin{rem}
We remind the reader that $\Num(\PP\X_k, n,m)$ has been defined in \cite{BM13}.
The definition 
is analogous to the definition of  
$\Num(\A_1 \PP \X_k, n, m)$ as given in definition \ref{up_number_defn} in this paper. 
\end{rem}

\hf\hf Now we describe the basic organization of our paper. In section \ref{algorithm_for_numbers} 
we state the explicit algorithm to obtain the numbers $\Num(\A_1\X_k, n)$ in our MAIN THEOREM. 
In section \ref{summary_notation_def} we recapitulate all the spaces, vector bundles 
and sections of vector bundles we encountered in the process of enumerating curves with one singular point. 
In section \ref{bundle_sections} we introduce 
a few new notation needed for this paper and 
write down the relevant sections that are transverse to the zero set. 
The proof of why the sections are transverse to the zero set can be found in \cite{BM_Detail}. 
In section \ref{closure_of_spaces} we stratify the space $\ov{\ov{\A}_1 \circ \PP \X}_k$ as described in \eqref{stratification_general}. Along the way we also compute the \textit{order} to which a certain section vanishes around certain points (i.e., the contribution of the section to the Euler class of a bundle). Finally, using the splitting principal, in section \ref{Euler_class_computation} we compute the 
Euler class of the relevant bundles and obtain the recursive formula similar to \eqref{sample_computation} above. 

\begin{ack}\nonumber 
\textup{One of the crucial results of this paper is to compute the closure of relevant spaces, i.e.,  
to stratify the space $\ov{\ov{\A}_1 \circ \PP \X}_k$ as described in \eqref{stratification_general}. A key step here is to observe that certain sections are transverse to the zero set and utilize them to describe the neighborhood of a point. The second author is indebted to Aleksey Zinger for sharing his understanding of transversality and explaining this crucial idea to him (i.e., how to  describe the neighborhood of a point using transversality of bundle sections). In addition, the second author is also grateful to Aleksey Zinger for suggesting several non trivial low degree checks to verify our formulas. One of those low degree checks proved to be crucial in figuring out a mistake the second author had made earlier.  \\ 
\hf \hf The authors are grateful to Dennis Sullivan for sharing his perspective on this problem and indicating its connection to other areas of mathematics.}
\end{ack}


\section{Algorithm}  
\label{algorithm_for_numbers}
\hf\hf We now give an algorithm to compute the numbers $\Num(\A_1\X_k,n)$. 
Equations \eqref{algoa1a1}-\eqref{algope6a1} are recursive formulas for $\Num(\A_1\PP\X_k,n,m)$ in terms of 
$\Num(\A_1 \PP \X_{k-1}, n^{\prime}, m^{\prime})$ and  
$\Num(\PP\X_{k+1},n, m)$. In \cite{BM13} we had obtained
an explicit formula for $\Num(\PP\X_{k+1},n, m)$. Finally, using 
Lemma \ref{up_down_equation}, we get our desired numbers $\Num(\A_1\X_k, n)$. We have implemented this algorithm in a Mathematica program to obtain the final answers. The program is available on our web page \url{https://www.sites.google.com/site/ritwik371/home}. We prove the formulas in section \ref{Euler_class_computation}.  \\
\hf \hf First we note that 
using the ring structure of $H^*(\D \times \P^2 \times \P T\P^2; \mathbb{Z})$, it is easy to see that 
for every singularity  type $\X_k$ we have
\begin{align}
\Num(\A_1\PP \X_k,n,m)
& = -3\Num(\A_1\PP \X_k,n+1,m-1)
-3\Num(\A_1\PP \X_k,n+2,m-2)   \qquad\forall ~~m\ge2. \label{ringp}
\end{align}
We now give recursive formulas for $\Num(\A_1\A_1, n)$ 
and  $\Num(\A_1\PP \X_k, n, m)$:    
\begin{eqnarray}
\Num(\A_1\A_1, n) &=& \Num(\A_1,0) \times \Num(\A_1, n)   \nonumber \\ 
                  & & -\big( \Num(\A_1, n) + d \Num(\A_1, n+1)  + 3 \Num(\A_2, n)\big)  \label{algoa1a1} \\
\Num(\A_1\PP \A_2, n, 0) & = & 2 \Num (\A_1\A_1,n) + 2(d-3) \Num(\A_1\A_1,n+1)  \nonumber  \\ 
                         & & - 2\Num(\PP \A_3, n, 0)\label{algopa20a1}\\
\Num(\A_1\PP \A_2, n, 1) & = & \Num(\A_1\A_1,n) + (2d-9) \Num(\A_1\A_1,n+1) + (d^2-9d+18) \Num(\A_1\A_1,n+2) \nonumber \\ 
                     & & - 2\Num(\PP \A_3, n, 1) - 3 \Num(\D_4,n)  \label{algopa21a1}\\
\Num(\A_1\PP \A_3, n, m) & = & \Num(\A_1\PP \A_2, n, m ) +  3\Num(\A_1\PP \A_2, n, m+1) + d\Num(\A_1\PP \A_2, n+1, m) \nonumber \\ 
                         & &   -2\Num(\PP \A_4,n,m)   \label{algopa3a1}
\end{eqnarray}
\begin{eqnarray}
\Num( \A_1 \PP \A_4, n, m) & = &  2\Num( \A_1\PP \A_3, n, m) + 2\Num(\A_1\PP \A_3, n, m+1) + (2d-6)\Num( \A_1\PP \A_3, n+1, m) \nonumber \\ 
                           &   &  -2\Num(\PP \A_5, n,m) \label{algopa4a1} \\
\Num(\A_1\PP \A_5, n, m) & = & 3\Num(\A_1\PP \A_4, n, m) + \Num(\A_1\PP \A_4, n, m+1) + (3d -12)\Num( \A_1\PP \A_4, n+1, m)\nonumber\\
&     & -2\Num(\A_1\PP \D_5, n, m) -2\Num(\PP \A_6, n, m)- 5\Num(\PP \E_6, n, m) \label{algopa5a1}\\ 
\Num( \A_1\PP\A_6, n, m) & = & 4\Num(\A_1 \PP\A_5, n, m) +0\Num(\PP \A_1\A_5, n, m+1) +  (4d -18)\Num(\PP \A_1\A_5, n+1, m) \nonumber \\   
              &     & -4\Num(\A_1\PP \D_6, n, m) - 3\Num(\A_1\PP \E_6, n, m)\nonumber\\
              &     & -2\Num(\PP \A_7, n, m) -6\Num(\PP \E_7, n, m)  \label{algopa6a1}\\
\Num(\A_1\PP \D_4, n, 0) & = & \Num(\A_1\PP \A_3, n, 0) -2\Num(\A_1\PP \A_3, n, 1) + (d-6)\Num(\A_1\PP \A_3, n+1, 0) \nonumber \\
                         &   & -2\Num(\D_5,n) \label{algopd4a1}   \\ 
\Num(\A_1\PP \D_4, n, 1) & = & \Num(\A_1 \D_4, n, 0) + (d-9)\Num(\A_1 \D_4, n+1, 0)  \label{algopd4a1_lambda} \\
\Num(\A_1\PP \D_5, n, m) & = & \Num(\A_1\PP \D_4, n, m) + \Num(\A_1\PP \D_4, n, m+1) + (d-3)\Num(\A_1\PP \D_4, n+1, m) \nonumber \\ 
                         & & -2\Num(\PP \D_6,n,m) \label{algopd5a1} \\
\Num(\A_1\PP \D_6, n, m) & = & \Num(\A_1\PP \D_5, n, m) + 4\Num(\A_1\PP \D_5, n, m+1) + d\Num(\A_1\PP \D_5, n+1, m)  \nonumber \\  
                         & &   -2\Num(\PP \D_7, n, m) - \Num(\PP \E_7,n,m) \label{algopd6a1} \\
\Num(\A_1\PP \E_6, n, m) & = & \Num(\A_1\PP \D_5, n, m) -\Num(\A_1\PP \D_5, n, m+1) + (d-6)\Num(\A_1\PP \D_5, n+1, m) \nonumber  \\ 
                         & &   -\Num(\PP \E_7, n,m) \label{algope6a1} 
\end{eqnarray}


\section{Review of definitions and notations for one singular point} 
\label{summary_notation_def}
\hf\hf We recall a few definitions and notation from \cite{BM13} so that our paper is self-contained.  
\subsection{The vector bundles involved}\label{summary_vector_bundle_definitions}
\hf\hf The first three of the vector bundles we will encounter, the tautological line bundles, have been defined in notations \ref{tau_bundle_defn} and \ref{tau_bundle_pv_defn}. 
Let $\pi:\D\times \P T\P^2\lra \D\times\P^2$ be the projection map. 
\begin{rem}
\label{an1_again}
We will make the abuse of notation of usually omitting the pullback maps 
$\pi_{\D}^*$ and $\pi_{\P^2}^*$.   
Our intended meaning should be clear when we say, for instance, $\gD^* \lra \D \times \P^2$. 
However, we will not omit to write the pullback via $\pi^*$. 
\end{rem}
We have the following bundles over $\D\times\P^2$ :
\begin{eqnarray*}
\DL_{\A_0} &:= & \gD^*\otimes \gP^{*d} \lra \D \times \P^2 \\
\DV_{\A_1} &:= & \gD^*\otimes \gP^{*d} \otimes T^*\P^2 \lra \D \times \P^2 \\
\DL_{\A_2} &:= & (\gD^* \otimes \gP^{*d} \otimes \Lambda^2 T^*\P^2)^{\otimes 2} 
\lra \D \times \P^2 \\ 
\DV_{\D_4} &:= & \gD^*\otimes \gP^{*d} \otimes 
\textnormal{Sym}^2 (T^*\P^2 \otimes T^*\P^2) \lra \D \times \P^2 \\
\DV_{\XC_8} &:= &  \gD^*\otimes \gP^{*d} \otimes 
\textnormal{Sym}^3 (T^*\P^2 \otimes T^*\P^2 \otimes T^*\P^2) 
\lra \D \times \P^2 
\end{eqnarray*}
Associated to the map $\pi$ there are pullback bundles 
\begin{eqnarray*}
\UL_{\AA_0} &:= &\pi^* \DL_{\A_0} \lra \D \times \P T\P^2 \\
\UV_{\AA_1} &:= & \pi^{*} \DV_{\A_1} \lra \D \times \P T\P^2 \\
\UV_{\DD_4} &:= & \pi^{*} \DV_{\D_4} \lra \D \times \P T\P^2 \\
\UV_{\XX_8} &:= & \pi^{*} \DV_{\XC_8} \lra \D \times \P T\P^2 \\
\UV_{\PP \A_2} &:= & \G^*\otimes \gD^*\otimes \gP^{*d} \otimes \pi^* T^*\P^2 
\lra  \D \times \P T\P^2 \\
\UV_{\PP \D_5} &:= & \G^{*2}\otimes \gD^*\otimes \gP^{*d} \otimes \pi^* T^*\P^2 
\lra  \D \times \P T\P^2.
\end{eqnarray*}
Finally, we have 
\begin{align*}
\UL_{\PP \D_4} &:= (T\P^2/\G)^{*2} \otimes \gD^* \otimes \gP^{*d} 
 \lra \D \times \P T\P^2  \\ 
 \UL_{\PP \D_5} &:= \G^{*2} \otimes (T\P^2/\G)^* \otimes 
 \gD^* \otimes \gP^{*d} \lra \D \times \P T\P^2 \\
\UL_{\PP \D_5^{\vee}} &:= \G^{*2} \otimes (T\P^2/\G)^{*4} \otimes 
 \gD^{*2} \otimes \gP^{*2d} \lra \D \times \P T\P^2 \\
 \UL_{\PP \D_6^{\vee}} &:= \G^{*8} \otimes (T\P^2/\G)^{*4} \otimes 
 \gD^{* 5} \otimes \gP^{*5d} \lra \D \times \P T\P^2 \\
\UL_{\PP \E_6} &:= \G^{*} \otimes (T\P^2/\G)^{*2} \otimes 
 \gD^* \otimes \gP^{*d} \lra \D \times \P T\P^2\\
\UL_{\PP \E_7} &:=  \G^{*4} \otimes \gD^* \otimes \gP^{*d}
\lra \D \times \P T\P^2 \\
\UL_{\PP \E_8} &:=  \G^{*3} \otimes (T\P^2/\G)^* \otimes  \gD^* \otimes \gP^{*d}
\lra \D \times \P T\P^2 \\
\UL_{\PP \XC_8} &:= (T\P^2/\G)^{*3} \otimes \gD^* \otimes \gP^{*d} 
\lra \D \times \P T\P^2 \\
\UL_{\J} &:= \G^{*9}\otimes(T\P^2/\G)^{*3}\otimes \gD^{*3} \otimes \gP^{*d} 
\lra \D \times \P T\P^2 \\
k \geq 3 \qquad \UL_{\PP \A_k} &:= \G^{*k} \otimes (T\P^2/\G)^{*(2k-6)} \otimes \gD^{*(k-2)} \otimes \gP^{*(d(k+1)-3d)}  
\lra \D \times \P T\P^2 \\  
k \geq 6 \qquad \UL_{\PP \D_k} & := \G^{*(k-2 + \epsilon_k)} \otimes (T\P^2/\G)^{*(2\epsilon_k)} \otimes \gD^{*(1+\epsilon_k)} \otimes \gP^{*(d(1+ \epsilon_k))} \lra \D \times \P T\P^2,
\end{align*}
where $\epsilon_6 =0$, $\epsilon_7=1$ and $\epsilon_8 =3$. 
The reason for defining these bundles will become clearer in section \ref{summary_sections_of_vector_bundle_definitions}, when we define 
sections of these bundles.\\
\hf \hf With the abuse of notation as explained in Remark \ref{an1_again}, the bundle $T\P^2/\G$ is the quotient of the bundles $V$ and $W$, 
where $V$ is the pullback of the tangent bundle $T\P^2\to \P^2$ via $\D\times\P T\P^2\stackrel{\pi}{\rightarrow}\D\times \P^2\to \P^2$ and $W$ 
is pullback of $\G\to \P T \P^2$ via $\D\times \P T \P^2\to \P T \P^2$.

\subsection{Sections of Vector Bundles}\label{summary_sections_of_vector_bundle_definitions}
\hf\hf Let us recall the definition of \textit{vertical derivative}. 
\begin{defn}
\label{vertical_derivative_defn}
\ni Let $\pi:V\longrightarrow M$ be a holomorphic vector bundle of rank $k$ and $s:M\lra V$ be a holomorphic section. Suppose $h: V|_{\U} \longrightarrow \U\times \mathbb{C}^{k}$ is a holomorphic trivialization of $V$ and $\pi_{1}, \pi_{2}: \U \times \C^{k} \longrightarrow \U, \C^{k}$ the projection maps.
Let 
\begin{align}
\label{section_local_coordinate}
\hati{s}&:= \pi_{2} \circ h \circ s.
\end{align}
For $q\in \U$, we define the {\it vertical derivative} of $s$ to be the $\C$-linear map  
\begin{align*}
\N s|_{q}: T_{q}M \longrightarrow V_{q}, \qquad  
\N s|_{q} & := (\pi_{2} \circ h)|_{V_{q}}^{-1} \circ d \hati{s}|_q,
\end{align*}
where $V_{q} = \pi^{-1}(q)$, the fibre at $q$. In particular, if $v \in T_q M$  
is given by a 
holomorphic map $\cu: \mathrm{B}_{\epsilon}(0) \lra M$ such that $\cu (0)=q$ and $\frac{\partial \cu}{\partial z}\big|_{z=0} = v$, then
\begin{align*} 
\nabla s|_q (v) & := (\pi_{2} \circ h)|^{-1}_{V_{q}} \circ 
\frac{\partial \hati{s}(\cu(z))}{\partial z}\bigg|_{z=0}
\end{align*}
were $\mathrm{B}_{\epsilon}$ is an open $\epsilon$-ball in $\C$ around the origin.$\footnote{Not every tangent vector is given by a holomorphic map; 
however combined with 
the fact that $\nabla s|_p $ is $\C$-linear, this definition determines $\nabla s|_p$ completely.} $ 
Finally, if $v, w \in T_q M$  are tangent vectors such that  
there exists a family of complex curves  $\cu: \mathrm{B}_\epsilon\times \mathrm{B}_\epsilon \lra M$ such that 
\bgd
\cu (0,0) =q, \qquad 
 \frac{\partial \cu(x,y)}{\partial x}\bigg|_{(0,0)}= v, \qquad 
\qquad \frac{\partial \cu(x,y)}{\partial y}\bigg|_{(0,0)} = w
\edd
then  
\begin{align}
\label{verder}
\nabla^{i+j} s|_q 
(\underbrace{v,\cdots v}_{\textnormal{$i$ times}}, \underbrace{w,\cdots w}_{\textnormal{$j$ times}}) & :=  
(\pi_{2} \circ h)\mid^{-1}_{V_{q}} \circ\left[
\frac{ \partial^{i+j} \hati{s}(\cu(x,y))}{\partial^i x \partial^j y}\right]\bigg|_{(0,0)}.
\end{align}
\end{defn}
\begin{rem}
In general the quantity in \eqref{verder} is not well defined, i.e., it depends on the 
trivialization and the curve $\cu$. In \cite{BM_Detail} we explain on what subspace 
this quantity is well defined. 
\end{rem}
\begin{rem}
\label{transverse_local}
The section $ s: M \lra V $  is transverse to the zero set if and only if the induced map
\begin{align} 
\label{section_local_coordinate_calculus}
\hatii{s} & := \hati{s} \circ \varphi_{\U}^{-1} : \C^m \lra \C^k 
\end{align}
is transverse to the zero set in the usual calculus sense, where $\varphi_{\U}: \U \lra \C^m $ is a coordinate chart and $\hati{s}$ is as defined in \eqref{section_local_coordinate}.
\end{rem}
\hf\hf Let $f:\P^2 \lra \gP^{*d}$ be a section and $\p \in \P^2$. 
We can think of $p$ as a non-zero vector in $\gP$ and $p^{\otimes d}$ a non-zero vector 
in $\gP^{\otimes d}$ $\footnote{Remember that $p$ is an element of $\C^3-0$ while $\p$ is the corresponding equivalence class in $\P^2$.}$. 
The quantity $\nabla f|_{\p}$ acts on a vector in  $\gP^{d}|_{\p}$ and produces an element of $T^*_{\p}\P^2$ . Let us denote this quantity as $\nabla f|_p$, i.e., 
\begin{align}
\nabla f|_p &:= \{\nabla f|_{\p}\}(p^{\otimes d}) \in T^*_{\p}\P^2.   
\end{align}
Notice that  $\nabla f|_{\p}$ is an element of the fibre of $T^*\P^{2} \otimes \gP^{*d}$ at $\p$ while $\nabla f|_{p}$ is an element of $T^*_{\p}\P^{2}$. \\
\hf\hf Now observe that $\pi^{*} T\P^2 \cong \G \oplus \pi^*T\P^2/\G \lra \P T\P^2$, where 
$\pi: \P T\P^2 \lra \P^2~$ is the projection map. Let us denote a vector in $\G$ by $v$ and  a vector 
in $\pi^*T \P^2/\G$ by $\w$. 
Given $\ff \in \D$ and $\p \in \P^2$, let 
\begin{align}
\label{abbreviation}
f_{ij} & := \nabla^{i+j} f|_p 
(\underbrace{v,\cdots v}_{\textnormal{$i$ times}}, \underbrace{w,\cdots w}_{\textnormal{$j$ times}}).
\end{align}
Note that $f_{ij}$ is a \textit{number}. 
In general $f_{ij}$ is not well defined; it depends on the trivialization and  the curve. 
Moreover it is also not well defined on the quotient space. 
Since our sections are not defined on the whole space, 
we will use the notation $ s:M \lrab V$ to indicate that $s$ is defined only on 
a subspace of $M$. \\
\hf \hf With this terminology, we now explicitly define the sections that we will 
encounter. 
\begin{eqnarray*}
\ds_{\A_0}: \D \times \P^2 \lra \DL_{\A_0}, & & \{\ds_{\A_0}(\ff, \p)\} (f \otimes p^{\otimes d}):= f(p) \\
\ds_{\A_1}: \D \times \P^2 \lrab \DV_{\A_1}, & & \{\ds_{\A_1}(\ff, \p)\}(f \otimes p^{\otimes d}):= \nabla f|_p  \\
\ds_{\D_4}: \D \times \P^2 \lrab \DV_{\D_4}, & & \{\ds_{\D_4}(\ff, \p)\}(f\otimes p^{\otimes d}):= \nabla^2 f|_p  \\ 
\ds_{\XC_8}: \D \times \P^2 \lrab \DV_{\XC_8}, & & \{\ds_{\XC_8}(\ff, \p)\}(f \otimes p^{\otimes d}):= \nabla^3 f|_p  \\ 
\ds_{\A_2} : \D \times \P^2 \lrab \DL_{\A_2}, & & \{\ds_{\A_2}(\ff, \p)\}(f \otimes p^{\otimes d}):= \textnormal{det}\, \nabla^2 f|_p \\
\us_{\AA_0}:  \D \times \P T\P^2 \lrab \UL_{\AA_0}, & & \us_{\AA_0}(\ff, \lp):= \ds_{\A_0}(\ff,\p)  \\
\us_{\AA_1}: \D \times \P T\P^2 \lrab \UV_{\AA_1}, & & \us_{\AA_1}(\ff,\lp):= \ds_{\A_1}(\ff, \p) \\
\us_{\DD_4}: \D \times \P T\P^2 \lrab \UV_{\DD_4}, & & \us_{\DD_4}(\ff,\lp):= \ds_{\D_4}(\ff, \p) \\
\us_{\XX_8}: \D \times \P T\P^2 \lrab \UV_{\XX_8}, & & \us_{\XX_8}(\ff,\lp):= \ds_{\XC_8}(\ff, \p).
\end{eqnarray*}
We also have
\begin{eqnarray*}
\us_{\PP \A_2}: \D \times \P T\P^2 \lrab \UV_{\PP \A_2}, & & \{\us_{\PP \A_2}(\ff,\lp)\}(f \otimes p^{\otimes d} \otimes v):= \nabla^2 f|_p (v,\cdot) \\  
\us_{\PP \D_5}^{\UV}: \D \times \P T\P^2 \lrab \UV_{\PP \D_5}, & & \{\us_{\PP \D_5}^{\UV}(\ff,\lp)\}(f \otimes p^{\otimes d} \otimes v^{\otimes 2}):= \nabla^3 f|_p (v,v,\cdot)  \\
\us_{\PP \D_4}: \D \times \P T\P^2 \lrab \UL_{\PP \D_4}, & &  \{\us_{\td_4}(\ff,\lp)\}(f \otimes p^{\otimes d} \otimes w^{\otimes 2}):= f_{02} \\
\us_{\PP \D_5}^{\UL}: \D \times \P T\P^2 \lrab \UL_{\PP \D_5}, & & \{\us_{\PP \D_5}^{\UL}(\ff,\lp)\}(f \otimes p^{\otimes d} \otimes v^{\otimes 2}\otimes w):= f_{21} \\
\us_{\PP \E_6}: \D \times \P T\P^2 \lrab \UL_{\PP \E_6}, & & \{\us_{\PP \E_6}(\ff,\lp)\}(f \otimes p^{\otimes d} \otimes v\otimes w^{\otimes 2}):= f_{12} \\ 
\us_{\PP \E_7}: \D \times \P T\P^2 \lrab \UL_{\PP \E_7}, & & \{\us_{\PP \E_7}(\ff,\lp)\}(f \otimes p^{\otimes d} \otimes v^{\otimes 4}):= f_{40} \\
\us_{\PP \E_8}: \D \times \P T\P^2 \lrab \UL_{\PP \E_8}, & & \{\us_{\PP \E_8}(\ff,\lp)\}(f \otimes p^{\otimes d} \otimes v^{\otimes 3} \otimes w ):= f_{31} \\
\us_{\PP \XC_8}: \D \times \P T\P^2 \lrab \UL_{\PP \XC_8}, & & \{\us_{\PP \XC_8}(f,\lp)\}(f \otimes p^{\otimes d} \otimes w^{\otimes 3}):= f_{03}.
\end{eqnarray*}
We also have sections of the following bundles:  $\us_{\PP \D_5^{\vee}}: \D \times \P T\P^2 \lrab \UL_{\PP \D_5^{\vee}}$ given by
\bge
\label{psi_d5_dual}
\{\us_{\PP \D_5^{\vee}}(\ff,\lp)\}( f^{\otimes 2} \otimes p^{\otimes 2 d} 
\otimes v^{\otimes 2} \otimes w^{\otimes 4}) := 3 f_{12}^2 - 4 f_{21} f_{03}, 
\ede
and  $\us_{\PP \D_6^{\vee}}: \D \times \P T\P^2 \lrab \UL_{\PP \D_6^{\vee}}$ at $(\ff,\lp)$ is given by
\bge
\label{psi_d6_dual}
(f^{\otimes 5} \otimes p^{\otimes 5 d} 
\otimes v^{\otimes 8} \otimes w^{\otimes 4})\mapsto \big(f_{12}^4 f_{40} - 8f_{12}^3 f_{21} f_{31} + 24 f_{12}^2 f_{21}^2 f_{22}-32 f_{12} f_{21}^3 f_{13} + 16 f_{21}^4 f_{04} \big),
\ede
and $\us_{\J}: \D \times \P T\P^2 \lrab \UL_{\J}$ given by
\bge
\label{J_psi}
\{\us_{\J}(\ff,\lp)\}( f^{\otimes 3} \otimes p^{\otimes d} \otimes v^{\otimes 9} \otimes w^{\otimes 3}) := \Big(- \frac{f_{31}^3}{8 } + \frac{3 f_{22} f_{31} f_{40}}{16 } - \frac{f_{13} f_{40}^2}{16}  \Big).
\ede
When $k\geq 3$ we have $\us_{\PP \A_k}: \D \times \P T\P^2 \lrab \UL_{\PP \A_k}$ given by
\bgd
\{\us_{\PP \A_k}(\ff,\lp) \}\big(f^{\otimes (k-2)} \otimes p^{\otimes d} \otimes v^{\otimes k} \otimes w^{\otimes (2k-6)}\big):= f_{02}^{k-3} \A^f_k.
\edd
Similarly, when $k \geq 6$ we have $\us_{\PP \D_k}: \D \times \P T\P^2 \lrab \UL_{\PP \D_k}$ given by
\bgd
\{\us_{\PP \D_k}(\ff,l_p)\}\big(f^{\otimes (1+ \epsilon_k)} \otimes p^{\otimes d(1+\epsilon_k)} \otimes v^{\otimes (k-2+\epsilon_k)}\otimes w^{\otimes (2\epsilon_k)}\big):= f_{12}^{\epsilon_k} \D^f_k,
\edd
where, 
$\epsilon_6 = 0$, $\epsilon_7 =1$ and 
$\epsilon_8=3$. The expressions for $\A^f_k$ (resp. $\D^f_k$) are given below explicitly in 
\eqref{Formula_Ak} (resp. \eqref{Formula_Dk}), till $k=7$ (resp. till $k=8$). \\
\hf\hf Here is an explicit  formula for  
$\A^f_k$ till $k=7$:  
\begin{align}
\label{Formula_Ak}
\A^f_3&= f_{30},\qquad
\A^f_4 = f_{40}-\frac{3 f_{21}^2}{f_{02}}, \qquad
\A^f_5= f_{50} -\frac{10 f_{21} f_{31}}{f_{02}} + 
\frac{15 f_{12} f_{21}^2}{f_{02}^2} \nonumber \\
\A^f_6 &= f_{60}- \f{ 15 f_{21} f_{41}}{f_{02}}-\f{10 f_{31}^2}{f_{02}} + \f{60 f_{12} f_{21} f_{31}}{f_{02}^2}
   +
   \f{45 f_{21}^2 f_{22}}{f_{02}^2} - \f{15 f_{03} f_{21}^3}{f_{02}^3}
   -\f{90 f_{12}^2 f_{21}^2}{f_{02}^3} \nonumber \\  
\A^f_7 &= f_{70} - \frac{21 f_{21} f_{51}}{f_{02}} 
- \frac{35 f_{31} f_{41}}{f_{02}} + \frac{105 f_{12} f_{21} f_{41}}{f_{02}^2} + \f{105 f_{21}^2 f_{32}}{f_{02}^2} + 
\f{70 f_{12} f_{31}^2}{f_{02}^2}+ \f{210 f_{21}f_{22}f_{31}}{f_{02}^2} \nonumber \\
&
-\f{105 f_{03} f_{21}^2 f_{31}}{ f_{02}^3}
-\f{420 f_{12}^2 f_{21} f_{31}}{f_{02}^3}
-\f{630 f_{12}f_{21}^2 f_{22}}{f_{02}^3}
-\f{105 f_{13} f_{21}^3}{f_{02}^3}
+ \f{315 f_{03} f_{12} f_{21}^3}{f_{02}^4}
+ \f{630 f_{12}^3 f_{21}^2}{f_{02}^4}.
\end{align} 
Here is an explicit formula for 
$\D^f_k$  till $k=8$:
\begin{align}
\label{Formula_Dk}
\D^f_6 &=  f_{40},\,\,
\D^f_7 =   f_{50} -\f{5 f_{31}^2}{3 f_{12}},\,\,
\D^f_8 = f_{60} + 
\frac{5 f_{03} f_{31} f_{50}}{3 f_{12}^2} 
-\frac{5 f_{31} f_{41}}{f_{12}} - \frac{10 f_{03} f_{31}^3}{3 f_{12}^3} 
+ \frac{5 f_{22} f_{31}^2}{f_{12}^2}.
\end{align}  

\subsection{The spaces involved}
\hf\hf We begin by explaining a terminology. If $l_{\p} \in \P T_{\p}\P^2$, then we say that $v \in l_{\p}$ if 
$v$ is a tangent vector in $T_{\p}\P^2$ and lies over the fibre of $l_{\p}$.
We now define the spaces that we will encounter.
\label{definition_of_px}
\begin{align*}
\X_k &:= \{ ( \ff,\p) \in \D \times \P^2~~~~~: \textnormal{$f$ has a singularity of type $\X_k$ at $\p$} \} \\
\hat{\X}_k &:= \{ (\ff,\lp) \in \D \times \P T\P^2: \textnormal{$f$ has a singularity of type $\X_k$ at $\p$} \} ~= \pi^{-1}(\X_k) \\
\textnormal{if $~k>1$ } \quad 
\PP \A_k &:= \{ (\ff,\lp) \in \D \times \P T\P^2: 
\textnormal{$f$ has a singularity of type $\A_k$ at $\p$},\\
& \qquad\qquad\qquad\qquad \qquad  \qquad \nabla^2 f|_p(v, \cdot) =0\,\,\textup{if}\,\,v \in  l_{\p}\}\\ 
\PP \D_4 &:= \{ (\ff,\lp) \in \D \times \P T\P^2: 
\textnormal{$f$ has a singularity of 
type $\D_4$ at $\p$}, \\ 
& \qquad \qquad \qquad \qquad \qquad  \qquad \nabla^3 f|_p(v,v,v) =0\,\,\textup{if}\,\,v  \in  l_{\p}\} \\
\textnormal{if $~k>4$} \qquad 
\PP \D_k &:= \{ (\ff,\lp) \in \D \times \P T\P^2: \textnormal{$f$ has a singularity of 
type $\D_k$ at $\p$} \\
& \qquad \qquad \qquad \qquad \qquad\qquad \nabla^3 f|_p(v, v, \cdot) =0\,\,\textup{if}\,\,v \in l_{\p}\} \\
\textnormal{if $k=6, 7$ or  $8$}\qquad 
\PP \E_k &:= \{ (\ff,\lp) \in \D \times \P T\P^2: \textnormal{$f$ has a singularity of 
type $\E_k$ at $\p$}\\
& \qquad \qquad \qquad \qquad \qquad\qquad \nabla^3 f|_p(v, v, \cdot) =0\,\,\textup{if}\,\,v\in l_{\p}\} \\
\textnormal{if $~k>4$} \qquad \PP\D_k^{\vee} &:= \{ (\ff,\lp) \in \D \times \P T\P^2: 
\textnormal{$f$ has a singularity of 
type $\D_k$ at $\p$}, \\
& \qquad \qquad \qquad \qquad \qquad  \qquad \nabla^3 f|_p(v,v,v) =0, ~~\nabla^3 f|_p(v,v,w) \neq 0 \\ 
& \qquad \qquad \qquad \qquad \qquad  \qquad  \,\,\textup{if}\,\,v  \in  l_{\p} ~~\textup{and} ~~w \in (T_{\p}\P^2)/l_{\p}\} 
\end{align*}
\ni We also need the definitions for a few other spaces which will make our computations convenient. 
\begin{align*}
\hat{\A}_1^{\#} := \{ (\ff, l_p) \in \D \times \P T\P^2 &: 
 f(p) =0, \nabla f|_p =0,  \nabla^2 f|_p(v, \cdot) \neq 0, \forall ~v\neq 0 \in l_{\p}  \}  \\
\hat{\D}_4^{\#} := \{ (\ff, \lp) \in \D \times \P T\P^2 &: 
f(p) =0, \nabla f|_p =0, \nabla^2 f|_p \equiv 0, \nabla^3 f|_p (v,v,v) \neq 0, \forall ~v \neq 0 \in l_{\p} \} \\
\hat{\D}_k^{\#\flat} := \{ (\ff, \lp) \in \D \times \P T\P^2 &: \textnormal{$f$ has a $\D_k$ singularity at $\p$},  ~~\nabla^3 f|_p (v,v,v) \neq 0, \forall ~v \neq 0 \in l_{\p},\,\,k\geq 4\} \\ 
\hat{\XC}_8^{\#} := \{ (\ff, \lp) \in \D \times \P T\P^2 &: 
f(p) =0, \nabla f|_p =0, \nabla^2 f|_p \equiv 0, \nabla^3 f|_p =0,  \\ 
                & \qquad \nabla^4 f|_p (v,v,v, v) \neq 0 ~\forall ~v \neq 0  \in l_{\p} \}    \\
\XX_{8}^{\# \flat} :=  \{ (\ff, \lp) \in \D \times \P T\P^2 &: (\ff, \lp) \in \hat{\XC}_8^{\#}, \us_{\J}(\ff, \lp)\neq 0, 
\textnormal{where $\us_{\J}$ is defined in \eqref{J_psi}}\}.                        
\end{align*}

\section{Transversality} 
\label{bundle_sections}

\hf\hf In this section we list down all the relevant bundle sections that are transverse to the zero set. We set up our notation first. Let us define the following projection maps: 
\begin{align*}
\pi_{1}&:\D \times \P^2_1 \times \P^2_2  \lra \D \times \P^2_1, \\
\pi_{2}&:\D \times \P^2_1 \times \P^2_2  \lra \D \times \P^2_2, \\
\pi_{1}&:\D \times \P^2_1 \times \P T \P^2_2  \lra \D \times \P^2_1, \\
\pi_{2}&:\D \times \P^2_1 \times \P T \P^2_2  \lra \D \times \P T \P^2_2. 
\end{align*}
Hence, given a vector bundle over $\D \times \P^2$ or $\D \times \P T\P^2$  we obtain a bundle over 
$\D \times \P^2 \times \P^2$ and $\D \times \P^2 \times \P T\P^2$ respectively, via the pullback maps.  \\
\hf\hf A section of a bundle over $\D \times \P^2$ or $\D \times \P T\P^2$ induces a section over 
the corresponding bundle over  $\D \times \P^2 \times \P^2$ and $\D \times \P^2 \times \P T\P^2$ respectively, via the pullback maps.

\begin{rem}
 To describe bundles over $\D  \times  \P^2$ or $\D \times \P T\P^2$, we follow the abuse of notation 
of omitting pullback maps (as mentioned in Remark \ref{an1_again}). 
However, to describe bundles over $\D \times \P^2 \times \P^2$ or $\D \times\P^2 \times  \P T\P^2$ we do  write the pullback maps. 
\end{rem}
\begin{lmm}
\label{tube_lemma}
Let $\pi:E \lra M$ be a fibre bundle with compact fibers. Let $X\subseteq E$ and $Y\subseteq M$. Then 
\begin{align}
 \pi (\ov{X}) & = \ov{\pi (X)}  \label{tube_lemma_X}\\
\pi^{-1} (\ov{Y}) &= \ov{\pi^{-1}(Y)}. \label{tube_lemma_Y}
\end{align}
\end{lmm}
\pf Since the fibers are compact, $\pi$ is a closed map (Tube Lemma). Combined with the fact that $\pi$ is continuous 
\eqref{tube_lemma_X} follows. Secondly, equation \eqref{tube_lemma_Y} holds for the trivial bundle, hence it also for 
an arbitrary fibre bundle, since this is a local statement. \qed  
\begin{prp}
\label{ift_ml}
The sections of the vector bundles
\begin{align*}
\pi_2^*\ds_{\A_0}:\ov{\A}_1 \times \P^2 - \Delta \ov{\A}_1   \lra \pi_2^* \DL_{\A_0}, \qquad \pi_2^*\ds_{\A_1}: \pi_2^*\ds_{\A_0}^{-1}(0) \lra \pi_2^*\DV_{\A_1}  
\end{align*}
are transverse to the zero set if $d \geq 3$. 
\end{prp}

\begin{prp}
\label{A2_Condition_prp}
The section of the vector bundle $\pi_2^*\us_{\PP \A_2}:  \ov{\A}_1 \circ \ov{\hat{\A}}_1 \lra \pi_2^*\UV_{\PP \A_2}$ is transverse to the zero set, provided $d \geq 4$.  
\end{prp}

\begin{prp}
\label{D4_Condition_prp}
The sections of the vector bundles 
\begin{align*}
\pi_2^*\us_{\PP \A_3}: \ov{\A}_1 \circ \ov{\PP \A}_2 \lra \pi_2^*\UL_{\PP \A_3}, \quad \pi_2^*\us_{\PP \D_4}: \pi_2^*\us_{\PP \A_3}^{-1} (0) \lra \pi_2^*\UL_{\PP \D_4},  
\quad \pi_2^*\us_{\PP \D_5}^{\UL}: \pi_2^*\us_{\PP \D_4}^{-1} (0) \lra \pi_2^*\UL_{\PP \D_5} 
\end{align*} 
are transverse to the zero set provided $d \geq 5$.
\end{prp}

\begin{prp}
\label{A3_Condition_prp}
If $i \geq 4$, then the sections of the vector bundles 
\begin{align*}
\pi_2^*\us_{\PP \A_3}&:\ov{\A}_1 \circ \ov{\PP \A}_2 \lra \pi_2^*\UL_{\PP \A_3}, 
\qquad \pi_2^*\us_{\PP \A_4}:\pi_2^*\us_{\PP \A_{3}}^{-1} (0) - \pi_2^*\us_{\PP \D_4}^{-1}(0) \lra \pi_2^*\UL_{\PP \A_4},  \ldots, \\ 
\pi_2^*\us_{\PP \A_{i}}&:\pi_2^*\us_{\PP \A_{i-1}}^{-1} (0) - \pi_2^*\us_{\PP \D_4}^{-1}(0) \lra \pi_2^*\UL_{\PP \A_{i}} 
\end{align*} 
are transverse to the zero set provided $d\geq i+2$.
\end{prp}

\begin{prp}
\label{D6_Condition_prp}
If $i \geq 6$, then the sections of the vector bundles 
\begin{align*}
\pi_2^*\us_{\PP \D_6}&:\ov{\A}_1 \circ \ov{\PP \D}_5 \lra \pi_2^*\UL_{\PP \D_6}, 
\qquad \pi_2^*\us_{\PP \D_7}:\pi_2^*\us_{\PP \D_{6}}^{-1} (0) - \pi_2^*\us_{\PP \E_6}^{-1}(0) \lra \pi_2^*\UL_{\PP \D_7},  \ldots, \\ 
\pi_2^*\us_{\PP \D_{i}}&:\pi_2^*\us_{\PP \D_{i-1}}^{-1} (0) - \pi_2^*\us_{\PP \E_6}^{-1}(0) \lra \pi_2^*\UL_{\PP \D_{i}} 
\end{align*} 
are transverse to the zero set provided $d\geq i+2$.
\end{prp}

\begin{prp}
\label{E6_Condition_prp}
The section of the vector bundle $\pi_2^*\us_{\PP \E_6}:\ov{\A}_1 \circ \ov{\PP \D}_5 \lra \pi_2^*\UL_{\PP \E_6}$ is transverse to the zero set provided $d\geq 5$.
\end{prp}

\begin{prp}
\label{PD4_Condition_prp}
The section of the vector bundle $\pi_2^*\us_{\PP \A_3}:\ov{\A}_1 \circ \ov{\hat{\D}}_4 \lra \pi_2^*\UL_{\PP \A_3}$ is transverse to the zero set provided $d\geq 5$.
\end{prp}

\pf 
We have omitted the proofs here; they can be found in \cite{BM_Detail}. 
They are similar to the way we prove transversality of bundle sections in \cite{BM13}.
\qed  

\section{Closure and Euler class contribution} 
\label{closure_of_spaces}

\hf\hf In this section we compute closure of a singularity with one $\A_1$-node. Along the way 
we also compute how much a certain section contributes to the Euler class of a bundle. But 
first, let us recapitulate what we know about the closure of one singular point 
which was proved in \cite{BM13}.
\begin{lmm}
\label{cl}
Let $\X_k$ be a singularity of type $\A_k$, $\D_k$, 
$\E_k$ or $\XC_8$. Then the closures are given by :
\begin{enumerate} 
\item \label{A0cl} $\ov{\A}_0 = \A_0 \cup \ov{\A}_1$ 
\qquad if  $d \geq 3$.
\item \label{A1cl}
$\ov{\hat{\A}_1} = \ov{\hat{\A}^{\#}_1} = \hat{\A}_1^{\#} \cup 
\ov{\PP \A}_2$ \qquad if $d \geq 3$.
\item \label{D4_cl_no_direction}$ \ov{\hat \D^{\#}_4} = \hat \D^{\#}_4 \cup \ov{\PP \D}_4 $ \qquad if $d \geq 3$.
\item \label{D4cl}$ \ov{\PP \D}_4 = \PP \D_4 \cup \ov{\PP \D}_5 \cup \ov{\PP \D_5^{\vee}}$ 
\qquad if $d \geq 4$. 
\item \label{E6cl}$\ov{\PP \E}_6 = 
\PP \E_6 \cup \ov{\PP \E}_7 \cup \ov{\hat{\XC}^{\#}_8}$ 
\qquad if $d \geq 4$. 
\item \label{D5cl}$\ov{\PP \D}_5 = \PP \D_5 \cup \ov{\PP \D}_6 \cup \ov{\PP \E}_6$ \qquad if $d \geq 4$.
\item \label{D6cl}$\ov{\PP \D}_6 = \PP \D_6 \cup \ov{\PP \D}_7 \cup \ov{\PP \E}_7$ \qquad if $d \geq 5$.
\item \label{A2cl}$\ov{\PP \A}_2 = \PP \A_2 \cup 
\ov{\PP \A}_3 \cup \ov{\hat{\D}_4^{\#}} $ 
\qquad if $d \geq 4$.
\item \label{A3cl}$\ov{\PP \A}_3 = \PP \A_3 \cup  \ov{\PP \A}_4 \cup \ov{\PP \D}_4$ \qquad if $d \geq 5$.
\item \label{A4cl}$\ov{\PP \A}_4 = \PP \A_4 \cup  \ov{\PP \A}_5 \cup \ov{\PP \D}_5$ \qquad if $d \geq 6$.
\item \label{A5cl}$\ov{\PP \A}_5 = \PP \A_5 \cup  \ov{\PP \A}_6 \cup \ov{\PP \D}_6 \cup \ov{\PP \E}_6 $ \qquad if $d \geq 7$.
\item \label{A6cl}$\ov{\PP \A}_6 = \PP \A_6 \cup  \ov{\PP \A}_7 \cup 
\ov{\PP \D}_7 \cup \ov{\PP \E}_7 \cup \ov{\hat{\XC}_8^{\# \flat}}$  
 \qquad if $d \geq 8$.
\end{enumerate}
\end{lmm}

\ni Let us now state a few facts about the closure of one singular point that will be required in this paper. 
These facts were not explicitly stated in \cite{BM13} because it was not needed in that paper.  

\begin{lmm}
\label{Dk_sharp_closure}
We have the following equality (or inclusion) of sets 
\begin{enumerate}
\item \label{a1_cl_new} $\ov{\A}_1 = \A_1 \cup \ov{\A}_2$ \qquad if ~$d \geq 2$. \vspace*{-0.1cm} 
 \item \label{d4_new} $\ov{\hat{\D}^{\#}_4} = \ov{\hat{\D}}_4$ \qquad if ~$d \geq 3$. \vspace*{-0.1cm}
 \item \label{dk_new} $\ov{\hat{\D}^{\#\flat}_k} = \ov{\hat{\D}}_k$ \qquad if ~$k\geq 4$ ~~ and ~~$d \geq 3$. 
 \vspace*{-0.1cm}
 \item \label{d5_pa3_zero} $\big\{ (\ff, \lp) \in \ov{\hat{\D}}_5: \us_{\PP \A_3}(\ff, \lp) =0 \big\} 
 = \ov{\PP \D}_5 \cup \ov{\PP \D^{\vee}_5}$  \qquad if ~$d \geq 3$. \vspace*{-0.1cm} 
 \item \label{pd6_pd7_pd8_closure} 
 $\big\{(\ff, \lp) \in  \ov{\PP \D}_6: \us_{\PP \E_6}(\ff, \lp) \neq 0 \big\} 
 =\PP \D_6 \cup \PP \D_7\cup \big\{(\ff, \lp) \in \ov{\PP \D}_8: \us_{\PP \E_6}(\ff, \lp) \neq 0 \big\}$
 \item \label{pe6_subset_of_cl_pd5_dual} $\ov{\PP \E}_6  \subset \ov{\PP \D_5^{\vee}}$ \qquad if ~$d \geq 3$. 
\end{enumerate}
\end{lmm}
The proofs are straightforward; the details can be found in \cite{BM_Detail}. 
Hence we stated Lemma \ref{Dk_sharp_closure}, statements \ref{d4_new} and \ref{dk_new} separately to avoid confusion.

\hf \hf Before stating the main results of this section, 
let us define three more spaces which will be required 
while formulating some of the Lemmas: 
\begin{align}
\Delta \mp &:= \{ (\ff, \p, \lp) \in \ov{\ov{\A}_1 \circ \PP \A}_4: 
\pi_2^*\us_{\PP \D_4}(\ff, \p, \lp) =0, 
~\pi_2^*\us_{\PP \E_6}(\ff, \p, \lp) \neq 0  \}, \nonumber \\
\Delta \mq &:= \{ (\ff, \p, \lp) \in \ov{\ov{\A}_1 \circ \PP \A}_5: 
\pi_2^*\us_{\PP \D_4}(\ff, \p, \lp) =0, \nonumber 
~\pi_2^*\us_{\PP \E_6}(\ff, \p, \lp) \neq 0  \}, \\
\Delta \mr & := \{ (\ff, \p, \lp) \in \ov{\ov{\A}_1 \circ \PP \D}_4: 
\pi_2^*\us_{\PP \D_5}(\ff, \p, \lp) \neq 0  \}.\label{new_defn_delta}
\end{align} 
We are now ready to state the main Lemmas. 

\begin{lmm}
\label{cl_two_pt}
Let $\X_k$ be a singularity of type $\A_k$, $\D_k$, 
$\E_k$. Then their closures with one $\A_1$-node are given by :
\begin{enumerate}
\item \label{a1a_minus_1_cl} $\ov{\ov{\A}_1 \circ (\D \times \P^2)} = \ov{\A}_1 \circ (\D \times \P^2) \du \Delta \ov{\A}_1 $ 
\qquad if  $d \geq 1$. 
\item \label{a1a1_up_cl} $\ov{\ov{\A}_1 \circ  \hat{\A}^{\#}_1} = \ov{\A}_1\circ \hat{\A}^{\#}_1 \du  \ov{\A}_1 \circ  (\ov{\hat{\A}_1^{\#}}- \hat{\A}_1^{\#}  ) \du \Delta \ov{\hat{\A}}_3$ 
\qquad if  $d \geq 3$. 
\item \label{a1_pa2_cl} $\ov{\ov{\A}_1 \circ  \PP \A}_2 = \ov{\A}_1 \circ  \PP \A_2 \du \ov{\A}_1 \circ (\ov{\PP \A}_2- \PP \A_2) \du \Big(\Delta \ov{\PP \A}_4 \cup 
\Delta \ov{\hat{\D}^{\#\flat}_5} \Big )$  \qquad if  $d \geq 4$.
\item \label{a1_pa3_cl} $\ov{\ov{\A}_1 \circ  \PP \A}_3 = \ov{\A}_1 \circ  \PP \A_3 \du \ov{\A}_1 \circ (\ov{\PP \A}_3- \PP \A_3) \du 
\Big( \Delta \ov{\PP \A}_5 \cup \Delta \ov{\PP \D^{\vee}_5} \Big)$ \qquad if  $d \geq 5$.  
\item \label{a1_pa4_cl} $\ov{\ov{\A}_1 \circ  \PP \A}_4 = \ov{\A}_1 \circ  \PP \A_4 \du \ov{\A}_1 \circ (\ov{\PP \A}_4- \PP \A_4) \du 
\Big( \Delta \ov{\PP \A}_6 \cup  
\Delta \ov{\mp} \cup \Delta \ov{\PP \E}_6  \Big) $ \qquad if  $d \geq 6$. 
\item \label{a1_pa5_cl} $\ov{\ov{\A}_1 \circ  \PP \A}_5 = \ov{\A}_1 \circ  \PP \A_5 \du \ov{\A}_1 \circ (\ov{\PP \A}_5- \PP \A_5) \du 
\Big( \Delta \ov{\PP \A}_7 \cup \Delta  \ov{\mq} \cup \Delta \ov{\PP \E}_7 \Big)$ \qquad if  $d \geq 7$. 
\item \label{a1_pd4_cl} $\ov{\ov{\A}_1 \circ  \PP \D}_4 = \ov{\A}_1 \circ  \PP \D_4 \du \ov{\A}_1 \circ (\ov{\PP \D}_4- \PP \D_4) \du 
\Big( \Delta \ov{\mr} \cup \Delta \ov{\PP \D}_6 \Big) $  \qquad if  $d \geq 4$.
\item \label{a1_d4_cl} $\ov{\ov{\A}_1 \circ  \hat{\D}_4} = \ov{\A}_1 \circ  \hat{\D}_4 \du \ov{\A}_1 \circ (\ov{\hat{\D}}_4- \hat{\D}_4) \du   
\Big( \Delta \ov{\hat{\D}}_6 \Big) $  \qquad if  $d \geq 4$.
\item \label{a1_pd5_cl} $\ov{\ov{\A}_1 \circ  \PP \D}_5 = \ov{\A}_1 \circ  \PP \D_5 \du \ov{\A}_1 \circ (\ov{\PP \D}_5- \PP \D_5) \du 
\Big( \Delta \ov{\PP \D}_7 \cup \Delta \ov{\PP \E}_7 \Big)$ \qquad if  $d \geq 5$.  
\end{enumerate}
\end{lmm}

\begin{rem} Although the statement of 
Lemma \ref{cl_two_pt} (\ref{a1a_minus_1_cl}) is trivial we give an elaborate proof for two reasons. 
Firstly, along the way, we prove a few other statements that will be required later. 
Secondly, as a corollary we also compute the contribution of certain sections to the Euler class of relevant bundles. 
\end{rem}


\textbf{Proof of Lemma \ref{cl_two_pt} (\ref{a1a_minus_1_cl}):} It suffices to show  that 
\begin{align}
\{ (\ff, \p, \p) \in \ov{\ov{\A}_1 \circ (\D \times \P^2)}\} &= \Delta \ov{\A}_{1}. \label{a1_a0_is_a1}  
\end{align} 
Clearly the lhs\footnote{We shall use {\it lhs} to denote {\it left hand side} and {\it rhs} to denote {\it right hand side} of an equation.} of \eqref{a1_a0_is_a1} is a subset of its rhs. To show the converse we will prove the 
following two claims simultaneously: 
\begin{align}
 \ov{\ov{\A}_1 \circ (\D \times \P^2)} & \supset \Delta (\A_1 \du \A_2),  \label{a1_du_a2_is_subset_of_a1_a0}\\
 \ov{\ov{\A}_1 \circ  \A}_1  \cap \Delta (\A_1 \du \A_2)  & = \varnothing. \label{a1_du_a2_intersect_a1_a1_is_empty} 
\end{align}
Since $\ov{\ov{\A}_1 \circ (\D \times \P^2)}$ 
is a closed set, \eqref{a1_du_a2_is_subset_of_a1_a0} 
implies that the rhs of \eqref{a1_a0_is_a1} is a subset of its lhs.\footnote{In fact
the full strength of \eqref{a1_du_a2_is_subset_of_a1_a0} is not really needed; 
we simply need that $\ov{\ov{\A}_1 \circ (\D \times \P^2)} \supset \Delta \A_1$.}
Moreover, \eqref{a1_du_a2_intersect_a1_a1_is_empty} is not required at all 
for the proof of Lemma \ref{cl_two_pt} \eqref{a1a_minus_1_cl}. 
However, these statements will be required in the proofs of  
Lemma \ref{cl_two_pt} (\ref{a1a1_up_cl}).
Furthermore, in the process of proving \eqref{a1_du_a2_is_subset_of_a1_a0} and \eqref{a1_du_a2_intersect_a1_a1_is_empty}, 
we will also be computing the contribution of certain sections to the Euler class of relevant bundles as a corollary.

\begin{claim}
\label{a1_a1_closure_intersect_a1_or_a2_empty_equations_claim}
Let $(\ff,\p, \p) \in \Delta (\A_1 \cup \A_{2})$.
Then there exist solutions 
$$(\ff(t_1, t_2), \p(t_1, t_2), \p(t_1)  ) \in \D \times \P^2 \times \P^2 $$ 
sufficiently close to $(\ff,\p, \p)$ to the set of equations
\begin{align}
\pi_1^*\ds_{\A_0}( \ff(t_1, t_2), \p(t_1, t_2), \p(t_1)) & = 0, ~~\pi_1^*\ds_{\A_1}( \ff(t_1, t_2), \p(t_1, t_2), \p(t_1) ) = 0, ~~\p(t_1, t_2) \neq \p(t_1).  \label{a1_a1_sharp_closure_functional_eqn_a1_plus_a2} 
\end{align}
Furthermore, any such solution sufficiently close to $(\ff, \p, \p)$ satisfies 
\begin{align}
\Big( \pi_2^*\ds_{\A_0}(\ff(t_1, t_2), \p(t_1, t_2), \p(t_1)), ~\pi_2^*\ds_{\A_1}(\ff(t_1, t_2), \p(t_1, t_2), \p(t_1)) \Big) \neq (0,0). \label{a1_a1_closure_intersect_a1_is_empty_equation}
\end{align} 
In particular, $(\ff(t_1, t_2), \p(t_1, t_2), \p(t_1))$ does not lie in $\ov{\ov{\A}_1 \circ \A}_1$.
\end{claim}
It is easy to see that claim \ref{a1_a1_closure_intersect_a1_or_a2_empty_equations_claim} proves statements \eqref{a1_du_a2_is_subset_of_a1_a0} and 
\eqref{a1_du_a2_intersect_a1_a1_is_empty} 
simultaneously. \\   

\pf  Choose homogeneous coordinates $[\mathrm{X}: \mathrm{Y}: \mathrm{Z}]$ 
so that $\p = [0:0:1]$ and  let $\mathcal{U}_{\p}$ be a sufficiently 
small neighbourhood of $\p$ inside $\P^2$.
Denote   
$\pi_{x}, \pi_y : \mathcal{U}_{\p} \lra \C $ 
 to be the projection maps given by 
$$\pi_{x}([\mathrm{X}: \mathrm{Y}:\mathrm{Z}]) := \mathrm{X}/\mathrm{Z} \qquad \textnormal{and} \qquad \pi_{y}([\mathrm{X}:\mathrm{Y}:\mathrm{Z}]) := \mathrm{Y}/\mathrm{Z},$$ 
and $v, w: \mathcal{U}_{\p} \lra T \P^2$ the vector fields dual to the one forms $d\pi_x$ and $d\pi_y $ respectively. Let $ (\ff (t_1, t_2),  \p(t_1)) \in \D \times \P^2 $ be an arbitrary point that is close to 
$( \ff, \p)$ and let $\p(t_1, t_2)$ be a point in $\P^2$ that is close to $\p(t_1)$. Let 
\begin{align*}
\p(t_1) := [x_{t_1}:y_{t_1}:1] \in \P^2,  \qquad p(t_1) &:= (x_{t_1},y_{t_1},1) \in \C^3,\\ 
\p(t_1, t_2) := [x_{t_1} + x_{t_2}: ~y_{t_1}+y_{t_2}: 1] \in \P^2, \qquad p(t_1, t_2) &:= (x_{t_1} + x_{t_2}, ~y_{t_1}+y_{t_2}, 1) \in \C^3, \\
\ff(t_1, t_2)  \in \P^{\delta_d}, \qquad f(t_1, t_2) & \in \C^{\delta_d +1}.  
\end{align*}
Define the following numbers: 
\begin{align*}
f_{ij}(t_1, t_2) & := \{ \nabla^{i+j} f(t_1, t_2)|_{\p(t_1)} 
(\underbrace{v,\cdots v}_{\textnormal{$i$ times}}, \underbrace{w,\cdots w}_{\textnormal{$j$ times}}) \}(p(t_1)^{\otimes d}),\\
\FF   &:= f_{00}(t_1, t_2) + f_{10}(t_1, t_2) \xt + f_{01}(t_1, t_2) \yt +  \sum_{i+j=2}\big({\textstyle \frac{f_{ij}(t_1, t_2)}{i ! j !}} x_{t_2}^i y_{t_2}^j\big) + \ldots,  \\  
\Fx &:= f_{10}(t_1, t_2)  +  f_{20}(t_1, t_2) \xt + f_{11}(t_1, t_2) y_{t_2} + \ldots, \\ 
\Fy &:= f_{01}(t_1, t_2) + f_{11}(t_1, t_2) x_{t_2} + f_{02}(t_1, t_2) y_{t_2} + \ldots.   
\end{align*}
It is easy to see that  
\begin{align*}
\{ \pi_1^* \ds_{\A_0}(\ff(t_1, t_2), \p(t_1, t_2),  \p(t_1)) \} (f(t_1, t_2) \otimes p(t_1, t_2)^{\otimes d} ) & = \FF, \\
\{ \pi_1^* \ds_{\A_1}(\ff(t_1, t_2), \p(t_1, t_2),  \p(t_1)) \} (f(t_1, t_2) \otimes p(t_1, t_2)^{\otimes d} \otimes v) & = \Fx, \\ 
\{ \pi_1^* \ds_{\A_1}(\ff(t_1, t_2), \p(t_1, t_2),  \p(t_1)) \} (f(t_1, t_2) \otimes p(t_1, t_2)^{\otimes d} \otimes w) &=  \Fy.
\end{align*}
We now observe  that 
\eqref{a1_a1_sharp_closure_functional_eqn_a1_plus_a2}  
has a solution if and only if 
the following set of equations has a solution 
\begin{align}
\mathrm{F} = 0, \qquad \Fx = 0, \qquad \Fy = 0,  \qquad (x_{t_2}, y_{t_2}) \neq (0, 0) \qquad \textnormal{(but small)}. \label{eval_f1} 
\end{align} 
Note that in equation \eqref{a1_a1_sharp_closure_functional_eqn_a1_plus_a2} the equality holds as 
\textit{functionals}, while in equation \eqref{eval_f1}, the equality holds as \textit{numbers}. 
We will now show that  
\eqref{eval_f1} (and hence \eqref{a1_a1_sharp_closure_functional_eqn_a1_plus_a2}) 
has solutions whenever $(\ff, \p, \p) \in \Delta (\A_1 \cup \A_2)$. 
Furthermore, for all those solutions,  
\eqref{a1_a1_closure_intersect_a1_is_empty_equation} holds. \\
\hf \hf First let us assume 
$(\ff, \p, \p) \in \Delta \A_1$. 
It is obvious that solutions to \eqref{eval_f1} exist; we can solve for $f_{10}(t_1, t_2)$ and $f_{01}(t_1, t_2)$ 
using $\Fx =0$ and $\Fy =0$ and 
and then solve for 
$f_{00}(t_1, t_2)$ using $\FF =0$. 
To show that \eqref{a1_a1_closure_intersect_a1_is_empty_equation} holds it suffices to show that 
if $(\xt, \yt)$ is small but non zero,  then 
\begin{align}
\big(f_{00}(t_1, t_2), ~f_{10}(t_1, t_2), ~f_{01}(t_1, t_2) \big) &\neq (0,0,0). \label{a1_a1_closure_intersect_a1_is_empty_equation_numbers}
\end{align}
Observe that  
\eqref{eval_f1}  implies that 
\begin{equation}
\label{f_ij_matrix_equation}
\left( \begin{array}{c} f_{10}(t_1, t_2) \\ f_{01}(t_1, t_2) \end{array} \right) = 
-\begin{pmatrix} f_{20}(t_1,t_2) & f_{11}(t_1, t_2)  \\ f_{11}(t_1, t_2)  & f_{02}(t_1, t_2) 
\end{pmatrix} \left(\begin{array}{c} x_{t_2} \\ y_{t_2} \end{array} \right)+ \left (\begin{array}{c} \mathrm{E}_1(\x_{t_2}, \y_{t_2}) \\ \mathrm{E}_2(\x_{t_2}, \y_{t_2})   
\end{array} \right) 
\end{equation}
where $\mathrm{E}_i(\x_{t_2}, \y_{t_2})$ are second order in $(\x_{t_2}, \y_{t_2})$. 
Since $(\ff,\p, \p) \in \Delta \A_1$, the matrix
\begin{equation*}
\mathrm{M} := \left( \begin{array}{cc}
f_{20}(t_1,t_2) & f_{11}(t_1, t_2)    \\
f_{11}(t_1, t_2) & f_{02}(t_1, t_2)  \end{array} \right)
\end{equation*}
is invertible if $\ff(t_1, t_2)$ is sufficiently close to $\ff$. 
Equation \eqref{f_ij_matrix_equation} now implies that  
that if $(\xt, \yt)$ is small but non zero, then  
$f_{10}(t_1, t_2)$ and $f_{01}(t_1, t_2)$ can not both be zero. 
Hence \eqref{a1_a1_closure_intersect_a1_is_empty_equation_numbers}  holds and hence \eqref{a1_a1_closure_intersect_a1_is_empty_equation} holds. \\
\hf \hf Next, let $(\ff, \p, \p) \in \Delta \A_2$. 
Since $\Delta \A_2 \subset \Delta \ov{\A}_1$, we conclude that solutions to 
\eqref{a1_a1_sharp_closure_functional_eqn_a1_plus_a2} exist; we only need to show that 
\eqref{a1_a1_closure_intersect_a1_is_empty_equation} holds. 
Observe that 
$f_{20}$ and $f_{02}$ can not both be zero; assume $f_{02} \neq 0$ 
(and hence $f_{02}(t_1, t_2) \neq 0$). 
Write $\mathrm{F}$ as 
\bgd
\mathrm{F} = f_{00}(t_1, t_2) +f_{10}(t_1, t_2) \xt + f_{01}(t_1, t_2) \yt +  \Z_0(x_{t_2}) + \Z_1(x_{t_2})y_{t_2} + \Z_2(x_{t_2}) y_{t_2}^2 + \ldots
\edd
where $\Z_2(0) \neq 0.$ We claim that there exists a unique holomorphic function $\Y(x_{t_2})$, 
vanishing at the origin,  such that 
after we make a change of coordinates $y_{t_2} = \hat{y}_{t_2} + \Y(x_{t_2})$, 
the function $\mathrm{F}$ is given by 
\bgd
\mathrm{F} = f_{00}(t_1, t_2) +f_{10}(t_1, t_2) \xt + f_{01}(t_1, t_2) \hat{y}_{t_2} +  f_{01}(t_1, t_2)\Y(x_{t_2}) +  
 \hat{\Z}_0(x_{t_2}) + \hat{\Z}_2(x_{t_2}) \hat{y}_{t_2}^2 + \hat{\Z}_3(x_{t_2}) \hat{y}_{t_2}^3 + \ldots 
\edd
for some $\hat \Z_k(x_{t_2})$ (i.e., $\hat{\Z}_1(x_{t_2}) \equiv 0$). 
This is possible if $\Y(x_{t_2})$ satisfies the identity
\begin{align}
\Z_1(x_{t_2}) + 2 \Z_2(x_{t_2}) \Y(\xt) + 3 \Z_3(x_{t_2}) \Y(\xt)^2 + \ldots \equiv 0.  \label{psconvgg2}
\end{align}
Since $\Z_2(0) \neq 0$, $\Y(x_{t_2})$ exists by the Implicit Function Theorem and 
we can compute $\Y(x_{t_2})$ explicitly  
as a power series using \eqref{psconvgg2} and then
compute $\hat \Z_{0}(x_{t_2})$. Hence, 
\begin{align*} 
\mathrm{F} & = f_{00}(t_1, t_2) +f_{10}(t_1, t_2) x_{t_2} + f_{01}(t_1, t_2) \hat{y}_{t_2} +  \varphi(x_{t_2}, \hat{y}_{t_2})  \hat{y}_{t_2}^2  
+f_{01}(t_1, t_2)\Y(x_{t_2}) \\ 
    & + \underbrace{{\textstyle \frac{\B^{f(t_1, t_2)}_2}{2!}x_{t_2}^2 + \frac{\B^{f(t_1, t_2)}_3}{3!} x_{t_2}^3 + \mathrm{O}(x_{t_2}^4)}}_{\hat \Z_{0}(x_{t_2})},
 \end{align*}
 where $\varphi(0,0) \neq 0$ and
\begin{align*}
 \B^{f(t_1, t_2)}_2  & := f_{20}(t_1, t_2) - \frac{f_{11}(t_1, t_2)^2}{f_{02}(t_1, t_2)}, 
\qquad \textup{and}\\
\B^{f(t_1, t_2)}_{3} & :=f_{30}(t_1, t_2)-\frac{3 f_{11}(t_1, t_2) f_{21}(t_1, t_2)}
{f_{02}(t_1, t_2)} + \frac{3 f_{11}(t_1, t_2)^2 f_{12}(t_1, t_2)}{f_{02}(t_1, t_2)^2}- \frac{f_{11}(t_1, t_2)^3 
f_{03}(t_1, t_2)}{f_{02}(t_1, t_2)^3}\neq 0.  
\end{align*}
The last inequality holds because $(\ff, \p, \p) \in \Delta \A_2$.
In these new coordinates $\hat{y}_{t_2}$ 
and $x_{t_2}$, 
equation \eqref{eval_f1} is equivalent to  
\begin{align}
\mathrm{F} = f_{00}(t_1, t_2) + f_{10}(t_1, t_2) x_{t_2} + f_{01}(t_1, t_2) \hat{y}_{t_2} + \varphi(x_{t_2}, \hat{y}_{t_2})  \hat{y}_{t_2}^2 + f_{01}(t_1, t_2)\Y(x_{t_2}) + &   \nonumber \\
    +\frac{\B^{f(t_1, t_2)}_2}{2!}x_{t_2}^2 + \frac{\B^{f(t_1, t_2)}_3}{3!} x_{t_2}^3  + \mathrm{O}(x_{t_2}^4) & = 0,  \label{closure_a1_a1_is_not_a2_F} \\
f_{10}(t_1, t_2)+ f_{01} (t_1, t_2) \mathrm{B}^{\prime}(x_{t_2})+ \varphi_{x_{t_2}}(x_{t_2}, \hat{y}_{t_2})  \hat{y}_{t_2}^2  + \B^{f(t_1, t_2)}_2 x_{t_2} + \frac{\B^{f(t_1, t_2)}_3}{2!} x_{t_2}^2   + 
\mathrm{O}(x_{t_2}^3) & = 0, 
\label{closure_a1_a1_is_not_a2_Fx}\\
f_{01}(t_1, t_2)+ 2 \hat{y}_{t_2}  \varphi(x_{t_2}, \hat{y}_{t_2}) +  \hat{y}_{t_2}^2 \varphi_{\hat{y}_{t_2}}(x_{t_2}, \hat{y}_{t_2}) & =0, \nonumber \\ 
(\xt, \hat{y}_{t_2})  \neq (0,0) \qquad \textnormal{but small}.& \label{closure_a1_a1_is_not_a2_Fy}
\end{align}
Let us clarify a point of confusion: we are claiming that \eqref{eval_f1} has a solution if and only if the equation 
\begin{align}
\mathrm{F}& = 0, \qquad \Fx = 0, \qquad \FF_{\hat{y}_{t_2}} = 0,  \qquad (x_{t_2}, \hat{y}_{t_2}) \neq (0, 0) \qquad \textnormal{(but small)} \label{eval_f1_hat} 
\end{align}
has a solution. We are \textit{not} claiming that the partial derivatives in the old 
coordinates are individually equal to the partial derivatives in the new coordinates.  
Since $\varphi(0,0) \neq 0$ we can use \eqref{closure_a1_a1_is_not_a2_Fy} to solve for $\hat{y}_{t_2}$ in terms of $x$ and 
$f_{01}(t_1, t_2)$ to get 
\begin{align}
 \hat{y}_{t_2} & = f_{01} (t_1, t_2) \mathrm{E}(x_{t_2}, f_{01}(t_1, t_2)), \label{closure_a1_a1_is_not_a2_F_hat_y}
\end{align}
where $\mathrm{E}(x_{t_2}, f_{01}(t_1, t_2))$ is a  holomorphic 
function of $(x_{t_2}, f_{01}(t_1, t_2))$.
Using \eqref{closure_a1_a1_is_not_a2_F_hat_y}, \eqref{closure_a1_a1_is_not_a2_Fx} and \eqref{closure_a1_a1_is_not_a2_F}  
we get (by eliminating $\B^{f(t_1, t_2)}_2$ and $\hat{y}_{t_2}$)    
\begin{align}
\mathrm{F} = -\frac{\B^{f(t_1, t_2)}_3}{12} x_{t_2}^3 + \mathrm{O}(x_{t_2}^4) + f_{00}(t_1, t_2)+f_{10}(t_1, t_2) \mathrm{E}_{1} (x_{t_2}, f_{10}(t_1, t_2), f_{01}(t_1, t_2))  &  \nonumber \\ 
+f_{01}(t_1, t_2) \mathrm{E}_{2} (x_{t_2}, f_{10}(t_1, t_2), f_{01}(t_1, t_2))  &  = 0 \label{closure_a1_a1_is_not_a2_F_eliminated}
\end{align}
where $\mathrm{E}_i (x_{t_2}, f_{10}(t_1, t_2), f_{01}(t_1, t_2))$ is a  holomorphic 
function of $(x_{t_2}, f_{10}(t_1, t_2), f_{01}(t_1, t_2))$. 
Since solutions to equation \eqref{eval_f1}
satisfy  \eqref{closure_a1_a1_is_not_a2_F_eliminated}, we conclude 
that $f_{00}(t_1, t_2)$,
$f_{10}(t_1, t_2)$ and $f_{01}(t_1, t_2)$ can not all be zero.  
If they were all zero 
then $\mathrm{F}$ could not be zero for small but non zero $x_{t_2}$ (by \eqref{closure_a1_a1_is_not_a2_F_eliminated}), 
since 
$\B^{f(t_1, t_2)}_3 \neq 0$ (this is where we are using $(\ff, \p, \p) \in \Delta \A_2$). 
This implies \eqref{a1_a1_closure_intersect_a1_is_empty_equation} holds 
as functionals. This completes the proof of Lemma \ref{cl_two_pt} \eqref{a1a_minus_1_cl}. \qed  \\


\hf \hf Before proceeding to the proof of Lemma \ref{cl_two_pt} (\ref{a1a1_up_cl}),  
let us prove  a corollary which will be needed in the proof of equation \eqref{algoa1a1} in section \ref{Euler_class_computation}. The proof of this corollary follows from the setup of the preceding proof, hence we prove it here.
\begin{cor}
\label{a1_section_contrib_from_a1_and_a2}
Let $\WL \lra \D \times \P^2 \times \P^2$ 
be a vector bundle  
such that the rank of $\WL$  is 
same as the dimension of $\Delta \A_2$\footnote{This dimension is also one less than the dimension of $ \Delta \A_{1}$, equal to $\delta_d-2$.} 
and $\Q: \D \times \P^2 \times \P^2 \lra \WL$ 
a generic smooth section.  
Then the contribution of the section 
\[ \pi_2^* \ds_{\A_{0}} \oplus \pi_2^* \ds_{\A_{1}}  \oplus \Q: \ov{\A}_1 \times \P^2 \lra \pi_2^* \DL_{\A_0} \oplus  \pi_2^* \DV_{\A_1} \oplus \WL \]
to the Euler class from the points of $\Delta \A_{1}$ is given by 
\begin{align}
\mathcal{C}_{\Delta \A_1} (\pi_2^* \ds_{\A_{0}} \oplus \pi_2^* \ds_{\A_{1}} \oplus \Q) = \big\langle e( \pi_2^* \DL_{\A_0} \oplus \WL), ~[\Delta \ov{\A}_1] \big\rangle.  \label{contrib_from_delta_a1_psi_a1}
\end{align}
Secondly, 
if  $(\ff,\p, \p) \in \Delta \A_{2} \cap \Q^{-1}(0)$, 
then this section vanishes on $(\ff, \p,  \p)$ with a multiplicity of $3$. 
\end{cor}

\begin{rem}
When we use the phrase ``number of zeros''  (resp. ``number of solutions'') our intended meaning is number of zeros counted with a sign (resp. the number of solutions counted with a sign).
\end{rem}

\textbf{Proof of Corollary  \ref{a1_section_contrib_from_a1_and_a2}:} The contribution of $\pi_2^* \ds_{\A_{0}}\oplus \pi_2^* \ds_{\A_{1}} \oplus \Q$ 
to the the Euler class 
from the points of $\Delta\A_1$
is 
the number of solutions of 
\begin{align}
\pi_2^*\ds_{\A_0} (\ff(t_1, t_2), \p(t_1, t_2), \p(t_1)) &=  \nu_0(\ff(t_1, t_2), \p(t_1, t_2), \p(t_1)), \nonumber \\  
\pi_2^*\ds_{\A_1} (\ff(t_1, t_2), \p(t_1, t_2), \p(t_1)) &=  \nu_1(\ff(t_1, t_2), \p(t_1, t_2), \p(t_1)), \nonumber  \\ 
\Q(\ff(t_1, t_2), \p(t_1, t_2), \p(t_1)) &= 0, ~~(\ff(t_1, t_2), \p(t_1, t_2), \p(t_1)) \in \mathcal{U}_{\mathcal{K}} \subset \ov{\A}_1 \times \P^2, \label{psi_A1_contribution_from_A1}
\end{align}
where $\mathcal{K}$ is a sufficiently 
large compact subset of $\Delta \A_1$, $\mathcal{U}_{\mathcal{K}}$ is a sufficiently small neighborhood of $\mathcal{K}$ inside $\ov{\A}_1 \times \P^2$, and $\nu_0$ and $\nu_1$ are generic smooth perturbations. We do not need to perturb $\Q$ since it is already generic. We will convert the functional equation \eqref{psi_A1_contribution_from_A1} into an equation that involves equality of numbers. 
Let $\theta_1, \theta_2, \ldots, \theta_{\delta_d-2}$ form a basis for $\WL^*$ at the 
point $(\ff(t_1, t_2), \p(t_1, t_2), \p(t_1))$. 
Define 
\begin{align}
\xi_{0} (\ff(t_1, t_2), \p(t_1), (\xt, \yt)) &:=  \{ \nu_0(\ff(t_1, t_2), \p(t_1, t_2), \p(t_1)) \}(f(t_1, t_2) \otimes p(t_1)^{\otimes d}) \in \C, \nonumber \\
\xi_{1x} (\ff(t_1, t_2), \p(t_1), (\xt, \yt)) &:=  \{ \nu_1(\ff(t_1, t_2), \p(t_1, t_2), \p(t_1)) \}(f(t_1, t_2) \otimes p(t_1)^{\otimes d} \otimes v) \in \C,\nonumber  \\
\xi_{1y} (\ff(t_1, t_2), \p(t_1), (\xt, \yt)) &:=  \{ \nu_1(\ff(t_1, t_2), \p(t_1, t_2), \p(t_1)) \}(f(t_1, t_2) \otimes p(t_1)^{\otimes d} \otimes w)\in \C, \nonumber \\
\R(\ff(t_1, t_2), \p(t_1), (\xt, \yt)) &:= \{ \Q(\ff(t_1, t_2), \p(t_1, t_2), \p(t_1)) \}(\theta_1 \oplus \theta_2 \oplus \ldots \oplus \theta_{\delta_{d-2}}) \in \C^{\delta_d-2}. \label{number_defn}
\end{align}
Consider now the following set of equations 
\begin{align}
\xi_0 (\ff(t_1, t_2), \p(t_1), (\xt, \yt)) + \xt \xi_{1x} (\ff(t_1, t_2), \p(t_1), (\xt, \yt)) &  \nonumber \\ 
  \qquad \qquad + \yt \xi_{1y} (\ff(t_1, t_2), \p(t_1), (\xt, \yt)) & \nonumber \\ 
+ \frac{f_{20}(t_1, t_2)}{2} x_{t_2}^2 + f_{11}(t_1, t_2) x_{t_2} y_{t_2} +  
\frac{f_{02}(t_1, t_2)}{2} y_{t_2}^2 + \ldots &=0  \label{eval_perturbed_number_again} \\ 
\R(\ff(t_1, t_2), \p(t_1), (\xt, \yt)) &=0 \label{relation_number_generic} \\ 
\xi_{1x} (\ff(t_1, t_2), \p(t_1), (\xt, \yt))  +  f_{20}(t_1, t_2) \xt + f_{11}(t_1, t_2) y_{t_2} + \ldots &=0 \label{fx_perturbed_number} \\ 
\xi_{1y} (\ff(t_1, t_2), \p(t_1), (\xt, \yt))  + f_{11}(t_1, t_2) x_{t_2} + f_{02}(t_1, t_2) y_{t_2} + \ldots &= 0 \label{fy_perturbed_number} \\ 
\textnormal{$(\xt, \yt)=$  small}, \qquad \big|f_{20}(t_1, t_2) f_{02}(t_1, t_2) - f_{11}(t_1, t_2)^2\big| & > \mathrm{C}   \label{psi_A1_contribution_from_A1_numbers}
\end{align}
where $\mathrm{C}$ is a small positive constant. Observe that the number of solutions 
of \eqref{psi_A1_contribution_from_A1} is same as the number of solutions of 
\eqref{eval_perturbed_number_again}-\eqref{psi_A1_contribution_from_A1_numbers}. 
Let $\mathrm{N}$ be the number of solutions $(\ff, \p)$ of     
\begin{align}
\xi_{0}(\ff, \p, 0, 0) =0, ~~\R(\ff,\p, 0,0) =0, ~~f_{00}=0, ~~ (f_{10}, ~f_{01}) = (0,0).  \label{eval_perturbed_number} 
\end{align}
Observe that this number is same as 
the number of solutions 
\begin{align}
\nu_0(\ff, \p, \p) =0, ~~\Q(\ff, \p, \p) =0, ~~(\ff, \p, \p) \in \Delta \ov{\A}_1.  \label{psi_a0_nu_0_Q_c0}
\end{align}
Hence 
\bge
\mathrm{N} = \big\langle e( \pi_2^* \DL_{\A_0} \oplus \WL), ~[\Delta \ov{\A}_1] \big\rangle. \label{N_equal_Euler}
\ede
Since $\Q$ is generic, all solutions of \eqref{psi_a0_nu_0_Q_c0} 
(and hence \eqref{eval_perturbed_number})
belong to $\Delta \A_1$, i.e,  
\begin{align}
f_{20}  f_{02} - f_{11}^2 \neq 0 \implies \big|f_{20} f_{02} - f_{11}^2\big| > \mathrm{C}. \label{ift_hypothesis_det_hess}
\end{align}
Let $(\ff, \p)$ be a solution of \eqref{eval_perturbed_number}. 
Since the sections 
\[ (\pi_2^*\ds_{\A_0}+\nu_0) \oplus \Q: \ov{\A}_0\times\{\p\} \lra \pi_2^* \DL_{\A_0} \oplus \WL, \qquad  (\pi_2^*\ds_{\A_0} +\nu_0) \oplus \Q: \ov{\A}_0 \times \P^2 \lra \pi_2^* \DL_{\A_0} \oplus \WL \]
are transverse to the zero set (for any $\p \in \P^2$)
we conclude that given an $(\xt, \yt)$ sufficiently small, there exists a unique 
$(\ff(t_1, t_2), \p(t_1)) \in \ov{\A}_0$
close to  $(\ff, \p)$, such that $(\ff(t_1, t_2), \p(t_1), (\xt, \yt))$ 
solves \eqref{eval_perturbed_number_again} and \eqref{relation_number_generic}.  
Plugging in this value in \eqref{fx_perturbed_number} and \eqref{fy_perturbed_number} and 
using \eqref{ift_hypothesis_det_hess}, we conclude 
there is  a unique $(\xt, \yt)$ that solves \eqref{fx_perturbed_number} and \eqref{fy_perturbed_number} (provided the norm of $\xi_{1x}$ and $\xi_{1y}$ is sufficiently small). 
Hence, there is a one to one correspondence between the number solutions of  \eqref{eval_perturbed_number} and the solutions of 
\eqref{eval_perturbed_number_again}-\eqref{psi_A1_contribution_from_A1_numbers}. Equation \eqref{N_equal_Euler} now proves 
\eqref{contrib_from_delta_a1_psi_a1}. \\ 
\hf \hf Next, suppose $(\ff,\p, \p) \in \Delta \A_{2} \cap \Q^{-1}(0)$. 
Since $f_{20}$ and $f_{02}$ are not both zero, let us assume $f_{02} \neq 0$.
The contribution of the section 
\[ \pi_2^* \ds_{\A_{0}} \oplus \pi_2^* \ds_{\A_{1}} \oplus \Q: \ov{\A}_1 \times \P^2 \lra \pi_2^* \DL_{\A_0} \oplus \pi_2^* \DV_{\A_1} \oplus \WL \]
to the Euler class is the number of solutions of  
\begin{align}
\pi_2^*\ds_{\A_0} (\ff(t_1, t_2), \p(t_1, t_2), \p(t_1))&=  \nu_0, ~~\pi_2^*\ds_{\A_1} (\ff(t_1, t_2), \p(t_1, t_2), \p(t_1))=  \nu_1 \nonumber \\ 
\Q(\ff(t_1, t_2), \p(t_1, t_2), \p(t_1)) &= 0, \qquad (\ff(t_1, t_2), \p(t_1, t_2), \p(t_1)) \in \mathcal{U}_{(\ff, \p,\p)} \subset \ov{\A}_1 \times \P^2 \label{contrib_psi_a0_plus_a1_from_cusps}
\end{align}
where $\mathcal{U}_{(\ff, \p,\p)}$ is a sufficiently small neighbourhood of $(\ff, \p, \p)$ inside 
$\ov{\A}_1 \times \P^2$ and $\nu_0$ and $\nu_1$
are generic smooth perturbations. Let $\xi_0$, $\xi_{1x}$, $\xi_{1y}$ and $\R$ be defined 
as in \eqref{number_defn}. 
The number of solutions of \eqref{contrib_psi_a0_plus_a1_from_cusps} 
(which is a functional equation), is equal to the number of 
solutions of \eqref{eval_perturbed_number_again}-\eqref{fy_perturbed_number} and 
\begin{align}
\textnormal{$(\xt, \yt)=$  small}, \qquad  f_{02}(t_1, t_2)\B^{f(t_1, t_2)}_2 := f_{02}(t_1, t_2)f_{20}(t_1, t_2) - f_{11}(t_1, t_2) = \textnormal{small}.   \label{psi_A1_contribution_from_A2_numbers}
\end{align}
Since the section
\[ \ds_{\A_2} \oplus \Q : \Delta \ov{\A}_1 \lra \DL_{\A_2} \oplus \WL  \]
is transverse to the zero set, we conclude that 
there is a unique $(\ff(t_1, t_2), \p(t_1)) \in \Delta \ov{\A}_1 \cap \Q^{-1}(0)$ 
close to $(\ff, \p)$ with a specified value of $\B^{f(t_1, t_2)}_2$.\footnote{This is true provided $\B^{f(t_1, t_2)}_2$ is sufficiently small.} 
In other words, we can express all the 
$f_{ij}(t_1, t_2)$ in terms of $\B^{f(t_1, t_2)}_2$. 
Plug in this expression for $f_{ij}(t_1, t_2)$ 
in  
\eqref{eval_perturbed_number_again}, \eqref{fx_perturbed_number} and \eqref{fy_perturbed_number}. 
Since $~(\ff(t_1, t_2), \p(t_1, t_2), \p(t_1))$  
is close to $\Delta \A_2$,   
we conclude that after a change of coordinates, the set of equations  
\eqref{eval_perturbed_number_again}, \eqref{fx_perturbed_number} and \eqref{fy_perturbed_number}
%
is equivalent to   
\eqref{closure_a1_a1_is_not_a2_F}, 
\eqref{closure_a1_a1_is_not_a2_Fx} and \eqref{closure_a1_a1_is_not_a2_Fy}, 
with  
\begin{align}
f_{00}(t_1, t_2) & = \xi_{0}, \qquad  f_{10}(t_1, t_2) = \xi_{1x}, \qquad f_{01}(t_1, t_2) = \xi_{1y}. \label{psi_a1_contribution_from_a2}
\end{align}
Hence, \eqref{closure_a1_a1_is_not_a2_F_eliminated} holds which combined with \eqref{psi_a1_contribution_from_a2} 
implies that the multiplicity is $3$ to one in $x_{t_2}$. 
Given $x_{t_2}$, we can solve for $\hat{y}_{t_2}$ uniquely using \eqref{closure_a1_a1_is_not_a2_F_hat_y}. 
And given $x_{t_2}$ and $\hat{y}_{t_2}$, we can uniquely solve for  $\B^{f(t_1, t_2)}_2$  
using \eqref{closure_a1_a1_is_not_a2_Fx}, provided $\xi_{1x}$ and $\xi_{1y}$ are sufficiently small. 
Hence, the total multiplicity is $3$. \qed \\

\textbf{Proof of Lemma \ref{cl_two_pt} (\ref{a1a1_up_cl}):} It suffices to show  that 
\bge
\Big\{ (\ff, \p, \lp) \in \ov{\ov{\A}_1 \circ \hat{\A}^{\#}_1} \Big\} = \Delta \ov{\hat{\A}}_{3}. \label{a1_hat_a1_sharp_is_a3_sharp}  
\ede
By using Lemma \ref{tube_lemma} with Lemma \ref{Dk_sharp_closure} \eqref{a1_cl_new} we get $\ov{\hat{\A}}_1=\hat{\A}_1\cup\ov{\hat{\A}}_2$. 
Now we use Lemma \ref{tube_lemma} again with Lemma \ref{cl} \eqref{A2cl} to get $\ov{\hat\A}_2=\hat{\A}_2\sqcup \ov{\hat\A}_3 \cup \ov {\hat{\D}}_4$. 
By Lemma \ref{cl} \eqref{A3cl} and \ref{tube_lemma}, we conclude that $\ov{\hat{\D}}_4 $ is a subset of $\ov{\hat{\A}}_3$. 
Hence $\ov{\hat\A}_2=\hat{\A}_2\sqcup \ov{\hat\A}_3$. 
This implies that
\bgd
\Delta \ov{\hat{\A}}_1=\Delta \hat{\A}_1\sqcup \Delta\hat{\A}_2\sqcup \Delta\ov{\hat{\A}}_3.
\edd
First we observe that the lhs  of \eqref{a1_hat_a1_sharp_is_a3_sharp} 
is a subset of its rhs. To see this, observe that  by \eqref{a1_du_a2_intersect_a1_a1_is_empty}
\begin{align*}
\{ (\ff, \p, \p) \in \ov{\ov{\A}_1 \circ \ov{\A}_1} \} \cap \Delta (\A_1 \cup \A_{2}) & = \varnothing \\
\implies \{ (\ff, \p, \lp) \in \ov{\ov{\A}_1 \circ \ov{\hat{\A}}_1} \} \cap \Delta (\hat{\A}_1 \cup \hat{\A}_{2}) & = \varnothing \\
\implies \{ (\ff, \p, \lp) \in \ov{\ov{\A}_1 \circ \hat{\A}^{\#}_1} \} \cap \Delta (\hat{\A}_1 \cup \hat{\A}_{2}) &= \varnothing  \qquad \textnormal{since} ~\ov{\hat{\A}^{\#}_1} = \ov{\hat{\A}}_1\,\,(\textup{Lemma}\, \,\ref{cl}\,\eqref{A1cl}). 
\end{align*} 

\ni Next we will show that the rhs  of \eqref{a1_hat_a1_sharp_is_a3_sharp} is a 
subset of its lhs. We will simultaneously prove the following two statements. 
\begin{align}
\{ (\ff, \p, \lp) \in \ov{\ov{\A}_1 \circ \hat{\A}_1^{\#}} \} &\supset \Delta (\hat{\A}_{3} \du \hat{\D}_4^{\#\flat}),  \label{a3_is_subset_of_a1_a1_closure_down_stairs}   \\
\{ (\ff, \p, \lp) \in \ov{\ov{\A}_1 \circ \PP \A_2} \} \cap \Delta (\hat{\A}_{3} \du \hat{\D}_4^{\#\flat}) &= \varnothing. \label{a3_interesct_a1_pa2_is_empty_set_equation}
\end{align} 

\ni Since $\ov{\ov{\A}_1 \circ \hat{\A}_1^{\#}}$ is a closed set, \eqref{a3_is_subset_of_a1_a1_closure_down_stairs} implies that 
the rhs  of \eqref{a1_hat_a1_sharp_is_a3_sharp} is a subset of its lhs.\footnote{As before, we do not need the full strength of 
\eqref{a3_is_subset_of_a1_a1_closure_down_stairs}. We simply need that 
$\{ (\ff, \p, \lp) \in \ov{\ov{\A}_1 \circ \hat{\A}_1^{\#}} \} \supset \Delta \hat{\A}_{3}$.}

\begin{claim}
\label{claim_a3_subset_of_a1_a1_closure}
Let $~(\ff,\p, \lp) \in \Delta (\hat{\A}_{3} \du \hat{\D}_4^{\#\flat})$.
Then there exist solutions 
$$ (\ff(t_1, t_2), \p(t_1, t_2),  l_{\p(t_1)}) \in \ov{(\D \times \P^2) \circ \hat{\A}_1^{\#}}$$ 
close to  $(\ff, \p, \lp)$ to the set of equations 
\begin{align}
\label{a3_is_subset_of_a1_a1_closure_down_stairs_equation}
\pi_1^* \ds_{\A_0}(\ff(t_1, t_2), \p(t_1, t_2),  l_{\p(t_1)}) & = 0, ~\pi_1^* \ds_{\A_1}(\ff(t_1, t_2), \p(t_1, t_2),  l_{\p(t_1)})= 0, ~~\p(t_1, t_2) \neq \p(t_1). 
\end{align}
Moreover, whenever such a solution is sufficiently close to $(\ff,\p, \lp)$, it lies in $\ov{\A}_1 \circ \hat{\A}_1^{\#}$, i.e., 
\begin{align}
\pi_2^* \us_{\PP \A_2}(\ff(t_1, t_2), \p(t_1, t_2),  \p(t_1)) & \neq 0. \label{pa3_intersect_a1_a2_is_empty}
\end{align}
In particular, $(\ff(t_1, t_2), \p(t_1, t_2),  \p(t_1))$ does not lie in $\ov{\A}_1 \circ \PP \A_2$.
\end{claim}
It is easy to see that claim \ref{claim_a3_subset_of_a1_a1_closure} proves  \eqref{a3_is_subset_of_a1_a1_closure_down_stairs} and \eqref{a3_interesct_a1_pa2_is_empty_set_equation} 
simultaneously. \\

\pf Choose homogeneous coordinates $[\mathrm{X}: \mathrm{Y}: \mathrm{Z}]$ so that 
$\p = [0:0:1]$ and let ~$\mathcal{U}_{\p}$,   
$\pi_x$, $\pi_y$,  $v$, $w$, 
$x_{t_1}$, $y_{t_1}$, $x_{t_2}$, $y_{t_2}$,
$f_{ij}(t_1, t_2)$
be exactly the same as defined in the 
proof of claim \ref{a1_a1_closure_intersect_a1_or_a2_empty_equations_claim}.
Take
\[(\ff (t_1, t_2),  l_{\p(t_1)}) \in \ov{\hat{\A}^{\#}_1}\] 
to be a point that is close to 
$( \ff, l_{\p})$ and $l_{\p(t_1, t_2)}$ a point 
in $\P T\P^2$ that is close to $l_{\p(t_1)}$.
Without loss of generality, we can assume that   
\[ v + \eta w \in l_{\p}, ~~v + \eta_{t_1} w \in l_{\p(t_1)} ~~\textnormal{and} ~~v + (\eta_{t_1} + \eta_{t_2}) w \in l_{\p(t_1, t_2)} \] 
for some complex 
numbers $\eta$, $\eta_{t_1}$ and $\eta_{t_1} +\eta_{t_2}$ close to each other. 
Let the numbers 
$\FF$, $\Fx$ and $\Fy$ be the same as in the proof of  claim \ref{a1_a1_closure_intersect_a1_or_a2_empty_equations_claim}.
Since  $(\ff (t_1, t_2),  l_{\p(t_1)}) \in \ov{\hat{\A}^{\#}_1}$, we conclude that 
\[ f_{00}(t_1, t_2) =f_{10}(t_1, t_2) = f_{01}(t_1, t_2) =0.\] 
The functional equation   
\eqref{a3_is_subset_of_a1_a1_closure_down_stairs_equation}  
has a solution if and only if the following has a numerical solution: 
\begin{align}
\mathrm{F} = 0, \qquad \Fx = 0, \qquad \Fy = 0,  \qquad (x_{t_2}, y_{t_2}) \neq (0, 0) \qquad \textnormal{(but small)}. \label{eval_f1_again_d4_hat} 
\end{align} 
First let us assume $(\ff, \p, \lp) \in \Delta \hat{\A}_3$. 
It is easy to see that $f_{20}$ and $f_{02}$ can not both be zero; 
let us assume $f_{02}$ is non zero.  
Following the same argument as in the proof of claim \ref{a1_a1_closure_intersect_a1_or_a2_empty_equations_claim}, 
we make a change of coordinates and write $\mathrm{F}$ as 
\begin{align*}
\mathrm{F} & = \hat{\hat{y}}_{t_2}^2 +\frac{\B^{f(t_1, t_2)}_{2}}{2!} x_{t_2}^2 + \frac{\B^{f(t_1, t_2)}_{3}}{3!} x_{t_2}^3 + 
\frac{\B^{f(t_1, t_2)}_{4}}{4!} x_{t_2}^4 + \mathrm{O}(x_{t_2}^5), \\ 
\textnormal{where} & \qquad \hat{\hat{y}}_{t_2} := \sqrt{\varphi ( x_{t_2}, \hat{y}_{t_2})} \hat{y}_{t_2}, \qquad \B^{f(t_1, t_2)}_2 = f_{20}(t_1, t_2) - \frac{f_{11}(t_1, t_2)^2}{f_{02}(t_1, t_2)}.  
\end{align*}
It is easy to see that the only solutions to \eqref{eval_f1_again_d4_hat}  (in terms of the new coordinates) are
\bge
\B_{2}^{f(t_1, t_2)}  = \ {\textstyle \frac{\B_{4}^{f(t_1, t_2)}}{12}} x_{t_2}^2  + \mathrm{O}(x_{t_2}^3), ~\B_{3}^{f(t_1, t_2)} = {\textstyle \frac{-\B_{4}^{f(t_1, t_2)}}{2}} x_{t_2}+\mathrm{O}(x_{t_2}^2), ~\hat{\hat{y}}_{t_2} =0, ~x_{t_2} \neq 0 ~\textnormal{(but small).} \label{b20_equation} 
\ede
It remains to show that these solutions satisfy \eqref{pa3_intersect_a1_a2_is_empty}. 
First consider the case when $(\ff, \p,  \lp) \notin \PP \A_3 $, i.e., 
$\pi_2^* \us_{\PP \A_2}(\ff, \p,  \lp)  \neq 0$.
Then \eqref{pa3_intersect_a1_a2_is_empty} is obviously true, since the section $\pi_2^* \us_{\PP \A_2}$ is \textit{continuous}. 
Next, consider the case when $(\ff, \p,  \lp) \in \PP \A_3$. 
Define the numbers 
\begin{align}
\mathrm{J}_1 &:= f_{20}(t_1, t_2) + \eta_{t_1} f_{11}(t_1, t_2), ~~ \mathrm{J}_2:= f_{11}(t_1, t_2) + \eta_{t_1} f_{02}(t_1, t_2). \label{b2_can_not_be_zero}
\end{align}
Observe that 
\begin{align}
\{ \pi_2^* \us_{\PP \A_2}(\ff(t_1, t_2), \p(t_1, t_2),  l_{\p(t_1)} ) \} (f(t_1, t_2) \otimes p(t_1)^{\otimes d} \otimes (v + \eta_{t_1} w) \otimes v) &= \mathrm{J}_1, \nonumber  \\   
\{ \pi_2^* \us_{\PP \A_2}(\ff(t_1, t_2), \p(t_1, t_2),  l_{\p(t_1)} ) \} (f(t_1, t_2) \otimes p(t_1)^{\otimes d} \otimes (v + \eta_{t_1} w) \otimes w) &= \mathrm{J}_2.  
\label{psi_pa2_number_form}
\end{align}
If \eqref{pa3_intersect_a1_a2_is_empty} were false, then $\mathrm{J}_1$ and $\mathrm{J}_2$ would vanish (by \eqref{psi_pa2_number_form}). 
Equations \eqref{b2_can_not_be_zero} and \eqref{b20_equation} imply
\begin{align}
\frac{\B^{f(t_1, t_2)}_{4}}{12} x_{t_2}^2 + \mathrm{O}(x_{t_2}^3) &= \mathrm{J}_1 - \frac{f_{11}(t_1, t_2)}{f_{02}(t_1, t_2)} \mathrm{J}_2.   \label{b2_can_not_be_zero_again}
\end{align}
Since $(\ff, \p, \lp) \in \hat{\A}_3$, 
we conclude  $\B^{f(t_1, t_2)}_{4} \neq 0$. Hence, 
\eqref{b2_can_not_be_zero_again} implies that 
$\mathrm{J}_1$ and $\mathrm{J}_2$ can not both be zero. 
Hence \eqref{pa3_intersect_a1_a2_is_empty} holds.  \\
\hf \hf Next, let us assume that $(\ff, \p, \lp) \in \Delta \hat{\D}_4^{\#\flat}$. 
Define the following number 
\begin{align*}
\mathrm{G}&:= x_{t_2} \mathrm{F}_{x_{t_2}} + y_{t_2} \mathrm{F}_{y_{t_2}} -2 \mathrm{F} \\
          & = \frac{f_{30}(t_1, t_2)}{6} x_{t_2}^3 + \frac{f_{21}(t_1, t_2)}{2} x_{t_2}^2 y_{t_2} +  \frac{f_{12}(t_1, t_2)}{2} x_{t_2} y_{t_2}^2 + \frac{f_{03}(t_1, t_2)}{6} y_{t_2}^3 + \ldots.
\end{align*}
Note that the cubic term of $\mathrm{G}$ is same as the cubic term of $\FF$. 
Since $(\ff, \p, \lp) \in \Delta\hat{\D}_4$  
we conclude, using the same argument as in \cite{BM13}
that there exists a change of coordinates  
\begin{align*}
x_{t_2} & = \hat{x}_{t_2} +\mathrm{E}_1 (\hat{x}_{t_2}, \hat{y}_{t_2}), \qquad y_{t_2} = \hat{y}_{t_2} + \mathrm{E}_2 (\hat{x}_{t_2}, \hat{y}_{t_2}),  
\end{align*}
(where $\mathrm{E}_i(\hat{x}_{t_2}, \hat{y}_{t_2})$ are second order in $\hat{x}_{t_2}$ and $\hat{y}_{t_2}$)
so that  $\mathrm{G}$ is given by 
\begin{align}
\mathrm{G} & =\frac{f_{30}(t_1, t_2)}{6} \hat{x}_{t_2}^3 + \frac{f_{21}(t_1, t_2)}{2} \hat{x}_{t_2}^2 \hat{y}_{t_2} +  \frac{f_{12}(t_1, t_2)}{2} \hat{x}_{t_2} \hat{y}_{t_2}^2 + 
\frac{f_{03}(t_1, t_2)}{6} \hat{y}_{t_2}^3.  \label{new_G}
\end{align}
Since $(\ff, \p, \lp) \in \Delta\hat{\D}_4$, 
there 
are three possibilities to consider;  
\begin{align}
f_{30}(t_1, t_2) & \neq 0  \qquad \textnormal{or} \qquad f_{03}(t_1, t_2) \neq 0 ~~ \textnormal{or} \nonumber \\ 
f_{30}(t_1, t_2) & = f_{03}(t_1, t_2) =0, ~~~~\textnormal{but} ~~~~f_{21}(t_1, t_2) \neq 0 ~~\textnormal{and} ~~f_{12}(t_1, t_2) \neq 0. \label{fij_cases_D4_hat}
\end{align}
Let us assume $f_{30}(t_1, t_2)\neq 0$. 
Since $(\ff, \p, \lp) \in \Delta\hat{\D}_4^{\#\flat}$, 
equation \eqref{new_G} now can be written as   
\begin{align}
\mathrm{G}&=  
\frac{f_{30}(t_1,t_2)}{6}(\hat{x}_{t_2} - \mathrm{A}_1 \hat{y}_{t_2}) (\hat{x}_{t_2} - \mathrm{A}_2 \hat{y}_{t_2}) (\hat{x}_{t_2} - \mathrm{A}_3 \hat{y}_{t_2})
\end{align}
where $\mathrm{A}_i$ are complex numbers such that 
\begin{align}
 \mathrm{A}_1 & \neq \mathrm{A}_2 \neq \mathrm{A}_3\neq \mathrm{A}_1 ~~\textnormal{and} ~~\eta \neq \frac{1}{\mathrm{A}_1}, ~\frac{1}{\mathrm{A}_2} ~~\textnormal{or} ~\frac{1}{\mathrm{A}_3}. 
 \label{eta_neq_d4_sharp}
\end{align}
To see why the last  inequality is true,   
note that if $\eta$ is either $\frac{1}{\mathrm{A}_1}$, 
$\frac{1}{\mathrm{A}_2}$ or $\frac{1}{\mathrm{A}_3}$ then 
\[ \nabla^3 f|_{\p} (v+ \eta w, v + \eta w, v + \eta w) =0. \]
Since 
$(\ff, \p, \lp) \in \Delta \hat{\D}_{4}^{\#\flat}$ 
the 
last inequality of \eqref{eta_neq_d4_sharp}  holds. 
Hence,   
\eqref{a3_is_subset_of_a1_a1_closure_down_stairs_equation}  
has a solution if and only if 
\begin{align}
\mathrm{G} = &\,\,  \frac{f_{30}(t_1,t_2)}{6}(\hat{x}_{t_2} - \mathrm{A}_1 \hat{y}_{t_2}) (\hat{x}_{t_2} - \mathrm{A}_2 \hat{y}_{t_2}) (\hat{x}_{t_2} - \mathrm{A}_3 \hat{y}_{t_2}) = 0, \nonumber \\ 
\mathrm{F}_{x_{t_2}} =&\,\,  \hat{x}_{t_2} f_{20}(t_1, t_2) + \hat{y}_{t_2} f_{11}(t_1, t_2) \nonumber \\ 
                      & + \frac{f_{30}(t_1, t_2)}{6} \Big( 3 \hat{x}_{t_2}^2 
                      -2 \hat{x}_{t_2} \hat{y}_{t_2}\big ({\textstyle \sum_{i=1}^3\mathrm{A}_i}\big) + 
                                                                (\mathrm{A}_1\mathrm{A}_2 + \mathrm{A}_1\mathrm{A}_3+ \mathrm{A}_2\mathrm{A}_3 ) 
                                                                \hat{y}_{t_2}^2 \Big)+ 
                                                                \mathrm{E}_3(\hat{x}_{t_2}, \hat{y}_{t_2}) =0, \nonumber \\
\mathrm{F}_{y_{t_2}} =&\,\, \hat{x}_{t_2} f_{11}(t_1, t_2) + \hat{y}_{t_2} f_{02}(t_1, t_2) \nonumber \\ 
                      & -\frac{f_{30}(t_1, t_2)}{6}\Big(\hat{x}_{t_2}^2\big({\textstyle\sum_{i=1}^3  \mathrm{A}_i}\big) 
                      -2 \hat{x}_{t_2} \hat{y}_{t_2}\big({\textstyle \sum_{i\neq j}\mathrm{A}_i\mathrm{A}_j}\big)  
                                                                +3\hat{y}_{t_2}^2\mathrm{A}_1\mathrm{A}_2 \mathrm{A}_3 \Big) + \mathrm{E}_4(\hat{x}_{t_2}, \hat{y}_{t_2}) =0 \nonumber \\  
(\hat{x}_{t_2}, \hat{y}_{t_2}) \neq &\,\, (0,0) \qquad \textnormal{(but small)}  \label{hat_D4_neighbourhood_inside_a1_hat_a1}
\end{align}
has a solution, where $\mathrm{E}_i(\hat{x}_{t_2}, \hat{y}_{t_2})$ are third order in $(\hat{x}_{t_2}, \hat{y}_{t_2})$. \\  
\hf \hf To avoid confusion let us clarify one point; in the above equation 
$\mathrm{F}_{x_{t_2}}$ 
and $\mathrm{F}_{y_{t_2}}$  
are simply  
expressed in terms of the new coordinates 
$\hat{x}_{t_2}$ and $\hat{y}_{t_2}$. They are still the partial derivatives of 
$\FF$ with respect to $\x_{t_2}$ and $\y_{t_2}$;  
they are \textit{not} $\mathrm{F}_{\hat{x}_{t_2}}$ 
and $\mathrm{F}_{\hat{y}_{t_2}}$, the  partial derivatives of  $\mathrm{F}$ with respect to 
$\hat{x}_{t_2}$ and $\hat{y}_{t_2}$. Now we will construct the solutions to \eqref{hat_D4_neighbourhood_inside_a1_hat_a1}. 
There are three solutions; we will just give one of the solutions, the rest are similar.   
They are given by: 
\begin{align}
\hat{x}_{t_2}&= \mathrm{A}_1 \hat{y}_{t_2}, ~~\hat{y}_{t_2} \neq 0  ~~\textnormal{(but small)}, ~~f_{20}(t_1, t_2) = \textnormal{small}, \nonumber \\
f_{02}(t_1, t_2) & = \frac{f_{30}(t_1, t_2)}{6}\Big( 
2 \mathrm{A}_1^3 - 2 \mathrm{A}_1^2 \mathrm{A}_2 - \mathrm{A}_1^2 \mathrm{A}_3 + 2 \mathrm{A}_1 \mathrm{A}_2\mathrm{A}_3  \Big) \hat{y}_{t_2} 
+ \mathrm{A}_1^2 f_{20}(t_1, t_2) + 
\mathrm{E}_5(\hat{y}_{t_2}), \nonumber \\  
f_{11}(t_1, t_2) & = \frac{f_{30}(t_1, t_2)}{6}\Big( -\mathrm{A}_1^2 +\mathrm{A}_1 \mathrm{A}_2 +\mathrm{A}_1 \mathrm{A}_3 - \mathrm{A}_2 \mathrm{A}_3 \Big)\hat{y}_{t_2} 
-\mathrm{A}_1 f_{20}(t_1, t_2) + 
\mathrm{E}_6(\hat{y}_{t_2}), \label{pa2_section_around_hat_d4_sharp} 
\end{align}
where $\mathrm{E}_i(\hat{y}_{t_2})$ are second order in $\hat{y}_{t_2}$ and independent of $f_{20}(t_1, t_2)$.  
It remains to show that \eqref{pa3_intersect_a1_a2_is_empty} holds.
Equation \eqref{pa2_section_around_hat_d4_sharp} implies that 
\begin{align}
 & (1 -\eta_{t_1} \mathrm{A}_1)\mathrm{J}_2 + (\mathrm{A}_1 - \mathrm{A}_1^2 \eta_{t_1}) \mathrm{J}_1 = \beta \hat{y}_{t_2} + \mathrm{O}(\hat{y}_{t_2}^2), \label{pa2_section_around_hat_d4_sharp_again}\\
\textnormal{where} \qquad  \beta &:= \frac{f_{30}(t_1, t_2)}{6}(\mathrm{A}_1 - \mathrm{A}_2) (\mathrm{A}_1 - \mathrm{A}_3)(-1 + \mathrm{A}_1 \eta_{t_1})^2, \nonumber 
\end{align}
and $\mathrm{J}_1$ and $\mathrm{J}_2$ are as defined in \eqref{b2_can_not_be_zero}.
Note that the rhs of \eqref{pa2_section_around_hat_d4_sharp_again} 
is independent of $f_{20}(t_1, t_2)$. By \eqref{eta_neq_d4_sharp}, 
$\beta \neq 0$. Hence, by \eqref{pa2_section_around_hat_d4_sharp_again} $\mathrm{J}_1$ and $\mathrm{J}_2$ 
can not both vanish. As a result, \eqref{pa3_intersect_a1_a2_is_empty} holds. 
Similar argument holds for the other two solutions of \eqref{hat_D4_neighbourhood_inside_a1_hat_a1}. \\
\hf \hf A similar argument will go through if 
any of the 
other two cases of \eqref{fij_cases_D4_hat} holds. \qed 

\begin{cor}
\label{pa2_section_mult_around_pa3}
Let $\W \lra \D \times \P^2 \times \P T\P^2$ be a vector bundle such that 
the rank of $\W$ is equal to dimension of $ \Delta \PP \A_{3}$ and 
$\Q: \D \times \P^2 \times \P T\P^2 \lra \W$  a \textit{generic} 
smooth section. Suppose $(\ff, \p, \lp) \in \Delta \PP \A_3$. 
Then the section 
$$ \pi_2^*\us_{\PP \A_{2}} \oplus \Q: \ov{\ov{\A}_1 \circ \hat{\A}_1^{\#}}  \lra \pi_2^* \UV_{\PP \A_2} \oplus \W$$ 
vanishes around $(\ff, \p, \lp)$ with a multiplicity of $2$.   
\end{cor}
\pf Since $\Q$ is generic, the sections 
\[ \us_{\PP \A_2} \oplus \Q: \Delta \ov{\hat{\A}^{\#}_1} \lra \UV_{\PP \A_2} \oplus \W, 
~~ \us_{\PP \A_3}: \us_{\PP \A_2}^{-1}(0) \lra \UL_{\PP \A_3} \]
are transverse to the zero set. 
Hence, there exists a unique 
$(\ff(t_1, t_2), l_{p(t_1)}) \in  \Delta \ov{\hat{\A}^{\#}_1}$ close to $(\ff, \lp)$ for a specified 
value of $\mathrm{J}_1$, $\mathrm{J}_2$ and $f_{30}(t_1, t_2)$. \footnote{Provided 
$\mathrm{J}_1$, $\mathrm{J}_2$ and $f_{30}(t_1, t_2)$ are sufficiently small.} 
In other words we can express all the $f_{ij}(t_1, t_2)$ in terms of 
$\mathrm{J}_1$, $\mathrm{J}_2$ and $f_{30}(t_1, t_2)$. 
Since $\B^{f(t_1, t_2)}_4 \neq 0$, 
equation \eqref{b2_can_not_be_zero_again} implies that 
the number of 
solutions to the set of equations 
\[ \mathrm{J}_1 = \xi_1, \qquad \mathrm{J}_2 = \xi_2 \] 
is $2$, where $\xi_i$ is a small perturbation. \qed   

\begin{cor}
\label{pa2_section_mult_around_hat_d4}
Let $\W \lra \D \times \P^2 \times \P T\P^2$ be a vector bundle such that 
the rank of $\W$ is equal to dimension of $ \Delta \hat{\D}_{4}^{\#\flat}$ and 
$\Q: \D \times \P^2 \times \P T\P^2 \lra \W$  a \textit{generic} 
smooth section.  Suppose $(\ff, \p, \p) \in \Delta \hat{\D}_4^{\#\flat}$. 
Then the section 
$$ \pi_2^*\us_{\PP \D_{4}} \oplus \Q: \ov{\ov{\A}_1 \circ \hat{\A}_1^{\#}}  \lra \pi_2^* \UL_{\PP \D_4} \oplus \W$$ 
vanishes around $(\ff, \p, \lp)$ with a multiplicity of $3$.   
\end{cor}

\pf Since $\Q$ is generic, the sections 
\[ \us_{\PP \A_2} \oplus \Q: \Delta \ov{\hat{\A}^{\#}_1} \lra \UV_{\PP \A_2} \oplus \W, 
~~ \us_{\PP \D_4}: \us_{\PP \A_2}^{-1}(0) \lra \UL_{\PP \D_4} \]
are transverse to the zero set. 
Hence, there exists a unique 
$(\ff(t_1, t_2), l_{p(t_1)}) \in  \Delta \ov{\hat{\A}^{\#}_1}$ close to $(\ff, \lp)$ for a specified 
value of $\mathrm{J}_1$, $\mathrm{J}_2$ and $f_{02}(t_1, t_2)$. \footnote{Provided 
$\mathrm{J}_1$, $\mathrm{J}_2$ and $f_{02}(t_1, t_2)$ are sufficiently small.} 
In other words we can express all the $f_{ij}(t_1, t_2)$ in terms of 
$\mathrm{J}_1$, $\mathrm{J}_2$ and $f_{02}(t_1, t_2)$. Since $\beta \neq 0$,  
equation 
\eqref{pa2_section_around_hat_d4_sharp_again} 
implies that 
the number of 
solutions to the set of equations 
\[ \mathrm{J}_1 = \xi_1, \qquad \mathrm{J}_2 = \xi_2 \] 
is $1$, where $\xi_i$ is a small perturbation.  
Since there are a total of $3$ solutions to \eqref{hat_D4_neighbourhood_inside_a1_hat_a1}, 
the total multiplicity is $3$. \qed   \\


\textbf{Proof of Lemma \ref{cl_two_pt} (\ref{a1_pa2_cl}):} It suffices to prove the following two statements in view of \eqref{pak2_is_subset_of_a1_and_pak} :  
\begin{align}
\{ (\ff, \p, l_{\p}) \in \ov{\ov{\A}_1\circ \PP \A}_2: ~\pi_2^*\us_{\PP \D_4}( \ff, \p, l_{\p} ) \neq 0  \} &=  
\{ (\ff, \p, l_{\p}) \in \Delta \ov{\PP \A}_{4}: ~\pi_2^*\us_{\PP \D_4}( \ff, \p, l_{\p} ) \neq 0  \}
\label{closure_a1_pa2_f02_not_zero}  \\
\{ (\ff, \p, l_{\p}) \in \ov{\ov{\A}_1 \circ \PP \A}_2: ~\pi_2^*\us_{\PP \D_4}( \ff, \p, l_{\p} ) = 0  \} &=  \Delta \ov{\hat{\D}^{\#\flat}_5}   \label{a1_pa2_d5}
 \end{align}
Let us directly prove a more general version of \eqref{closure_a1_pa2_f02_not_zero}: 
\begin{lmm}
\label{closure_a1_pak_f02_not_zero}
If $k \geq 2$, then  
\begin{align*}
\{ (\ff, \p, l_{\p}) \in \ov{\ov{\A}_1 \circ  \PP \A}_k: ~\pi_2^*\us_{\PP \D_4}( \ff, \p, l_{\p} ) \neq 0  \} &= 
\{ (\ff, \p, l_{\p}) \in \Delta \ov{\PP \A}_{k+2}: ~\pi_2^*\us_{\PP \D_4}( \ff, \p, l_{\p} ) \neq 0  \}.  
\end{align*}
\end{lmm}
Note that \eqref{closure_a1_pa2_f02_not_zero} is a special case of Lemma \ref{closure_a1_pak_f02_not_zero}; take $k=2$.\\

\pf We will prove the following two facts simultaneously: 
\begin{align}
\{ (\ff, \p, \lp) \in \ov{\ov{\A}_1 \circ \PP \A}_{k} \} \supset \Delta \PP \A_{k+2} \qquad & \forall ~k \geq 2, \label{pak2_is_subset_of_a1_and_pak} \\
\{ (\ff, \p, \lp) \in \ov{\ov{\A}_1 \circ \PP \A}_{k+1} \} \cap \Delta \PP \A_{k+2} = \varnothing \qquad  & \forall ~k \geq 1.  \label{pak2_intersect_a1_and_pa1k+1_is_empty}
\end{align}
It is easy to see that \eqref{pak2_is_subset_of_a1_and_pak} and \eqref{pak2_intersect_a1_and_pa1k+1_is_empty} imply 
Lemma \ref{closure_a1_pak_f02_not_zero}. We will now prove the following claim: 
\begin{claim}
\label{claim_a4_closure_simultaneous}
Let $~(\ff,\p, \lp) \in \Delta \PP \A_{k+2}$ and $ k\geq 2$.
Then there exists a solution 
$$ (\ff(t_1, t_2), \p(t_1, t_2),  l_{\p(t_1)} ) \in \ov{ (\D \times \P^2) \circ \PP \A}_2$$ 
\textit{near} $(\ff, \p, \lp)$ to the set of equations
\begin{align}
\label{closure_a1_ak_hessian_not_zero}
\pi_1^* \ds_{\A_0}(\ff(t_1, t_2), \p(t_1, t_2),  l_{\p(t_1)}) & = 0, ~\pi_1^* \ds_{\A_1}(\ff(t_1, t_2), \p(t_1, t_2),  l_{\p(t_1)}) = 0, \nonumber \\
\pi_2^* \us_{\PP \A_3}(\ff(t_1, t_2), \p(t_1, t_2),  l_{\p(t_1)}) & = 0, \ldots, \pi_2^* \us_{\PP \A_k}(\ff(t_1, t_2), \p(t_1, t_2),  l_{\p(t_1)}) = 0, \nonumber \\
~\pi_2^* \us_{\PP \D_4}(\ff(t_1, t_2), \p(t_1, t_2),  l_{\p(t_1)}) & \neq  0, ~\p(t_1, t_2) \neq \p(t_1). 
\end{align}
Moreover, \textit{any} solution $(\ff(t_1, t_2), \p(t_1, t_2),  l_{\p(t_1)})$ sufficiently close to $(\ff, \p, \lp)$ 
lies in $ \ov{\A}_1 \circ \PP \A_k$, i.e.,
\begin{align}
\pi_2^* \us_{\PP \A_{k+1}}(\ff(t_1, t_2), \p(t_1, t_2),  l_{\p(t_1)}) \neq  0. \label{psi_pa_k_plus_1_does_not_vanish}
\end{align}
In particular $(\ff(t_1, t_2), \p(t_1, t_2),  l_{\p(t_1)})$ \textit{does not} lie in $ \ov{\ov{\A}_1 \circ \PP \A}_{k+1} $.
\end{claim}
\ni It is easy to see that claim \ref{claim_a4_closure_simultaneous} implies \eqref{pak2_is_subset_of_a1_and_pak} and \eqref{pak2_intersect_a1_and_pa1k+1_is_empty} 
simultaneously for all $k\geq 2$.
The fact that \eqref{pak2_intersect_a1_and_pa1k+1_is_empty} holds for $k=1$ follows from \eqref{a3_interesct_a1_pa2_is_empty_set_equation} 
(since $\PP \A_3$ is a subset of $\hat{\A}_3$.)\\

\pf Choose homogeneous coordinates $[\mathrm{X}: \mathrm{Y}: \mathrm{Z}]$ so that 
$\p = [0:0:1]$ and let ~$\mathcal{U}_{\p}$,   
$\pi_x$, $\pi_y$,  
$x_{t_1}$, $y_{t_1}$, $x_{t_2}$, $y_{t_2}$
be exactly the same as defined in the 
proof of claim \ref{a1_a1_closure_intersect_a1_or_a2_empty_equations_claim}.
Let $v_1, w:\mathcal{U}_{\p} \lra T\P^2$ be vectors dual to the one 
forms $d\pi_x$ and $d\pi_y$ respectively.
Take
\[(\ff (t_1, t_2),  l_{\p(t_1)}) \in \ov{\PP \A}_2\] 
to be a point that is close to 
$( \ff, l_{\p})$ and $l_{\p(t_1, t_2)}$ a point 
in $\P T\P^2$ that is close to $l_{\p(t_1)}$.
Without loss of generality, we can assume that   
\[ v:= v_1 + \eta w \in l_{\p}, ~~v_1 + \eta_{t_1} w \in l_{\p(t_1)} ~~\textnormal{and} ~~v + (\eta_{t_1} + \eta_{t_2}) w \in l_{\p(t_1, t_2)} \] 
for some complex 
numbers $\eta$, $\eta_{t_1}$ and $\eta_{t_2}$ close to each other. Let 
\begin{align*}
f_{ij}(t_1, t_2) & := \nabla^{i+j} f(t_1, t_2)|_{p(t_1)} 
(\underbrace{v,\cdots v}_{\textnormal{$i$ times}}, \underbrace{w,\cdots w}_{\textnormal{$j$ times}}).
\end{align*}
The numbers 
$\FF$, $\Fx$ and $\Fy$ are the same as in the proof of  claim \ref{a1_a1_closure_intersect_a1_or_a2_empty_equations_claim}. 
Since  $(\ff (t_1, t_2),  l_{\p(t_1)}) \in \ov{\PP \A}_2$, we conclude that 
\[ f_{00}(t_1, t_2) = f_{10}(t_1, t_2) = f_{01}(t_1, t_2) =f_{20}(t_1, t_2) = f_{11}(t_1, t_2)=0.\]
Moreover, since $(\ff, \lp) \in \PP \A_{k+2}$ we conclude that 
$f_{02}$ and $\A^{f}_{k+3}$ are non zero.
Hence $f_{02}(t_1, t_2)$ and $\A^{f(t_1, t_2)}_{k+3}$ are non zero  
if $\ff(t_1, t_2)$ is sufficiently close to $\ff$.
Since $f_{02}(t_1, t_2) \neq 0$ , 
following the same argument as in the proof of claim \ref{a1_a1_closure_intersect_a1_or_a2_empty_equations_claim}, 
we can make a change of coordinates to write $\mathrm{F}$ as   
\begin{align*} 
\mathrm{F}&= \hat{\hat{y}}_{t_2}^2 + \frac{\A^{f(t_1, t_2)}_3}{3!}x_{t_2}^3 + \frac{\A^{f(t_1, t_2)}_4}{4!} x_{t_2}^4 + \ldots 
\end{align*}
The functional equation   \eqref{closure_a1_ak_hessian_not_zero}  has a solution if and only if 
the following set of equations has a solution (as numbers): 
\begin{align}
\label{closure_a1_ak_hessian_not_zero_numbers}
\hat{\hat{y}}_{t_2}^2 + \frac{\A^{f(t_1, t_2)}_3}{3!}x_{t_2}^3 + \frac{\A^{f(t_1, t_2)}_4}{4!} x_{t_2}^4 + \ldots &=0, \qquad 2 \hat{\hat{y}}_{t_2} = 0, \nonumber \\
\frac{\A^{f(t_1, t_2)}_3}{2!}x_{t_2}^3 + \frac{\A^{f(t_1, t_2)}_4}{3!} x_{t_2}^4 + \ldots &= 0, \qquad \A^{f(t_1, t_2)}_3, \ldots, \A^{f(t_1, t_2)}_k = 0, \nonumber  \\
(\hat{\hat{y}}_{t_2}, x_{t_2}) & \neq (0,0)   \qquad \textnormal{(but small).}
\end{align}
It is easy to see that the solutions to \eqref{closure_a1_ak_hessian_not_zero_numbers} 
exist given by   
\begin{align}
\A^{f(t_1, t_2)}_3,& \ldots, \A^{f(t_1, t_2)}_k = 0, \nonumber \\
\A^{f(t_1, t_2)}_{k+1} &=  \frac{\A_{k+3}^{f(t_1, t_2)}}{(k+2)(k+3)} x_{t_2}^2 + \mathrm{O}(x_{t_2}^3),  \label{pak+2_inside_a1_and_pak_solution}\\ 
\A^{f(t_1, t_2)}_{k+2} &= -\frac{2 \A^{f(t_1, t_2)}_{k+3}}{(k+3)} x_{t_2} + \mathrm{O}(x_{t_2}^2), 
\qquad \hat{\hat{y}}_{t_2} = 0, \qquad x_{t_2}  \neq 0 \qquad \textnormal{(but small).} \nonumber 
\end{align}
By \eqref{pak+2_inside_a1_and_pak_solution}, it immediately follows that \eqref{psi_pa_k_plus_1_does_not_vanish} holds (since $\A^{f(t_1, t_2)}_{k+3} \neq 0$). \qed 

\begin{cor}
\label{a1_pak_mult_is_2_Hess_neq_0}
Let $\W \lra \D \times \P^2 \times \P T\P^2$ be a vector bundle such that 
the rank of $\W$ is same as the dimension of $ \Delta \PP \A_{k+2}$ and 
$\Q: \D \times \P^2 \times \P T\P^2 \lra \W$  a \textit{generic} 
smooth section. Suppose $(\ff,\p, \lp) \in \Delta \PP \A_{k+2} \cap \Q^{-1}(0)$. 
Then the section $$ \pi_2^\ast\us_{\PP \A_{k+1}} \oplus \Q: \Delta \ov{\PP \A}_{k} \lra \pi_2^* \UL_{\PP \A_{k+1}} \oplus \W$$
vanishes around $(\ff, \p,  \lp)$ with a multiplicity of $2$.
\end{cor}
\pf  This follows from the fact that the sections 
\begin{align*}
\pi_2^\ast\us_{\PP \A_{i}}:& \Delta \ov{\PP \A}_{i-1} - \pi_2^\ast\us_{\PP \D_4}^{-1}(0) \lra \pi_2^\ast \UL_{\PP \A_{i}}
\end{align*}
are transverse to the zero set for all $3 \leq i \leq k+2$, the fact that 
$\Q$ is generic and \eqref{pak+2_inside_a1_and_pak_solution}. 
The proof is now similar to that of Corollary \ref{a1_section_contrib_from_a1_and_a2}, 
\ref{pa2_section_mult_around_pa3} and \ref{pa2_section_mult_around_hat_d4}. \qed \\

\hf \hf Next let us prove \eqref{a1_pa2_d5}. First we will prove the following two facts: 
\begin{align}
\ov{\ov{\A}_1 \circ \PP \A_2} \cap \PP \D_4 &= \varnothing,\label{a1_pa2_intersect_pd4_empty} \\ 
\ov{\ov{\A}_1 \circ \PP \A_3} \cap \PP \D_5 &= \varnothing. \label{a1_pa3_intersect_pd5_empty} 
\end{align}
\ni Although \eqref{a1_pa3_intersect_pd5_empty} is not needed to prove \eqref{a1_pa2_d5}, we will prove these two 
statements in one go since their proofs are very similar. 

\begin{claim}
\label{a1_pa2_intersect_pd4_empty_claim}
Let $(\ff, \p, \lp) \in \Delta \PP \D_4$. Then there exist no solutions 
\[ (\ff(t_1, t_2), \p(t_1, t_2), l_{\p(t_1)}) \in \ov{\ov{\A}_1 \circ \PP \A}_2 \]
near $(\ff, \p, \lp)$ to the set of equations 
\begin{align}
\pi_1^* \ds_{\A_0}(\ff(t_1, t_2), \p(t_1, t_2),  l_{\p(t_1)}) & = 0, ~\pi_1^* \ds_{\A_1}(\ff(t_1, t_2), \p(t_1, t_2),  l_{\p(t_1)}) = 0. 
\label{a1_pa2_intersect_pd4_empty_equation}
\end{align}
Secondly, let $(\ff, \p, \lp) \in \Delta \PP \D_5$. Then there exist no solutions 
\[ (\ff(t_1, t_2), \p(t_1, t_2), l_{\p(t_1)}) \in \ov{\ov{\A}_1 \circ \PP \A}_3 \]
near $(\ff, \p, \lp)$ to the set of equations 
\begin{align}
\pi_1^* \ds_{\A_0}(\ff(t_1, t_2), \p(t_1, t_2),  l_{\p(t_1)}) & = 0, ~\pi_1^* \ds_{\A_1}(\ff(t_1, t_2), \p(t_1, t_2),  l_{\p(t_1)}) = 0. 
\label{a1_pa3_intersect_pd5_empty_equation}
\end{align}
\end{claim}
\ni It is easy to see that  claim \ref{a1_pa2_intersect_pd4_empty_claim} proves \eqref{a1_pa2_intersect_pd4_empty} 
and \eqref{a1_pa3_intersect_pd5_empty} .\\

\pf For the first part, choose homogeneous coordinates $[\mathrm{X}: \mathrm{Y}: \mathrm{Z}]$ so that 
$\p = [0:0:1]$ and let ~$\mathcal{U}_{\p}$,   
$\pi_x$, $\pi_y$,  $v_1$, $w$, $v$, $\eta$, $\eta_{t_1}$, $\eta_{t_2}$, 
$x_{t_1}$, $y_{t_1}$, $x_{t_2}$, $y_{t_2}$, 
$f_{ij}(t_1, t_2)$, $\FF$, $\Fx$ and $\Fy$
be exactly the same as defined in the 
proof of claim \ref{claim_a4_closure_simultaneous}.
Since $(\ff(t_1, t_2), l_{\p(t_1)}) \in \ov{\PP \A}_2$, we conclude that 
\[ f_{00}(t_1, t_2) = f_{10}(t_1, t_2) = f_{01}(t_1, t_2) =f_{20}(t_1, t_2) = f_{11}(t_1, t_2)=0.\]
The functional equation \eqref{a1_pa2_intersect_pd4_empty_equation}  has a solution if and only if 
the following set of equations has a solution (as numbers): 
\begin{align}
\mathrm{F} = 0, \qquad \Fx = 0, \qquad \Fy = 0,  \qquad (x_{t_2}, y_{t_2}) \neq (0, 0) \qquad \textnormal{(but small)}. \label{eval_f1_pd4} 
\end{align} 
For the convenience of the reader, let us rewrite the expression for $\mathrm{F}$: 
\begin{align*}
\mathrm{F} &: = \frac{f_{02}(t_1, t_2)}{2} y_{t_2}^2 + \frac{f_{30}(t_1, t_2)}{6} x_{t_2}^3 + \frac{f_{21}(t_1, t_2)}{2} x_{t_2}^2 y_{t_2}+ 
\frac{f_{12}(t_1, t_2)}{2} x_{t_2} y_{t_2}^2 + \frac{f_{03}(t_1, t_2)}{6} y_{t_2}^3 + \ldots.
\end{align*} 
Since $(\ff, \p, \lp) \in \PP \D_4$ we conclude that 
\begin{align}
f_{02} &=0, ~~f_{30} =0, ~~f_{21} \neq 0, ~~3f_{12}^2 - 4f_{21} f_{03} \neq 0. \label{pd4_nv_condition_equation}
\end{align}
To see why the two non vanishing conditions hold, first notice that since $(\ff, \lp) \in \hat{\D}_4$ we conclude that 
the cubic term in the Taylor expansion of $f$ has no repeated root. In other words 
\begin{align*}
\beta &:= f_{30}^2 f_{03}^2-6 f_{03} f_{12} f_{21} f_{30} + 4 f_{12}^3 f_{30} + 4 f_{03} f_{21}^3 - 3 f_{12}^2 f_{21}^2 \neq 0.   
\end{align*}
Since $(\ff, \lp) \in \PP \D_4$ we conclude $f_{30} =0$. Hence we get \eqref{pd4_nv_condition_equation}. 
Now we will show that \eqref{eval_f1_pd4} has no solutions. First of all we claim that 
$y_{t_2} \neq 0 $; we will justify that at the end. Assuming that, define 
$\mathrm{L} := \frac{x_{t_2}}{y_{t_2}}$. Substituting $\x_{t_2} = \mathrm{L} y_{t_2}$ in $\Fx =0 $ and using $y_{t_2} \neq 0$ and $f_{21}(t_1, t_2) \neq 0$ 
we can solve for $\mathrm{L}$ using the Implicit Function Theorem. That gives us 
\begin{align}
\mathrm{L} &= -\frac{f_{12}(t_1, t_2)}{2 f_{21}(t_1, t_2)} + y_{t_2} \mathrm{E}_1(y_{t_1}, f_{30}(t_1, t_2)) + f_{30}(t_1, t_2) \mathrm{E}_2(y_{t_1}, f_{30}(t_1, t_2)),  
\label{L_value_pd4_again}
\end{align}
where $\mathrm{E}_i(0,0) =0$. Using the value of $\mathrm{L}$ from \eqref{L_value_pd4_again}, 
and substituting $\x_{t_2} = \mathrm{L} y_{t_2}$ in $\mathrm{F} - \frac{y_{t_2}\Fy}{2} =0$, 
we conclude that as $(y_{t_2}, f_{30}(t_1, t_2))$ go to zero
\begin{align}
-\frac{f_{03}}{12} + \frac{f_{12}^2}{16 f_{21}} &=0. \label{pd5_dual_eqn}
\end{align}
It is easy to see that \eqref{pd5_dual_eqn} contradicts \eqref{pd4_nv_condition_equation}. 
It remains to show that $y_{t_2} \neq 0$. To see why that is so, consider the equation 
$\Fy =0$. It is easy to see that if $y_{t_2} =0$ then $f_{21}(t_1, t_2)$ will go to zero, contradicting 
\eqref{pd4_nv_condition_equation}. Hence \eqref{eval_f1_pd4} has no solutions. \\ 
\hf \hf  For the second part of the claim, we use the same set up except for one difference: we require $(\ff, \lp) \in \ov{\PP \A}_3$. 
Hence 
\[ f_{00}(t_1, t_2) = f_{10}(t_1, t_2) = f_{01}(t_1, t_2) =f_{20}(t_1, t_2) = f_{11}(t_1, t_2)= f_{30}(t_1, t_2)=0.\]
The functional equation \eqref{a1_pa3_intersect_pd5_empty_equation}  has a solution if and only if 
the following set of equations has a solution (as numbers): 
\begin{align}
\mathrm{F} = 0, \qquad \Fx = 0, \qquad \Fy = 0,  \qquad (x_{t_2}, y_{t_2}) \neq (0, 0) \qquad \textnormal{(but small)}. \label{eval_f1_pd5_again} 
\end{align} 
For the convenience of the reader, let us rewrite the expression for $\mathrm{F}$: 
\begin{align*}
\mathrm{F} &: = \frac{f_{02}(t_1, t_2)}{2} y_{t_2}^2 + \frac{f_{21}(t_1, t_2)}{2} x_{t_2}^2 y_{t_2}+ 
\frac{f_{12}(t_1, t_2)}{2} x_{t_2} y_{t_2}^2 + \frac{f_{03}(t_1, t_2)}{6} y_{t_2}^3 + \ldots.
\end{align*} 
Since $(\ff, \lp) \in \PP \D_5$ we conclude that 
\begin{align}
f_{02} &=0, ~~f_{30} =0, ~~f_{21} =0, ~~f_{40}\neq 0, ~~f_{12} \neq 0.  \label{pd5_nv_condition_equation_again}
\end{align}
We will now show that there are no solutions to \eqref{eval_f1_pd5_again}. First we claim that $y_{t_2} \neq 0$; we will justify that at the end. 
Assuming that, define $\mathrm{L} := \frac{x_{t_2}}{y_{t_2}}$. Substituting $\x_{t_2} = \mathrm{L} y_{t_2}$ in $\Fx =0 $ and using $y_{t_2} \neq 0$, we conclude that as $y_{t_2}$ and $f_{21}(t_1, t_2)$ 
go to zero, $f_{12}(t_1, t_2)$ goes to zero, contradicting \eqref{pd5_nv_condition_equation_again}. It remains to show that $y_{t_2} \neq 0$. Consider the equation $\Fx =0$. If $y_{t_2} =0$ then $f_{40}(t_1, t_2)$ would go to zero as $x_{t_2}$ goes to zero, 
contradicting \eqref{pd5_nv_condition_equation_again}. Hence \eqref{eval_f1_pd5_again} has no solutions. \qed \\

\hf \hf Now we return to the proof of \eqref{a1_pa2_d5}. First of all we observe that 
\eqref{a3_interesct_a1_pa2_is_empty_set_equation} and 
\eqref{a1_pa2_intersect_pd4_empty} imply that 
\begin{align}
\ov{\ov{\A}_1 \circ \PP \A}_2 \cap \Delta \hat{\D}_4 & = \varnothing.  \label{a1_pa2_intersect_hat_d4_is_empty}
\end{align}
Hence, the lhs of \eqref{a1_pa2_d5} is a subset of its 
rhs. This is because 
\begin{align}
 \ov{\hat{\D}}_4 & = \hat{\D}_4 \cup \ov{\hat{\D}}_5 \qquad \textnormal{and} \qquad \ov{\hat{\D}^{\#\flat}_5} = \ov{\hat{\D}}_5. \label{lot_of_closure_claims}
\end{align}
The first equality follows by applying Lemma \ref{tube_lemma} twice to Lemma \ref{cl} \eqref{D4cl} while the second is covered by Lemma \ref{Dk_sharp_closure}. To show 
rhs of \eqref{a1_pa2_d5} is a subset of its lhs, 
it suffices to show that   
\begin{align}
\{ (\ff, \p, \lp) \in \ov{\ov{\A}_1 \circ \PP \A_2}: ~\pi_2^*\us_{\PP \D_4}( \ff, \p, l_{\p} ) = 0\}
& \supset \Delta \hat{\D}_5^{\#\flat}.   \label{d5_hat+sharp_is_subset_of_a1_pa2}
\end{align}
\begin{claim}
\label{a1_pa2_closure_and_also_pa3}
Let $(\ff, \p, \lp) \in \Delta \hat{\D}_5^{\#\flat}$. Then there exists a solution 
$$ (\ff(t_1, t_2), \p(t_1, t_2),  l_{\p(t_1)} ) \in \ov{ (\D \times \P^2) \circ \PP \A}_2$$ 
\textit{near} $(\ff, \p, \lp)$ to the set of equations
\begin{align}
\pi_1^* \ds_{\A_0}(\ff(t_1, t_2), \p(t_1, t_2),  l_{\p(t_1)}) & = 0, ~\pi_1^* \ds_{\A_1}(\ff(t_1, t_2), \p(t_1, t_2),  l_{\p(t_1)}) = 0, ~~\p(t_1, t_2) \neq \p(t_1). \label{a1_pa2_closure_d5_sharp}
\end{align}
\end{claim}
\ni It is easy to see that claim \ref{a1_pa2_closure_and_also_pa3} 
implies \eqref{d5_hat+sharp_is_subset_of_a1_pa2}.\\ 

\pf Choose homogeneous coordinates $[\mathrm{X}: \mathrm{Y}: \mathrm{Z}]$ so that 
$\p = [0:0:1]$ and let ~$\mathcal{U}_{\p}$,   
$\pi_x$, $\pi_y$,  $v_1$, $w$, $v$, $\eta$, $\eta_{t_1}$, $\eta_{t_2}$, 
$x_{t_1}$, $y_{t_1}$, $x_{t_2}$, $y_{t_2}$, 
$f_{ij}(t_1, t_2)$, $\FF$, $\Fx$ and $\Fy$
be exactly the same as defined in the 
proof of claim \ref{claim_a4_closure_simultaneous}.
Since $(\ff(t_1, t_2), l_{\p(t_1)}) \in \ov{\PP \A}_2$, we conclude that 
\[ f_{00}(t_1, t_2) = f_{10}(t_1, t_2) = f_{01}(t_1, t_2) =f_{20}(t_1, t_2) = f_{11}(t_1, t_2)=0.\]
The functional equation \eqref{a1_pa2_closure_d5_sharp}  has a solution if and only if 
the following set of equations has a solution (as numbers): 
\begin{align}
\mathrm{F} = 0, \qquad \Fx = 0, \qquad \Fy = 0,  \qquad (x_{t_2}, y_{t_2}) \neq (0, 0) \qquad \textnormal{(but small)}. \label{eval_f1_d5_sharp} 
\end{align} 
For the convenience of the reader, let us rewrite the expression for $\mathrm{F}$: 
\begin{align*}
\mathrm{F} &: = {\textstyle \frac{f_{02}(t_1, t_2)}{2} y_{t_2}^2 + \frac{f_{30}(t_1, t_2)}{6} x_{t_2}^3 + \frac{f_{21}(t_1, t_2)}{2} x_{t_2}^2 y_{t_2}+ 
\frac{f_{12}(t_1, t_2)}{2} x_{t_2} y_{t_2}^2 + \frac{f_{03}(t_1, t_2)}{6} y_{t_2}^3 + \ldots}.
\end{align*}
Let us define  
\begin{align}
\beta_1 &:= f_{21}^2 - f_{12} f_{30},  
~~\beta_2^{\pm}:= -\frac{f_{03}}{12} -\frac{f_{21}^3}{6 f_{30}^2} + \frac{f_{12} f_{21}}{4 f_{30}} \pm
\sqrt{\beta_1} \Big( \frac{f_{21}^2}{6 f_{30}^2} - \frac{f_{12}}{6 f_{30}} \Big) \qquad \textnormal{and} \nonumber \\
\beta_3  & := f_{30}^2 f_{03}^2-6 f_{03} f_{12} f_{21} f_{30} + 4 f_{12}^3 f_{30} + 4 f_{03} f_{21}^3 - 3 f_{12}^2 f_{21}^2 = 4 f_{30}^2 \beta_{2}^{+} \beta_2^{-}. 
\label{beta_3_beta_2_plus_beta_2_minus}
\end{align}
Define $\beta_k(t_1, t_2)$ similarly with $f_{ij}$ replaced by $f_{ij}(t_1, t_2)$. 
Since $(\ff, \p, \lp) \in \hat{\D}_5$, the cubic 
\[ \Phi(\theta):= {\textstyle \frac{f_{30}}{6} \theta^3 + \frac{f_{21}}{2} \theta^2+ 
\frac{f_{12}}{2} \theta + \frac{f_{03}}{6}} \]
has a repeated root, 
but not all the three roots are the same. Hence we conclude  that 
\begin{align}
\beta_3 &=0, ~~\beta_1 \neq 0 ~~\textnormal{and} ~~f_{30} \neq 0. 
\end{align}
The last inequality follows from the fact that $(\ff, \p, \lp)$ belongs to  $\hat{\D}_5^{\#\flat}$ as opposed to
$\hat{\D}_5$.
We will now construct solutions to \eqref{eval_f1_d5_sharp}.    
Corresponding to each branch of $\sqrt{\beta_1(t_1, t_2)}$, the solutions are: 
\begin{align}
x_{t_2} &= \frac{-f_{21}(t_1, t_2) + \sqrt{\beta_1(t_1, t_2)}}{f_{30}(t_1, t_2)} y_{t_2} + \mathrm{O}(y_{t_2}^2), 
~~f_{02}(t_1, t_2) = \mathrm{O}(y_{t_2}), ~~\beta_2^{+}(t_1, t_2) = \mathrm{O}(y_{t_2}). \label{x_y_beta_soln_d5}  
\end{align}
Let us explain how we obtained these solutions. To obtain the value of $x_{t_2}$ we used $\Fx =0$. To obtain  the value of 
$f_{02}(t_1, t_2)$ we used $\Fy =0$ and the value of $x_{t_2}$ from the previous equation. Finally we used the fact that 
$2\FF - y_{t_2}\Fy =0$ and the value of $x_{t_2}$ to obtain $\beta_2^{+}(t_1, t_2)$. We get a similar solution 
for the other branch of $\sqrt{\beta_1(t_1, t_2)}$. 
By \eqref{beta_3_beta_2_plus_beta_2_minus} and \eqref{x_y_beta_soln_d5} , we conclude that 
as $y_{t_2}$ goes to zero, $f_{02}(t_1, t_2)$ and 
$\beta_3$ go to zero. Hence, the solutions in \eqref{x_y_beta_soln_d5} lie in 
$\ov{\A}_1\circ\ov{\PP \A}_2 $ and converge to a point 
$(\ff, \p , \lp)$  in  $\hat{\D}^{\#\flat}_5$.  \qed \\


\textbf{Proof of Lemma \ref{cl_two_pt} (\ref{a1_pa3_cl}):} By Lemma \ref{closure_a1_pak_f02_not_zero} ($k=3$), if we show that
\begin{align}
\{ (\ff, \p, l_{\p}) \in \ov{\ov{\A}_1 \circ \PP \A}_3: ~\pi_2^*\us_{\PP \D_4}( \ff, \p, l_{\p} ) = 0  \} &=  \Delta \ov{\PP \D_5^{\vee}} \cup 
\Delta \ov{\PP \D}_6  \label{a1_pa3_pd5_dual_eqn}
\end{align}
then we have
\bgd
\{ (\ff, \p, l_{\p}) \in \ov{\ov{\A}_1 \circ \PP \A}_3\}\subseteq \Delta\ov{\PP\A}_5 \cup \Delta\ov{\PP \D_5^{\vee}}\cup\Delta\ov{\PP \D}_6.
\edd
By Lemma \ref{cl} \eqref{A5cl}, we conclude that 
\bgd
\{ (\ff, \p, l_{\p}) \in \ov{\ov{\A}_1 \circ \PP \A}_3\}\subseteq \Delta\ov{\PP\A}_5 \cup \Delta\ov{\PP \D_5^{\vee}}.
\edd
On the other hand, by \eqref{pak2_is_subset_of_a1_and_pak} applied with $k=3$ and \eqref{a1_pa3_pd5_dual_eqn} we conclude that
\bgd
\{ (\ff, \p, l_{\p}) \in \ov{\ov{\A}_1 \circ \PP \A}_3\}\supseteq \Delta\ov{\PP\A}_5 \cup \Delta\ov{\PP \D_5^{\vee}}.
\edd
Thus, it suffices to prove \eqref{a1_pa3_pd5_dual_eqn}. Note that the intersection of the lhs of \eqref{a1_pa3_pd5_dual_eqn} with $\hat{\D}_4$ is empty. This follows from 
\eqref{a1_pa2_intersect_hat_d4_is_empty} and the fact that $\PP \A_3$ is a subset of $\ov{\PP \A}_2$ 
(see Lemma \ref{cl}). 
Equation \eqref{lot_of_closure_claims} 
and Lemma \ref{Dk_sharp_closure}, statement \ref{d5_pa3_zero}  now implies that the 
lhs of \eqref{a1_pa3_pd5_dual_eqn} is a subset of $\Delta \ov{\PP \D_5^{\vee}} \cup \Delta \ov{\PP \D}_5$.
However, the intersection of the lhs of \eqref{a1_pa3_pd5_dual_eqn} with $\Delta \PP \D_5$ is also empty by \eqref{a1_pa3_intersect_pd5_empty}. 
By Lemma \ref{cl}, statement \ref{D5cl} and Lemma \ref{Dk_sharp_closure}, statement \ref{pe6_subset_of_cl_pd5_dual} we have that 
\begin{align}
\ov{\PP \D}_5 & = \PP \D_5 \cup \ov{\PP \D}_6 \cup \ov{\PP \E}_6 \qquad \textnormal{and} \qquad \ov{\PP \E}_6  \subset \ov{\PP \D_5^{\vee}}.  \label{lot_of_closure_claims_again}
\end{align}
Hence the lhs of \eqref{a1_pa3_pd5_dual_eqn} is a subset of its rhs. 
To show the converse, we will first simultaneously prove the following three statements: 
\begin{align}
 \ov{\ov{\A}_1 \circ \PP \A}_3  & \supset \Delta \PP \D_5^{\vee},  \label{pd5_subset_a1_pa3}\\ 
 \ov{\ov{\A}_1 \circ \PP \D}_4 \cap \Delta \PP \D_5^{\vee} & = \varnothing,  \label{pd5_intersect_a1_pd4_is_empty}\\
\ov{\ov{\A}_1 \circ \PP \A}_4 \cap \Delta \PP \D_5^{\vee} & = \varnothing.  \label{pd5_intersect_a1_pa4_is_empty}
\end{align}
And then we will prove the following two statements simultaneously:
\begin{align}
\ov{\ov{\A}_1 \circ \PP \A}_3 & \supset \Delta \PP \D_6,   \label{a1_pa3_is_supsetof_pd6} \\  
\ov{\ov{\A}_1 \circ \PP \A}_4 \cap \Delta \PP \D_6 & = \varnothing. \label{a1_pa4_intersect_pd6_is_empty}  
\end{align}
Note that \eqref{pd5_subset_a1_pa3} and \eqref{a1_pa3_is_supsetof_pd6} imply rhs of 
\eqref{a1_pa3_pd5_dual_eqn} is a subset of its lhs, since $\ov{\ov{\A}_1 \circ \PP \A}_3$ is a closed set.  
\begin{claim}
\label{claim_pd5_subset_of_a1_pa3}
Let $(\ff,\p, \lp) \in \Delta \PP \D_5^{\vee}$.
Then there exists a solution 
$$ (\ff(t_1, t_2), \p(t_1, t_2),  l_{\p(t_1)} ) \in \ov{ (\D \times \P^2) \circ \PP \A}_3$$ 
\textit{near} $(\ff, \p, \lp)$ to the set of equations
\begin{align}
\pi_1^* \ds_{\A_0}(\ff(t_1, t_2), \p(t_1, t_2),  l_{\p(t_1)}) & = 0, ~\pi_1^* \ds_{\A_1}(\ff(t_1, t_2), \p(t_1, t_2),  l_{\p(t_1)}) = 0, ~~\p(t_1, t_2) \neq \p(t_1). \label{pd5_dual_limit_a1_pa3_functional_eqn}
\end{align}
Moreover, such solution  
lies in $\ov{\A}_1 \circ \PP \A_3$, i.e., 
\begin{align}
\pi_2^* \us_{\PP \D_4}(\ff(t_1, t_2), \p(t_1, t_2), l_{\p(t_1)}) & \neq 0, \label{psi_pd4_neq_0_a1_pa3_functional} \\
\pi_2^* \us_{\PP \A_4}(\ff(t_1, t_2), \p(t_1, t_2), l_{\p(t_1)}) & \neq 0. \label{psi_pa4_neq_0_a1_pa3_functional_new}
\end{align}
In particular, $(\ff(t_1, t_2), \p(t_1, t_2),  l_{\p(t_1)} )$ does not lie in $\ov{\ov{\A}_1 \circ \PP \D}_4$ or $\ov{\ov{\A}_1 \circ \PP \A}_4$. 
\end{claim}
\ni Note that claim \ref{claim_pd5_subset_of_a1_pa3} implies  \eqref{pd5_subset_a1_pa3}, \eqref{pd5_intersect_a1_pd4_is_empty} and 
\eqref{pd5_intersect_a1_pa4_is_empty} simultaneously. \\

\pf Choose homogeneous coordinates $[\mathrm{X}: \mathrm{Y}: \mathrm{Z}]$ so that 
$\p = [0:0:1]$ and let ~$\mathcal{U}_{\p}$,   
$\pi_x$, $\pi_y$,  $v_1$, $w$, $v$, $\eta$, $\eta_{t_1}$, $\eta_{t_2}$, 
$x_{t_1}$, $y_{t_1}$, $x_{t_2}$, $y_{t_2}$, 
$f_{ij}(t_1, t_2)$, $\FF$, $\Fx$ and $\Fy$
be exactly the same as defined in the 
proof of claim \ref{claim_a4_closure_simultaneous}, except for one difference:
we take $(\ff(t_1, t_2), l_{\p(t_1)})$ to be a point in $\ov{\PP \A}_3$.
Hence 
\[ f_{00}(t_1, t_2) = f_{10}(t_1, t_2) = f_{01}(t_1, t_2) =f_{11}(t_1, t_2)= f_{20}(t_1, t_2) = f_{30}(t_1, t_2)=0.\]
The functional equation \eqref{pd5_dual_limit_a1_pa3_functional_eqn}  
has a solution if and only if 
the following set of equations has a solution (as numbers): 
\begin{align}
\mathrm{F} = 0, \qquad \Fx = 0, \qquad \Fy = 0,  \qquad (x_{t_2}, y_{t_2}) \neq (0, 0) \qquad \textnormal{(but small)}. \label{eval_f1_d5_dual} 
\end{align} 
For the convenience of the reader, let us rewrite the expression for $\mathrm{F}$: 
\begin{align*}
\mathrm{F} &: = {\textstyle \frac{f_{02}(t_1, t_2)}{2} y_{t_2}^2 +\frac{f_{21}(t_1, t_2)}{2} x_{t_2}^2 y_{t_2}+ 
\frac{f_{12}(t_1, t_2)}{2} x_{t_2} y_{t_2}^2 + \frac{f_{03}(t_1, t_2)}{6} y_{t_2}^3 + \ldots}.
\end{align*}
Since $(\ff, \p, \lp) \in \Delta \PP \D_5^{\vee}$, we claim that 
\begin{align}
f_{21} & \neq 0, \qquad \beta_1 := \frac{(f_{12} \partial_x - 2 f_{21} \partial_y)^3 f}{2 f_{21}^2} = 3 f_{12}^2 -4f_{21} f_{03} =0  
\qquad \textnormal{and} \label{pd5_dual_nv1} \\
\beta_2&:= (f_{12} \partial_x - 2 f_{21} \partial_y)^4 f = 
f_{12}^4 f_{40} - 8f_{12}^3 f_{21} f_{31} + 24 f_{12}^2 f_{21}^2 f_{22}-32 f_{12} f_{21}^3 f_{13} + 16 f_{21}^4 f_{04} \neq 0. 
\label{pd5_dual_nv_2}
\end{align}
Let us justify this. Since $(\ff, \p) \in \D_5$ there exists a non zero vector $u = m_1v + m_2w$ such that 
\begin{align}
\nabla^3 f|_{\p} (u,u, v) & = m_1^2 f_{30} + 2 m_1 m_2 f_{21} + m_2^2 f_{12} =0,  \label{nabla_cube_f_1} \\
\nabla^3 f|_{\p} (u,u, w) &= m_1^2 f_{21} + 2 m_1 m_2 f_{12} + m_2^2 f_{03} =0, \label{nabla_cube_f_2} \\
\nabla^4 f|_{\p} (u,u,u,u) & \neq 0. \label{nabla_fourth_nv}
\end{align}
Since $(\ff, \lp) \in \PP \D_5^{\vee}$, we conclude \textit{by definition} that 
\begin{align}
f_{30} &=0, \qquad f_{21} \neq 0, \qquad m_2 \neq 0. \label{m_pd5_dual} 
\end{align}
(If $m_2 =0$ then $f_{21}$ would be zero). 
Equations \eqref{m_pd5_dual} and \eqref{nabla_cube_f_1} now imply that 
\begin{align}
m_1/ m_2 &= - f_{12}/(2 f_{21}). \label{m1_m2_value}
\end{align}
Equation \eqref{m1_m2_value} and \eqref{nabla_cube_f_2} implies \eqref{pd5_dual_nv1}. 
Finally, \eqref{nabla_fourth_nv} implies \eqref{pd5_dual_nv_2}.\\
\hf \hf We claim that solutions to \eqref{eval_f1_d5_dual} are given by  
\begin{align}
x_{t_2} = -\frac{f_{12}(t_1, t_2)}{2 f_{21}(t_1, t_2)} y_{t_2} + \mathrm{O}(y_{t_2}^2), \qquad 
\beta_1(t_1, t_2) &= -\frac{\beta_2(t_1, t_2)}{8 f_{21}(t_1, t_2)^3} y_{t_2} + \mathrm{O}(y_{t_2}^2) \qquad \textnormal{and} \nonumber \\
f_{02}(t_1, t_2) &= -\frac{\beta_2(t_1, t_2)}{192 f_{21}(t_1, t_2)^4} y_{t_2}^2 + \mathrm{O}(y_{t_2}^3). \label{d5_dual_multilicity}
\end{align}
Let us explain how we obtained these solutions. Assuming $\yt \neq 0$ (to be justified at the end) define $~\mathrm{L} := \frac{\xt}{\yt}.$
Using $x_{t_2} = \mathrm{L} y_{t_2}$ with $\yt\neq 0$ in the equation $\Fx=0$ we can solve for $\mathrm{L}$ via the Implicit Function Theorem. That gives us  
\begin{equation}
\mathrm{L} = \textstyle{\frac{-f_{12}(t_1, t_2)}{2 f_{21}(t_1, t_2)}}+ \Big(-\textstyle{\frac{f_{13}(t_1, t_2)}{6 f_{21}(t_1, t_2)}
+\frac{f_{12}(t_1, t_2) f_{22}(t_1, t_2)}{4 f_{21}(t_1, t_2)^2} - \frac{f_{12}(t_1, t_2)^2 f_{31}(t_1, t_2)}{8 f_{21}(t_1, t_2)^3} + 
\frac{f_{12}(t_1, t_2)^3 f_{40}(t_1, t_2)}{48 f_{21}(t_1, t_2)^4}} \Big) \yt + \mathrm{O}(\yt^2). \label{L_value_beta}
\end{equation}
Next, using the equation $2\FF - y_{t_2} \Fy =0$ and the fact that 
$\xt = \mathrm{L} \yt$ and \eqref{L_value_beta}, 
we obtain the expression for $\beta_1(t_1, t_2)$ in \eqref{d5_dual_multilicity}. 
Next, observe that 
\begin{align}
f_{30}(t_1, t_2) & = \frac{3 f_{12}(t_1, t_2)^2- \beta_1(t_1, t_2)}{4 f_{21}(t_1, t_2)}. \label{f30_beta1}
\end{align}
Finally, using the equation $\Fy=0$,  the fact that 
$\xt = \mathrm{L} \yt$, \eqref{L_value_beta}, the expression for $\beta_1(t_1, t_2)$ in \eqref{d5_dual_multilicity} and 
\eqref{f30_beta1} 
we obtain the expression for $f_{02}(t_1, t_2)$ in \eqref{d5_dual_multilicity}.  
It is now easy to see that \eqref{psi_pd4_neq_0_a1_pa3_functional} holds. It remains to show that $\yt \neq 0$. 
To see why that is so, suppose $\yt=0$. Then using the fact that $\Fy =0$, we conclude that 
$f_{21}(t_1, t_2) = \mathrm{O}(\xt)$. Hence $f_{21}(t_1, t_2)$ goes to zero as $\xt$ goes to zero, contradicting \eqref{pd5_dual_nv1}.\\ 
\hf \hf Finally, \eqref{psi_pa4_neq_0_a1_pa3_functional_new} is true because 
the section $\pi_2^*\us_{\PP \A_4}$ does not vanish on $\Delta \PP \D_5^{\vee}$. \qed 

\begin{cor}
\label{psi_pd4_section_vanishes_order_two_around_dual_d5}
Let $\W \lra \D \times \P^2 \times \P T\P^2$ be a vector bundle such that 
the rank of $\W$ is same as the dimension of $ \Delta \PP \D_5^{\vee}$ and 
$\Q: \D \times \P^2 \times \P T\P^2 \lra \W$  a \textit{generic} 
smooth section. Suppose $(\ff,\p, \lp) \in \Delta \PP \D_5^{\vee}\cap \Q^{-1}(0)$. 
Then the section $$ \pi_2^*\us_{\PP \D_4} \oplus \Q: \ov{\ov{\A}_1 \circ \PP\A}_3 \lra \pi_2^* (\UL_{\PP \D_4}) \oplus \W$$
vanishes around $(\ff, \p,  \lp)$ with a multiplicity of $2$.
\end{cor}
\pf This follows by observing that at $\Delta \PP \D_5^{\vee}$, the sections 
induced by $f_{02}$ and $\beta_1$  (the corresponding functionals) 
are transverse to the zero set over $\ov{\PP \A}_3$ 
\footnote{To see why; just take the partial derivative with respect to 
$f_{02}$ and $f_{03}$. Since $f_{21} \neq 0$, transversality follows.}, 
$\Q$ is generic and \eqref{d5_dual_multilicity}. \qed

\begin{claim}
\label{claim_a1_pa4_intersect_pd6_is_empty}
Let $(\ff,\p, \lp) \in \Delta \PP \D_6$.
Then there exist solutions 
$$ (\ff(t_1, t_2), \p(t_1, t_2),  l_{\p(t_1)} ) \in \ov{ (\D \times \P^2) \circ \PP \A}_3$$ 
\textit{near} $(\ff, \p, \lp)$ to the set of equations
\begin{align}
\pi_1^* \ds_{\A_0}(\ff(t_1, t_2), \p(t_1, t_2),  l_{\p(t_1)}) & = 0, ~\pi_1^* \ds_{\A_1}(\ff(t_1, t_2), \p(t_1, t_2),  l_{\p(t_1)}) = 0, \nonumber  \\ 
\pi_2^* \us_{\PP \D_4}(\ff(t_1, t_2), \p(t_1, t_2),  l_{\p(t_1)}) & \neq 0,  ~\p(t_1, t_2) \neq \p(t_1). \label{pd6_intersect_a1_pa4_is_empty_functional_eqn}
\end{align}
Moreover, any such solution sufficiently close to $(\ff, \p, \lp)$ lies in 
$\ov{\A}_1 \circ \PP \A_3$, i.e., 
\begin{align}
\pi_2^* \us_{\PP \A_4}(\ff(t_1, t_2), \p(t_1, t_2), l_{\p(t_1)}) \neq 0. \label{psi_pa4_neq_0_a1_pa3_functional}
\end{align}
In particular, $(\ff(t_1, t_2), \p(t_1, t_2),  l_{\p(t_1)} )$ does not lie in $\ov{\ov{\A}_1 \circ \PP \A}_4$. 
\end{claim}
\ni It is easy to see that claim \ref{claim_a1_pa4_intersect_pd6_is_empty} implies \eqref{a1_pa3_is_supsetof_pd6} and   
\eqref{a1_pa4_intersect_pd6_is_empty} simultaneously. \\

\pf Choose homogeneous coordinates $[\mathrm{X}: \mathrm{Y}: \mathrm{Z}]$ so that 
$\p = [0:0:1]$ and let ~$\mathcal{U}_{\p}$,   
$\pi_x$, $\pi_y$,  $v_1$, $w$, $v$, $\eta$, $\eta_{t_1}$, $\eta_{t_2}$, 
$x_{t_1}$, $y_{t_1}$, $x_{t_2}$, $y_{t_2}$, 
$f_{ij}(t_1, t_2)$, $\FF$, $\Fx$ and $\Fy$
be exactly the same as defined in the 
proof of claim \ref{claim_pd5_subset_of_a1_pa3}. 
Since $(\ff(t_1, t_2), l_{\p(t_1)} ) \in \Delta \ov{\PP \A}_3$ we conclude that 
\[ f_{00}(t_1, t_2) = f_{10}(t_1, t_2) = f_{01}(t_1, t_2) =f_{11}(t_1, t_2)= f_{20}(t_1, t_2) = f_{30}(t_1, t_2)=0.\]
Furthermore, since $(\ff, \p, \lp) \in \Delta \PP \D_6$, we conclude that 
\begin{align}
f_{21}, ~f_{40} &=0 ~~\textnormal{and} ~~f_{12}, ~\D^{f}_7 \neq 0. \label{pd6_vanish_non_vanish}
\end{align}
The functional equation \eqref{pd6_intersect_a1_pa4_is_empty_functional_eqn}has a solution if and only if the following has a numerical solution: 
\begin{align}
\mathrm{F} = 0, \qquad \Fx = 0, \qquad \Fy = 0,  \qquad (x_{t_2}, y_{t_2}) \neq (0, 0) \qquad \textnormal{(but small)}. \label{eval_f1_pd6} 
\end{align} 
For the convenience of the reader, let us rewrite the expression for $\mathrm{F}$: 
\begin{align*}
\mathrm{F} &: = {\textstyle \frac{f_{02}(t_1, t_2)}{2} y_{t_2}^2 + \frac{f_{21}(t_1, t_2)}{2} x_{t_2}^2 y_{t_2}+ \frac{f_{12}(t_1, t_2)}{2} x_{t_2} y_{t_2}^2 + \frac{f_{03}(t_1, t_2)}{6} y_{t_2}^3} + \ldots.
\end{align*}
We will now construct solutions to \eqref{eval_f1_pd6}. 
Let us define $\mathrm{G}:= \mathrm{F}-\frac{y_{t_2} \mathrm{F}_{y_{t_2}}}{2}-\frac{x_{t_2} \mathrm{F}_{x_{t_2}}}{4}$. Then
\begin{align*}
\mathrm{G}&:= -\frac{f_{03}(t_1, t_2)}{12} \yt^3 + \yt^3\mathrm{E}_1(\yt)\\
          & ~~~ -\xt \big({\textstyle \frac{f_{12}(t_1, t_2)}{8} \yt^2  +\frac{f_{31}(t_1, t_2)}{8} \xt^2 \yt +\frac{f_{50}(t_1, t_2)}{96} \xt^4  + \yt^2 \mathrm{E}_2(\xt, \yt)+\xt^2 \yt \mathrm{E}_3(\xt) + 
\xt^4 \mathrm{E}_4(\xt)} \big),
\end{align*}
where $\mathrm{E}_i$ is a holomorphic function vanishing at the origin. \\
\hf \hf We now make a change of variables by using a holomorphic function $\mathrm{C}(y_{t_2})$ such that by making the substitution $x_{t_2} = \hat{x}_{t_2} + \mathrm{C}(y_{t_2})$, the coefficients of $y_{t_2}^n$ are killed for all $n$ in $\mathrm{G}$. We may now make a change of coordinates $y_{t_2} = \hat{y}_{t_2} + \mathrm{B}(\hat{x}_{t_2})$ so that the coefficient of $\hat{x}_{t_2}\hat{y}_{t_2}^n$ is killed for all $n$ in $\mathrm{G}$. The existence of these functions follow from an identical argument as in 
in \cite{BM13}. After these changes, $\mathrm{G}$ is given by   
\begin{align*}
\mathrm{G} &=  {\textstyle -\frac{f_{12}(t_1, t_2)}{8} \hat{x}_{t_2} \Big( \hat{y}_{t_2}^2 + \frac{\D^{f(t_1, t_2)}_7}{60 f_{12}(t_1, t_2)} \hat{x}_{t_2}^4 + \mathrm{O}(\hat{x}_{t_2}^5)\Big)}. 
\end{align*}
Hence we are solving, in terms of the new variables $\hat{x}_{t_2}$ and $\hat{y}_{t_2}$, for 
\begin{align}
 \mathrm{G}&=0, ~~\mathrm{F}_{x_{t_2}} =0, ~~\mathrm{F}_{y_{t_2}} =0. \label{G_eq_0_Fx_eq_0_Fy_eq_0_D6}
\end{align}
Note that $\mathrm{F}_{x_{t_2}}$ and $\mathrm{F}_{y_{t_2}}$ are partials with respect to $x_{t_2}$ and $y_{t_2}$ 
expressed in the new coordinates, they are not the partials with respect to $\hat{x}_{t_2}$ and $\hat{y}_{t_2}$. 
Let us write $\mathrm{C}(\yt)$ and $\mathrm{B}(\hat{x}_{t_2})$ to second order: 
\begin{align*}
\mathrm{C}(\yt) = -& {\textstyle \frac{2 f_{03}(t_1, t_2)}{3 f_{12}(t_1, t_2)} \yt}\\ 
                + &  \big( {\textstyle -\frac{f_{04}(t_1, t_2)}{3 f_{12}(t_1, t_2)}  
+ \frac{2 f_{03}(t_1, t_2) f_{13} (t_1, t_2) }{3 f_{12}(t_1, t_2)^2 } - 
\frac{4 f_{03}(t_1, t_2)^2 f_{22} (t_1, t_2) }{9 f_{12}(t_1, t_2)^3 } + \frac{8 f_{03}(t_1, t_2)^3 f_{31} (t_1, t_2) }{81 f_{12}(t_1, t_2)^4}} \big) \yt^2\\ 
+ & \mathrm{O}(\yt^3), \\ 
\mathrm{B}(\hat{x}_{t_2}) = - &  {\textstyle \frac{f_{31}(t_1, t_2)}{6 f_{12}(t_1, t_2)} \hat{x}_{t_2}^2 + b_3 \hat{x}_{t_2}^3 + \mathrm{O}(\hat{x}_{t_2}^4)}.
\end{align*}
The exact value of $b_3$ is not important, 
but it will play a role in the 
subsequent calculation. 
The solutions to \eqref{G_eq_0_Fx_eq_0_Fy_eq_0_D6} are given by 
\begin{align}
\hat{y}_{t_2} &=  \alpha \hat{x}_{t_2}^2 +\mathrm{O}(\hat{x}_{t_2}^3) \label{pd6_hat_y}\\
f_{21}(t_1, t_2) &=\Big({\textstyle  - 2 \alpha f_{02}(t_1, t_2) + \frac{f_{02}(t_1, t_2) f_{31}(t_1, t_2)}{3 f_{12}(t_1, t_2)}} \Big) + \Big[{\textstyle -2 b_3 f_{02}(t_1, t_2) -\frac{8 \alpha^2 f_{02}(t_1, t_2) f_{03}(t_1, t_2)}{3 f_{12}(t_1, t_2)}} \nonumber\\
&  \,{\textstyle -2\alpha f_{12}(t_1, t_2)+\frac{8 \alpha f_{02}(t_1, t_2) f_{03}(t_1, t_2) f_{31}(t_1, t_2)}{9 f_{12}(t_1, t_2)^2} -\frac{2 f_{02}(t_1, t_2) f_{03}(t_1, t_2) f_{31}(t_1, t_2)^2}{27 f_{12}(t_1, t_2)^3}} \Big] \hat{x}_{t_2} \nonumber\\
& + \mathrm{O}(\hat{x}_{t_2}^2)  \label{pd6_f21}\\
f_{40}(t_1, t_2) &= \Big( 12 \alpha^2 f_{02}(t_1, t_2)- \frac{4 \alpha f_{02}(t_1, t_2) f_{31}(t_1, t_2)}{f_{12}(t_1, t_2)} + \frac{f_{02}(t_1, t_2) f_{31}(t_1, t_2)^2}{3 f_{12}(t_1, t_2)^2} \Big) \nonumber  \\
                 &+ \big({\textstyle 24 \alpha b_3 f_{02}(t_1, t_2)+\frac{32 \alpha^3 f_{02}(t_1, t_2) f_{03}(t_1, t_2)}{f_{12}(t_1, t_2)} + 24 \alpha^2 f_{12}(t_1, t_2)-4 \alpha f_{31}(t_1, t_2)} \nonumber \\
                 &{\textstyle -\frac{16 \alpha^2 f_{02}(t_1, t_2) f_{03}(t_1, t_2) f_{31}(t_1, t_2)}{f_{12}(t_1, t_2)^2}-\frac{4 b_3 f_{02}(t_1, t_2) f_{31}(t_1, t_2)}{f_{12}(t_1, t_2)}+ \frac{8 \alpha f_{02}(t_1, t_2) f_{03}(t_1, t_2) f_{31}(t_1, t_2)^2}{3 f_{12}(t_1, t_2)^3}}  \nonumber \\ 
                 & {\textstyle - \frac{4 f_{02}(t_1, t_2) f_{03}(t_1, t_2) f_{31}(t_1, t_2)^3}{27 f_{12}(t_1, t_2)^4}}\big) \hat{x}_{t_2} 
    + \mathrm{O}(\hat{x}_{t_2}^2), \label{pd6_f40}\\
\textnormal{where} ~~\alpha &:= \sqrt{-\frac{\D^{f(t_1, t_2)}_7}{60 f_{12}(t_1, t_2)}} ~~\textnormal{is a branch of the square root.}  \nonumber 
\end{align}
Note that each value of $\alpha$ corresponds to a different solution. Let us now explain how we obtained these solutions. 
First of all we claim that $\hat{x}_{t_2} \neq 0$; we will justify that at the end. Assuming that, we obtain 
\eqref{pd6_hat_y} from the fact $\mathrm{G}=0$. Next, we obtain \eqref{pd6_f21} from \eqref{pd6_hat_y} and using the fact that 
$\mathrm{F}_{y_{t_2}} =0$. Finally, we obtain \eqref{pd6_f40} from  \eqref{pd6_hat_y}, \eqref{pd6_f21} and using the fact that 
$\mathrm{F}_{x_{t_2}} =0$.   
Observe that equations \eqref{pd6_hat_y}, \eqref{pd6_f21} and \eqref{pd6_f40} imply that   
\begin{align}
f_{02}(t_1, t_2)\A^{f(t_1, t_2)}_4 &= \textstyle{\frac{1}{5}}\D^{f(t_1, t_2)}_7 f_{12}(t_1, t_2) \hat{x}_{t_2}^2 + \xt^2 \mathrm{E}(\xt, f_{02}(t_1, t_2)) \label{pd6_f02_pa4}
\end{align}
where $\mathrm{E}(0,0) =0$.
Hence, if $f_{02}(t_1, t_2)$ and $\hat{x}_{t_2}$ are small and non zero, $f_{02}(t_1, t_2)\A^{f(t_1, t_2)}_4$ is non zero. 
Hence \eqref{psi_pa4_neq_0_a1_pa3_functional} holds.  \\
\hf \hf It remains to show that $\hat{x}_{t_2} \neq 0$. If $\hat{x}_{t_2} =0$, then  
using the fact that $\mathrm{F}_{x_{t_2}} =0$ we get  
\begin{align*}
f_{12}(t_1, t_2) & = \frac{4 f_{03}(t_1, t_2) f_{21}(t_1, t_2)}{3 f_{12}(t_1, t_2)} + \mathrm{O}(\hat{y}_{t_2}).
\end{align*}
Hence $f_{12}(t_1, t_2)$ goes to zero as $f_{21}(t_1, t_2)$ and $\hat{y}_{t_2}$ go to zero, which  
contradicts \eqref{pd6_vanish_non_vanish}. \qed \\

\ni This proves Lemma \ref{cl_two_pt} (\ref{a1_pa3_cl}). Before proceeding further, note that  
\eqref{a1_pa3_intersect_pd5_empty} and \eqref{a1_pa4_intersect_pd6_is_empty} imply that 
\begin{align}
\Delta \mp & \subset \Delta \ov{\PP \D}_7. \label{one_a1_one_pa4_f02_zero_f12_not_zero_is_pd7}
\end{align}


\textbf{Proof of Lemma \ref{cl_two_pt} (\ref{a1_pa4_cl}):} By Lemma \ref{closure_a1_pdk_f12_not_zero} and \eqref{pak2_is_subset_of_a1_and_pak} for $k=4$, it suffices to show that 
\begin{align}
\{ (\ff, \p, l_{\p}) \in \ov{\ov{\A}_1 \circ \PP \A}_4: ~\pi_2^*\us_{\PP \D_4}( \ff, \p, l_{\p} ) = 0  \} &=  \Delta \ov{\mp} \cup \Delta \ov{\PP \E}_6 \label{a1_pa4_pd7s_pe6}
\end{align}
By the definition of $\Delta \mp$, to prove \eqref{a1_pa4_pd7s_pe6} it suffices to show that 
\begin{align}
\{ (\ff, \p, l_{\p})  \in \ov{\ov{\A}_1 \circ \PP \A}_4& : \pi_2^*\us_{\PP \D_4}( \ff, \p, l_{\p} ) = 0, ~\pi_2^*\us_{\PP \E_6}( \ff, \p, l_{\p} ) =0  \}  = \Delta \ov{\PP \E}_6. 
\label{one_a1_one_pa4_f02_zero_is_pe6}
\end{align}
It is clear that the lhs of \eqref{one_a1_one_pa4_f02_zero_is_pe6} is a subset of its rhs. 
To prove the converse, let us prove the following three facts simultaneously:
\begin{align}
\ov{\ov{\A}_1 \circ \PP \A}_4 & \supset \Delta \PP \E_6,    \label{a1_pa4_is_supsetof_pe6} \\  
\ov{\ov{\A}_1 \circ \PP \A}_5 \cap \Delta \PP \E_6 & = \varnothing,  \label{a1_pa5_intersect_pe6_is_empty} \\
\ov{\ov{\A}_1 \circ \PP \D}_4 \cap \Delta \PP \E_6 & = \varnothing.  \label{a1_pd4_intersect_pe6_is_empty} 
\end{align}
Note that since $\ov{\ov{\A}_1 \circ \PP \A}_4$ is a closed set,  
\eqref{a1_pa4_is_supsetof_pe6} implies  that the rhs of \eqref{one_a1_one_pa4_f02_zero_is_pe6}  is a subset of its lhs. 
We will need \eqref{a1_pd4_intersect_pe6_is_empty} later, since it follows from the present setup, we prove it here. 

\begin{claim}
\label{claim_a1_pa5_intersect_pe6_is_empty}
Let $(\ff,\p, \lp) \in \Delta \PP \E_6$.
Then there exist solutions 
$$ (\ff(t_1, t_2), \p(t_1, t_2),  l_{\p(t_1)} ) \in \ov{ (\D \times \P^2) \circ \PP \A}_3$$ 
\textit{near} $(\ff, \p, \lp)$ to the set of equations
\begin{align}
\pi_1^* \ds_{\A_0}(\ff(t_1, t_2), \p(t_1, t_2),  l_{\p(t_1)}) & = 0, ~\pi_1^* \ds_{\A_1}(\ff(t_1, t_2), \p(t_1, t_2),  l_{\p(t_1)}) = 0, \nonumber \\
\pi_2^* \us_{\PP \A_4}(\ff(t_1, t_2), \p(t_1, t_2),  l_{\p(t_1)}) & = 0, 
~\p(t_1, t_2) \neq \p(t_1). \label{pe6_intersect_a1_pa5_is_empty_functional_eqn}
\end{align}
Moreover, any such solution sufficiently close to $(\ff, \p, \lp)$ lies in 
$\ov{\A}_1 \circ \PP \A_4$, i.e., 
\begin{align}
\pi_2^* \us_{\PP \A_5}(\ff(t_1, t_2), \p(t_1, t_2), l_{\p(t_1)}) \neq 0,  \label{psi_pa5_neq_0_a1_pa4_functional} \\
\pi_2^* \us_{\PP \D_4}(\ff(t_1, t_2), \p(t_1, t_2), l_{\p(t_1)}) \neq 0. \label{psi_pd4_neq_0_a1_pa4_functional}
\end{align}
In particular, $(\ff(t_1, t_2), \p(t_1, t_2),  l_{\p(t_1)} )$ does not lie in $\ov{\ov{\A}_1 \circ \PP \A}_5$ or 
$\ov{\ov{\A}_1 \circ \PP \D}_4$. 
\end{claim}
\ni It is easy to see that claim \ref{claim_a1_pa5_intersect_pe6_is_empty} implies \eqref{a1_pa4_is_supsetof_pe6},  
\eqref{a1_pa5_intersect_pe6_is_empty} and \eqref{a1_pd4_intersect_pe6_is_empty} simultaneously. \\

\pf Choose homogeneous coordinates $[\mathrm{X}: \mathrm{Y}: \mathrm{Z}]$ so that 
$\p = [0:0:1]$ and let ~$\mathcal{U}_{\p}$,   
$\pi_x$, $\pi_y$,  $v_1$, $w$, $v$, $\eta$, $\eta_{t_1}$, $\eta_{t_2}$, 
$x_{t_1}$, $y_{t_1}$, $x_{t_2}$, $y_{t_2}$, 
$f_{ij}(t_1, t_2)$, $\FF$, $\Fx$ and $\Fy$
be exactly the same as defined in the 
proof of claim \ref{claim_pd5_subset_of_a1_pa3}.
Since $(\ff(t_1, t_2) , l_{\p(t_1)}) \in  \ov{\PP \A}_3$, we conclude that 
\[ f_{00}(t_1, t_2) = f_{10}(t_1, t_2) = f_{01}(t_1, t_2) =f_{11}(t_1, t_2)= f_{20}(t_1, t_2) = f_{30}(t_1, t_2)=0.\]
Moreover, since  $(\ff,\p, \lp) \in \Delta \PP \E_6$, we conclude that 
\begin{align}
f_{21}, f_{12} &=0, ~~f_{03}, f_{40} \neq 0. \label{pe6_condition_v_and_nv}
\end{align}
The functional equation \eqref{pe6_intersect_a1_pa5_is_empty_functional_eqn}  
has a solution if and only if 
the following set of equations has a solution (as numbers): 
\begin{align}
\mathrm{F} = 0, ~~\Fx = 0, ~~\Fy = 0, ~~f_{02}(t_1, t_2)\A^{f(t_1, t_2)}_4 =0, \quad (x_{t_2}, y_{t_2}) \neq (0, 0) 
~~\textnormal{(small)}. \label{eval_f1_pe6} 
\end{align}
Since $f_{02}(t_1, t_2)\A^{f(t_1, t_2)}_4 =0$ we conclude that 
$f_{02}(t_1, t_2) = \frac{3 f_{21}(t_1, t_2)^2}{f_{40}(t_1, t_2)}$. 
Hence 
\begin{align*}
\mathrm{F} & = \frac{3 f_{21}(t_1, t_2)^2}{2 f_{40}(t_1, t_2)} y_{t_2}^2 + \frac{f_{21}(t_1, t_2)}{2} x_{t_2}^2 y_{t_2}+ 
\frac{f_{12}(t_1, t_2)}{2} x_{t_2} y_{t_2}^2 + \frac{f_{03}(t_1, t_2)}{6} y_{t_2}^3 + \frac{f_{40}(t_1, t_2)}{24} x_{t_2}^4 \\ 
& + \frac{\mathcal{R}_{50}(x_{t_2})}{120} x_{t_2}^5+ 
\frac{\mathcal{R}_{31}(x_{t_2})}{6} x_{t_2}^3 y_{t_2} 
+\frac{\mathcal{R}_{22}(x_{t_2})}{4} x_{t_2}^2 y_{t_2}^2 
+\frac{\mathcal{R}_{13}(x_{t_2}, y_{t_2})}{4} x_{t_2} y_{t_2}^2 
+ \frac{\mathcal{R}_{04}(y_{t_2})}{24} y_{t_2}^4.  
\end{align*}
We will now eliminate $f_{12}(t_1, t_2)$ and $f_{21}(t_1, t_2)$ from \eqref{eval_f1_pe6}. 
First we  make a simple observation: let 
\begin{align*}
\mathrm{A}(\theta) &:= \mathrm{A}_0 + \mathrm{A}_1 \theta + \mathrm{A}_2 \theta^2, \quad \mathrm{B}(\theta) := \mathrm{B}_0 + \mathrm{B}_1 \theta, 
\quad p_1 := -\mathrm{A}_2 \mathrm{B}_1^2 , \quad p_2 := -\mathrm{A}_2^2 \mathrm{B}_0 + \mathrm{A}_1 \mathrm{A}_2 \mathrm{B}_1 +\mathrm{A}_2^2 \mathrm{B}_1 \theta.  
\end{align*}
Then $p_1 \mathrm{A}(\theta) + p_2\mathrm{B}(\theta)$ is independent of $\theta$. With this observation we will now proceed to define  
\begin{align*}
\mathrm{G}_1:=\mathrm{F} - x_{t_2} \mathrm{F}_{x_{t_2}}, ~~\mathrm{G}_2 := \mathrm{F} - \frac{y_{t_2} \mathrm{F}_{y_{t_2}}}{2}, ~~ \mathrm{G} := 
\mathrm{P}_1(t_1, t_2) \mathrm{G}_1 + \mathrm{P}_2(t_1, t_2) \mathrm{G}_2 
\end{align*}
where
\begin{align*}
\mathrm{P}_1 := - \frac{3 x^4 y^4}{32 f_{40}}, \qquad \mathrm{P}_2 &:= \frac{3 y^7 f_{03}}{16 f_{40}^2} + \frac{9 x^2 y^5 f_{21}}{16 f_{40}^2} - \frac{9 x^4 y^4}{32 f_{40}} + \frac{3 y^8 \mathcal{R}_{04}(y)}{32 f_{40}^2}
+\frac{3 x y^7 \mathcal{R}_{13}(x,y)}{16 f_{40}^2}  \\
             &-\frac{3 x^3 y^5 \mathcal{R}_{31}(x)}{16 f_{40}^2} - \frac{3 x^5 y^4 \mathcal{R}_{50}(x)}{160 f_{40}^2} 
                +\frac{3 y^9 \mathcal{R}_{04}^{(1)}(y)}{64 f_{40}^2}+\frac{3 x y^8 \mathcal{R}_{13}^{(0,1)}(x,y)}{16 f_{40}^2}.  
\end{align*}
The quantity $\mathrm{P}_k(t_1, t_2)$ is similarly defined, with $f_{ij}$, $x$ and $y$ 
replaced by $f_{ij}(t_1, t_2)$, 
$\xt$ and $\yt$.\\   
\hf \hf Note that $\mathrm{G}_1$ and $\mathrm{G}_2$ are independent of $f_{12}(t_1, t_2)$.  
Secondly, they are quadratic and linear in $f_{21}(t_1, t_2)$ respectively. 
Hence, using our previous observation, $\mathrm{G}$ is independent of $f_{12}(t_1, t_2)$ and $f_{21}(t_1, t_2)$. 
\footnote{Replace $\theta \lra f_{21}(t_1, t_2)$,  $\mathrm{A}(\theta) \lra \mathrm{G}_1$, $\mathrm{B}(\theta) \lra \mathrm{G}_2$, $p_1 \lra \mathrm{P}_1(t_1, t_2)$ and $p_2 \lra \mathrm{P}_2(t_1, t_2)$.}\\
\hf \hf 
We claim that $x_{t_2} \neq 0$ and $\yt \neq 0$; 
we will justify that at the end. 
Assuming this claim we conclude that  solving \eqref{eval_f1_pe6} 
is equivalent to solving: 
\begin{align}
\mathrm{G}=0, ~~\mathrm{G}_2 =0, ~~\mathrm{F} =0, 
\quad (x_{t_2}, y_{t_2}) \neq (0, 0) 
~~\textnormal{(small).} \label{eval_f1_pe6_modified}
\end{align}
Define $\mathrm{L}:= \xt^4 / \yt^3$. We will first solve for $\mathrm{L}$ in terms of $\xt$ and $\yt$ and then we will parametrize $(\xt, \yt)$.  
Notice that we can rewrite $\mathrm{G}$ in terms of $\xt$, $\yt$  and $\mathrm{L}$ in such a way that the highest power of 
$\xt$ is $3$; whenever there is an $\xt^4$ we replace it with $\mathrm{L} \yt^3$.
Hence  
\begin{align*}
\mathrm{G} &= \yt^{10}\Big(- \frac{f_{03}(t_1, t_2)^2}{64 f_{40}(t_1, t_2)^2} +\frac{f_{03}(t_1, t_2)}{64 f_{40}(t_1, t_2)} \mathrm{L}  + \mathrm{E}(\xt, \yt, \mathrm{L}) \Big)
\end{align*}
where $\mathrm{E}(0,0, \mathrm{L}) =0$. Hence $\mathrm{G} =0$ and $\yt \neq 0$ implies that 
\begin{align*}
 \Phi (\xt, \yt, \mathrm{L}) &:= - \frac{f_{03}(t_1, t_2)^2}{64 f_{40}(t_1, t_2)^2} +\frac{f_{03}(t_1, t_2)}{64 f_{40}(t_1, t_2)} \mathrm{L}  + \mathrm{E}(\xt, \yt, \mathrm{L}) =0.
\end{align*}
Hence, by the Implicit Function Theorem we conclude that 
\begin{align*}
\mathrm{L}(\xt, \yt) = \frac{f_{03}(t_1, t_2)}{f_{40}(t_1, t_2)} + \mathrm{E}_2(\xt, \yt)  
\end{align*}
where $\mathrm{E}_2(0, 0)$ is zero. Hence $(\xt, \yt)$ is parametrized by 
\begin{align*}
\yt &= u^4, \qquad \xt = \alpha u^3 + \mathrm{O}(u^4)  
\end{align*}
where $\alpha:= \sqrt[4]{\frac{f_{03}(t_1, t_2)}{f_{40}(t_1, t_2)}}$,  
a branch of the fourth root. Note that just one branch of the fourth root gives all the solutions, choosing another branch does not give us any more solutions.\footnote{Observe that a neighbourhood of the origin for the curve $\x^4 -y^3 =0$ is just one copy of $\C$ parametrized by $x = u^3$ and $y = u^4$.} We have
\begin{eqnarray}
f_{02}(t_1, t_2) & = & \frac{f_{03}(t_1, t_2)}{12} u^4 + \mathrm{O}(u^5)  \label{psi_pd4_neq_0_a1_pa4_number}\\
f_{02}(t_1, t_2)^2 \A^{f(t_1, t_2)}_5 & = &  -\frac{5 f_{03}(t_1, t_2)^2 f_{40}(t_1, t_2)}{18 \alpha} u^5 + \mathrm{O}(u^6) \label{pe6_multiplicity}
\end{eqnarray}
To arrive at these, use $\mathrm{G}_2 =0$ to get
\begin{align*}
f_{21}(t_1, t_2) &=  \frac{f_{03}(t_1, t_2)}{3 \alpha^2} u^2 + \mathrm{O}(u^3).
\end{align*}
We then use the fact that $f_{02}(t_1, t_2) = \frac{3 f_{21}(t_1, t_2)^2}{f_{40}(t_1, t_2)}$ to get \eqref{psi_pd4_neq_0_a1_pa4_number}. Finally using $\mathrm{F} =0$ we get that 
\begin{align*}
f_{12}(t_1, t_2) &= -\frac{2\alpha^3 f_{40}(t_1, t_2)}{3} u + \mathrm{O}(u^2)
\end{align*}
Plugging all this in we get \eqref{pe6_multiplicity}. Equations \eqref{psi_pd4_neq_0_a1_pa4_number} and \eqref{pe6_multiplicity} imply that \eqref{psi_pa5_neq_0_a1_pa4_functional} and \eqref{psi_pd4_neq_0_a1_pa4_functional} hold respectively. \\
\hf \hf It remains to show that that $x_{t_2} \neq 0$ and $\yt \neq 0$. 
If $y_{t_2} =0$ then using the fact that $\mathrm{F} =0$ we get  
that 
\begin{align*}
f_{40}(t_1, t_2) &= -\frac{\xt \mathcal{R}_{50}(\xt)}{5}.
\end{align*}
As $\xt$ goes to zero, $f_{40}(t_1, t_2)$ goes to zero, contradicting 
\eqref{pe6_condition_v_and_nv}. 
Similarly, if $\xt=0$, then using the fact that 
$\mathrm{F} - \frac{\yt \mathrm{F}_{\yt}}{2} =0$ we get that 
\begin{align*}
f_{03}(t_1, t_2) &= -\frac{2 \yt \mathcal{R}_{04}(\yt) + 
\yt^2 \mathcal{R}_{04}^{(1)}(\yt)}{4}
\end{align*}
As $\yt$ goes to zero, $f_{03}(t_1, t_2)$ goes to zero, contradicting 
\eqref{pe6_condition_v_and_nv}. 

\begin{cor}
\label{a1_pa4_mult_is_5_around_pe6}
Let $\W \lra \D \times \P^2 \times \P T\P^2$ be a vector bundle such that 
the rank of $\W$ is same as the dimension of $ \Delta \PP \E_{6}$ and 
$\Q: \D \times \P^2 \times \P T\P^2 \lra \W$  a \textit{generic} 
smooth section. Suppose $(\ff,\p, \lp) \in \Delta \PP \E_{6} \cap \Q^{-1}(0)$. 
Then the section $$ \pi_2^*\us_{\PP \A_{5}} \oplus \Q: \ov{\A_1 \circ \PP \A}_{4} \lra \pi_2^* (\UL_{\PP \A_{5}}) \oplus \W$$
vanishes around $(\ff, \p,  \lp)$ with a multiplicity of $5$.
\end{cor}
\pf Follows from the fact that the sections 
induced by $f_{02}$, $f_{21}$ and $f_{12}$ (the corresponding functionals) 
are transverse to the zero set over $\Delta \ov{\PP \A}_3$ \footnote{To see why, just take the partial derivatives with respect to 
$f_{02}$, $f_{21}$ and $f_{12}$}, the fact that  $\Q$ is generic and \eqref{pe6_multiplicity}. \qed \\

\ni This completes the proof of  Lemma \ref{cl_two_pt} (\ref{a1_pa4_cl}). \\


\textbf{Proof of Lemma \ref{cl_two_pt} (\ref{a1_pa5_cl}):} We have to show that 
\bgd
\{ (\ff, \p, l_{\p}) \in \ov{\ov{\A}_1 \circ \PP \A}_5\} =\Delta\ov{\PP\A}_7\cup\Delta\ov{\PP \D_8^s}\cup\Delta\ov{\PP \E}_7.
\edd
By Lemma \ref{closure_a1_pak_f02_not_zero} and \eqref{pak2_is_subset_of_a1_and_pak}, it is equivalent to showing that  
\begin{align}
\{ (\ff, \p, l_{\p}) \in \ov{\ov{\A}_1 \circ \PP \A}_5: ~\pi_2^*\us_{\PP \D_4}( \ff, \p, l_{\p} ) = 0  \} &=  \Delta \ov{\mq} \cup \Delta \ov{\PP \E}_7. \label{a1_pa5_pd8s_pe7}
\end{align}
By the definition of $\Delta \mq$, to prove \eqref{a1_pa5_pd8s_pe7} it suffices to show that 
\begin{align}
\{ (\ff, \p, l_{\p})  \in \ov{\ov{\A}_1 \circ \PP \A}_5 & : \pi_2^*\us_{\PP \D_4}( \ff, \p, l_{\p} ) = 0, ~\pi_2^*\us_{\PP \E_6}( \ff, \p, l_{\p} ) =0  \}  = \Delta \ov{\PP \E}_7. 
\label{one_a1_one_pa5_f02_zero_is_pe7}
\end{align}
Before we prove \eqref{one_a1_one_pa5_f02_zero_is_pe7}, we will first prove the following fact: 
\begin{align}
\Delta \mq & \subset \Delta \ov{\PP \D}_8. \label{one_a1_one_pa5_f02_zero_f12_not_zero_is_pd8}
\end{align}
Although \eqref{one_a1_one_pa5_f02_zero_f12_not_zero_is_pd8} 
is not required for the proof of Lemma \ref{cl_two_pt} (\ref{a1_pa5_cl}), it will be needed later. 
To prove \eqref{one_a1_one_pa5_f02_zero_f12_not_zero_is_pd8}, it suffices to show that 
\begin{align}
\ov{\ov{\A}_1 \circ \PP \A}_5 \cap \Delta \PP \D_7 & = \varnothing. \label{a1_pa5_intersect_pd7_is_empty}  
\end{align}
To see why, note that by \eqref{a1_pa4_intersect_pd6_is_empty} combined with $\PP \A_5\subset \ov{\PP \A}_4$ we have 
\bgd
\ov{\ov{\A}_1 \circ \PP \A}_5 \cap \Delta \PP \D_6=\varnothing.
\edd
Since $\Delta\ov{\PP \D}_8^s\subseteq \ov{\ov{\A}_1 \circ \PP \A}_5$,  
we conclude that \eqref{a1_pa5_intersect_pd7_is_empty} implies  
\bgd
\Delta\ov{\PP \D}_8^s \cap \big(\Delta \PP \D_6\cup \Delta \PP \D_7\big)=\varnothing.
\edd
On the other hand, by Lemma \ref{cl} (\ref{A5cl}) we 
know that $\Delta\ov{\PP \D}_8^s\subseteq \Delta \ov{\PP \D}_6$. 
Lemma \ref{Dk_sharp_closure} (\ref{pd6_pd7_pd8_closure}) now 
proves \eqref{one_a1_one_pa5_f02_zero_f12_not_zero_is_pd8} assuming \eqref{a1_pa5_intersect_pd7_is_empty}.
\begin{claim}
\label{claim_a1_pa5_intersect_pd7_is_empty}
Let $(\ff,\p, \lp) \in \Delta \PP \D_7$. 
Then there are no solutions
$$ (\ff(t_1, t_2), \p(t_1, t_2),  l_{\p(t_1)} ) \in \ov{ (\D \times \P^2) \circ \PP \A}_3$$ 
\textit{near} $(\ff, \p, \lp)$ to the set of equations
\begin{align}
\label{pd7_intersect_a1_pa5_is_empty_functional_eqn}
\pi_1^* \ds_{\A_0}(\ff(t_1, t_2), \p(t_1, t_2),  l_{\p(t_1)}) & = 0, ~\pi_1^* \ds_{\A_1}(\ff(t_1, t_2), \p(t_1, t_2),  l_{\p(t_1)}) = 0, \nonumber \\ 
\pi_2^* \us_{\PP \A_4}(\ff(t_1, t_2), \p(t_1, t_2),  l_{\p(t_1)}) & = 0, ~\pi_2^* \us_{\PP \A_5}(\ff(t_1, t_2), \p(t_1, t_2),  l_{\p(t_1)}) = 0,  \nonumber \\ 
\pi_1^* \us_{\PP \D_4}(\ff(t_1, t_2), \p(t_1, t_2),  l_{\p(t_1)}) & \neq 0, ~\p(t_1, t_2) \neq \p(t_1).
\end{align}
In particular, if $(\ff(t_1, t_2), \p(t_1, t_2),  l_{\p(t_1)} )$ is sufficiently close to 
$(\ff,\p, \lp)$, it does not lie in $\ov{\ov{\A}_1 \circ \PP \A}_5$. 
\end{claim}
\ni It is easy to see that claim \ref{claim_a1_pa5_intersect_pd7_is_empty} implies \eqref{a1_pa5_intersect_pd7_is_empty}. \\

\pf Choose homogeneous coordinates $[\mathrm{X}: \mathrm{Y}: \mathrm{Z}]$ so that 
$\p = [0:0:1]$ and let ~$\mathcal{U}_{\p}$,   
$\pi_x$, $\pi_y$,  $v_1$, $w$, $v$, $\eta$, $\eta_{t_1}$, $\eta_{t_2}$, 
$x_{t_1}$, $y_{t_1}$, $x_{t_2}$, $y_{t_2}$, 
$f_{ij}(t_1, t_2)$, $\FF$, $\Fx$ and $\Fy$
be exactly the same as defined in the 
proof of claim \ref{claim_pd5_subset_of_a1_pa3}. 
Since $(\ff(t_1, t_2), l_{\p(t_1)} ) \in \Delta \ov{\PP \A}_3$ we conclude that 
\[ f_{00}(t_1, t_2) = f_{10}(t_1, t_2) = f_{01}(t_1, t_2) =f_{11}(t_1, t_2)= f_{20}(t_1, t_2) = f_{30}(t_1, t_2)=0.\]
Furthermore, since $(\ff, \p, \lp) \in \Delta \PP \D_7$, we conclude that 
\begin{align}
f_{21}, ~f_{40}, ~\D_7^f &=0 ~~\textnormal{and} ~~f_{12}, ~\D^{f}_8 \neq 0. \label{pd7_vanish_non_vanish}
\end{align}
The functional equation \eqref{pd7_intersect_a1_pa5_is_empty_functional_eqn}  
has a solution if and only if 
the following set of equations has a solution (as numbers): 
\begin{align}
\mathrm{F} & = 0, \qquad \Fx = 0, \qquad \Fy = 0, \qquad \A^{f(t_1, t_2)}_4 =0, \qquad \A^{f(t_1, t_2)}_5 =0, \nonumber \\ 
(x_{t_2}, y_{t_2}) & \neq (0, 0), \qquad f_{02}(t_1, t_2) \neq 0  \qquad \textnormal{(but small)}. \label{eval_f1_pd7} 
\end{align} 
For the convenience of the reader, let us rewrite the expression for $\mathrm{F}$: 
\bge
\mathrm{F} : = {\textstyle \frac{f_{02}(t_1, t_2)}{2}} y_{t_2}^2 + {\textstyle \frac{f_{21}(t_1, t_2)}{2}} x_{t_2}^2 y_{t_2}+ {\textstyle \frac{f_{40}(t_1, t_2)}{24}}\xt^4 + {\textstyle \frac{f_{12}(t_1, t_2)}{2}} x_{t_2} y_{t_2}^2 +{\textstyle  \frac{f_{03}(t_1, t_2)}{6}} y_{t_2}^3 + \ldots. \label{F_recap_pd8}
\ede
We will now show that solutions to 
\eqref{eval_f1_pd7} can not exist.  
We will show that if $(f_{ij}(t_1, t_2), \xt, \yt)$ is a sequence converging 
to $(f,0,0) $
that satisfies \eqref{eval_f1_pd7}, then $\D^{f(t_1, t_2)}_8$ goes to zero (after passing to a subsequence).  
That would contradict \eqref{pd7_vanish_non_vanish}. 
First, we observe that any solution to the set of equations $\A^{f(t_1, t_2)}_4 =0$ and $\A^{f(t_1, t_2)}_5 =0$ 
is given by 
\begin{align}
f_{21}(t_1, t_2) & = {\textstyle \Big( \frac{f_{31}(t_1, t_2) \pm v }{ 3 f_{12}(t_1, t_2)} \Big)f_{02}(t_1, t_2), ~~f_{40}(t_1, t_2) = \frac{3 f_{21}(t_1, t_2)^2}{f_{02}(t_1, t_2)},~~\D^{f(t_1, t_2)}_7 = -\frac{5 v^2}{3 f_{12}(t_1, t_2)}}. \label{af4_af5_pd8}
\end{align}
Let us choose the $+v$ solution for $f_{21}(t_1, t_2)$, a similar argument will 
go through for the $-v$ solution. Now, using the value of $f_{40}(t_1, t_2)$ we observe that the first three terms in the expression for 
$\FF$ in \eqref{F_recap_pd8} form a perfect square. We will now rewrite the remaining part of $\FF$ by making a change of coordinates. 
Using an identical argument that is given in 
in \cite{BM13} and using \eqref{af4_af5_pd8}, we can make a change of coordinates 
\[ \x_{t_2} = \hat{x}_{t_2} + \mathrm{G}(\yt, v), \qquad \y_{t_2} = \hat{y}_{t_2} + \mathrm{B}(\hat{x}_{t_2}, v)  \]
so that $\mathrm{F}$ is given by 
\begin{align*}
\FF &= \frac{u}{2} \Big( \hat{y}_{t_2} + \mathrm{B}(\hat{x}_{t_2},v) + \frac{f_{31}(t_1, t_2)}{6 f_{12}(t_1, t_2)} (\hat{x}_{t_2} + \mathrm{J} )^2  + 
\frac{v (\hat{x}_{t_2} + \mathrm{J})}{6 f_{12}(t_1, t_2)}\Big)^2  \\ 
& -\frac{v^2 \hat{x}_{t_2}^5}{72 f_{12}(t_1, t_2)} + \frac{f_{12}(t_1, t_2)}{2} \hat{x}_{t_2} \hat{y}_{t_2}^2+ \frac{\D^{f(t_1, t_2)}_8}{720} \hat{x}_{t_2}^6 + 
\beta(\hat{x}_{t_2}) \hat{x}_{t_2}^7,
\end{align*}
where 
\[ \mathrm{J} := \mathrm{G}(\hat{y}_{t_2} + \mathrm{B}(\hat{x}_{t_2},v), v),\,\,\,u:=f_{02}(t_1,t_2).\]
Note that $\mathrm{B}$ is also a function of  $v$ because the coefficients of $\xt^n$ may depend   
on $f_{50}(t_1, t_2)$, which is equal to  $\D^{f(t_1, t_2)}_7 + \frac{5 f_{31}(t_1, t_2)^2}{3 f_{12}(t_1, t_2)} $. 
Let 
\[ z_{t_2} := \hat{y}_{t_2} + \frac{\hat{x}_{t_2}^2 v}{6 f_{12}(t_1, t_2)}.\]
Since 
\[\mathrm{G}(0,v) =0 \qquad \textnormal{and} \qquad  \mathrm{B}(\hat{x}_{t_2}, v) + \frac{f_{31}(t_1, t_2)}{6 f_{12}(t_1, t_2)} \hat{x}_{t_2}^2 = \mathrm{O}(\hat{x}_{t_2}^3) \]
we conclude that in the new coordinates $(\hat{x}_{t_2}, z_{t_2})$,  $\FF$ is given by  
\begin{align}
\FF &= \frac{u}{2} \Big( z_{t_2} + \alpha(\hat{x}_{t_2}, v) \hat{x}_{t_2} ^3 + \mathrm{E}(\hat{x}_{t_2}, z_{t_2}, v) z_{t_2} \Big)^2  \nonumber \\ 
& -\frac{v \hat{x}_{t_2}^3 z_{t_2}}{6 f_{12}(t_1, t_2)} + \frac{f_{12}(t_1, t_2)}{2} \hat{x}_{t_2} z_{t_2}^2+ \frac{\D^{f(t_1, t_2)}_8}{720} \hat{x}_{t_2}^6 + 
\beta(\hat{x}_{t_2}) \hat{x}_{t_2}^7, \qquad \textnormal{where} ~~\mathrm{E}(0,0,v) \equiv 0.\label{FF_new_form_pd8}
\end{align}
Hence, \eqref{eval_f1_pd7} has solutions if and only if the following set of equations has a solution 
\begin{align}
\FF &=0, \qquad \FF_{\hat{x}_{t_2}} =0, \qquad \FF_{z_{t_2}} =0, \qquad (\hat{x}_{t_2}, z_{t_2}) \neq (0,0), ~~u \neq 0, ~~v ~~\textnormal{small}. \label{f_eval_pd7_new}
\end{align}
We will now analyze solutions of \eqref{f_eval_pd7_new}. 
Notice that we can rewrite \eqref{f_eval_pd7_new} in the following way 
\begin{align}
p_0 + p_1 w + p_2 v &=0, ~~q_0 + q_1 w + q_2 v =0, ~~r_0 + r_1 w + r_2 v =0, \label{p0_q0_etc}
\end{align}
where 
\begin{align}
w &:= u (z_{t_2} + \alpha (\hat{x}_{t_2}, v) \hat{x}_{t_2}^3 + \mathrm{E}(\hat{x}_{t_2}, z_{t_2},  v)), 
~~\eta := z_{t_2} + \alpha (\hat{x}_{t_2}, v) \hat{x}_{t_2}^3 + \mathrm{E}(\hat{x}_{t_2}, z_{t_2},  v), \nonumber \\
p_0 &:=  \frac{f_{12}(t_1, t_2)}{2} \hat{x}_{t_2} z_{t_2}^2+ \frac{\D^{f(t_1, t_2)}_8}{720} \hat{x}_{t_2}^6 + 
\beta(\hat{x}_{t_2}) \hat{x}_{t_2}^7, ~~p_1 := \frac{\eta}{2}, ~~p_2 := -\frac{\hat{x}_{t_2}^3 z_{t_2}}{6}, \nonumber \\ 
q_0 &:= \frac{f_{12}(t_1, t_2)}{2}  z_{t_2}^2+ \frac{\D^{f(t_1, t_2)}_8}{120} \hat{x}_{t_2}^5 + 
(7 \beta(\hat{x}_{t_2}) + \hat{x}_{t_2} \beta^{\prime} (\hat{x}_{t_2})) \hat{x}_{t_2}^6, \nonumber \\ 
q_1 & := 3 \alpha (\hat{x}_{t_2}) \hat{x}_{t_2}^2 + \alpha^{\prime} (\hat{x}_{t_2}) \hat{x}_{t_2}^3 + \mathrm{E}_{{\hat{x}_{t_2}}}(\hat{x}_{t_2}, z_{t_2}, v), 
~~q_2:= -\frac{\hat{x}_{t_2}^2 z_{t_2}}{2}, \nonumber  \\ 
r_0& := f_{12}(t_1, t_2) \hat{x}_{t_2} z_{t_2}  ~~r_1 := 1+ \mathrm{E}(\hat{x}_{t_2}, z_{t_2}, v)+ z \mathrm{E}_{z_{t_2}}(\hat{x}_{t_2}, z_{t_2}, v)), 
~~r_2 := -\frac{\hat{x}_{t_2}^3}{6}. \label{p0_q0_defn}
\end{align}
Since \eqref{p0_q0_etc} holds, we conclude 
that 
\begin{align}
p_0 q_2 r_1 -p_0 q_1 r_2 + p_2 q_1 r_0 -p_1 q_2 r_0 -p_2q_0 r_1 + p_1 q_0 r_2 =0. \label{p0_q0_eliminated}  
\end{align}
Equations \eqref{p0_q0_eliminated} and \eqref{p0_q0_defn} now imply that  
\begin{align}
\Phi(\hat{x}_{t_2}, z_{t_2}) & : = z_{t_2}^3(f_{12}(t_1, t_2) + \mathrm{E}_2(\hat{x}_{t_2}, z_{t_2},  v)) 
+ z_{t_2}^2 \hat{x}_{t_2}^3(f_{12}(t_1, t_2) + \mathrm{E}_3(\hat{x}_{t_2}, z_{t_2}, v)) \nonumber \\  
 & ~~~~ +z_{t_2} \hat{x}_{t_2}^6(f_{12}(t_1, t_2) + \mathrm{E}_4(\hat{x}_{t_2}, z_{t_2}, v)) + 
 \hat{x}_{t_2}^9 \mathcal{R}_5(z_{t_2}, \hat{x}_{t_2}, v)=0,   \label{eliminate_u_v_pd8}
\end{align}
where $\mathrm{E}_i(0,0,v) =0$ and $\mathcal{R}_5(z_{t_2}, \hat{x}_{t_2}, v)$ is a holomorphic function. 
Now let us define $\mathrm{L} := z_{t_2}/ \hat{x}_{t_2}^3.$ Let $(f(t_1, t_2), \hat{x}_{t_2}, z_{t_2})$ 
be a sequence  converging to $(f, 0,0)$ 
that satisfies \eqref{f_eval_pd7_new} and such that $(\hat{x}_{t_2}, z_{t_2}) \neq (0,0)$.  
First we observe that $\hat{x}_{t_2} \neq 0$; this follows from 
\eqref{eliminate_u_v_pd8} and the fact that $f_{12} \neq 0$. 
Next, it is easy to see 
using \eqref{eliminate_u_v_pd8}, that
$\mathrm{L}$ is bounded, since $f_{12}\neq 0$. Hence, after passing to a subsequence $\mathrm{L}$ 
converges. Since $\FF =0$, we can easily see from \eqref{FF_new_form_pd8}, that 
as $(\hat{x}_{t_2}, z_{t_2}), u$ and $v$ go to zero, $\D^{f(t_1, t_2)}_8$ goes to zero. 
This contradicts \eqref{pd7_vanish_non_vanish}. \qed \\

\ni Now let us prove \eqref{one_a1_one_pa5_f02_zero_is_pe7}. It is clear that the lhs of \eqref{one_a1_one_pa5_f02_zero_is_pe7}  is a subset of its rhs. 
Next, we will prove the 
following two facts simultaneously:
\begin{align}
\ov{\ov{\A}_1 \circ \PP \A}_5 & \supset \Delta \PP \E_7,    \label{a1_pa5_is_supsetof_pe7} \\  
\ov{\ov{\A}_1 \circ \PP \A}_6 \cap \Delta \PP \E_7 & = \varnothing.  \label{a1_pa6_intersect_pe7_is_empty}  
\end{align}
Since $\ov{\ov{\A}_1 \circ \PP \A}_5$ is a closed set, \eqref{a1_pa5_is_supsetof_pe7} implies that the 
rhs of \eqref{one_a1_one_pa5_f02_zero_is_pe7} is a subset of its lhs. 

\begin{claim}
\label{claim_a1_pa6_intersect_pe7_is_empty}
Let $(\ff,\p, \lp) \in \Delta \PP \E_7$.
Then there exist solutions 
$$ (\ff(t_1, t_2), \p(t_1, t_2),  l_{\p(t_1)} ) \in \ov{ (\D \times \P^2) \circ \PP \A}_3$$ 
\textit{near} $(\ff, \p, \lp)$ to the set of equations
\begin{align}
\pi_1^* \ds_{\A_0}(\ff(t_1, t_2), \p(t_1, t_2),  l_{\p(t_1)}) & = 0, ~\pi_1^* \ds_{\A_1}(\ff(t_1, t_2), \p(t_1, t_2),  l_{\p(t_1)}) = 0, \nonumber \\
\pi_2^* \us_{\PP \A_4}(\ff(t_1, t_2), \p(t_1, t_2),  l_{\p(t_1)}) & = 0, ~~\pi_2^* \us_{\PP \A_5}(\ff(t_1, t_2), \p(t_1, t_2),  l_{\p(t_1)})  = 0, \nonumber \\
 ~\pi_2^* \us_{\PP \D_4}(\ff(t_1, t_2), \p(t_1, t_2),  l_{\p(t_1)}) & \neq 0, ~~\p(t_1, t_2) \neq \p(t_1). \label{pe7_intersect_a1_pa6_is_empty_functional_eqn} 
\end{align}
Moreover, any such solution sufficiently close to $(\ff, \p, \lp)$ lies in 
$\ov{\A}_1 \circ \PP \A_5$, i.e., 
\begin{align}
\pi_2^* \us_{\PP \A_6}(\ff(t_1, t_2), \p(t_1, t_2), l_{\p(t_1)}) \neq 0. \label{psi_pa6_neq_0_a1_pa5_functional}
\end{align}
In particular, $(\ff(t_1, t_2), \p(t_1, t_2),  l_{\p(t_1)} )$ does not lie in $\ov{\ov{\A}_1 \circ \PP \A}_6$. 
\end{claim}
\ni It is easy to see that claim \ref{claim_a1_pa6_intersect_pe7_is_empty} implies \eqref{a1_pa5_is_supsetof_pe7} and 
\eqref{a1_pa6_intersect_pe7_is_empty} simultaneously. \\ 
   
\pf Choose homogeneous coordinates $[\mathrm{X}: \mathrm{Y}: \mathrm{Z}]$ so that 
$\p = [0:0:1]$ and let ~$\mathcal{U}_{\p}$,   
$\pi_x$, $\pi_y$,  $v_1$, $w$, $v$, $\eta$, $\eta_{t_1}$, $\eta_{t_2}$, 
$x_{t_1}$, $y_{t_1}$, $x_{t_2}$, $y_{t_2}$, 
$f_{ij}(t_1, t_2)$, $\FF$, $\Fx$ and $\Fy$
be exactly the same as defined in the 
proof of claim \ref{claim_pd5_subset_of_a1_pa3}.
Since $(\ff(t_1, t_2), l_{\p(t_1)}) \in \ov{\PP \A}_3$, we conclude 
\[ f_{00}(t_1, t_2) = f_{10}(t_1, t_2) = f_{01}(t_1, t_2) =f_{11}(t_1, t_2)= f_{20}(t_1, t_2) = f_{30}(t_1, t_2)=0.\]
Moreover, since $(\ff, \p, \lp) \in \Delta \PP \E_7$, we conclude that
\begin{align}
f_{12}, f_{40} &=0, ~~f_{31}, f_{03} \neq 0. \label{pe7_nv}
\end{align}
The functional equation \eqref{pe6_intersect_a1_pa5_is_empty_functional_eqn}  
has a solution if and only if 
the following set of equations has a solution (as numbers): 
\begin{align}
\mathrm{F} & = 0, ~~\Fx = 0, ~~\Fy = 0, ~~\A^{f(t_1, t_2)}_4 =0, ~~\A^{f(t_1, t_2)}_5 =0, \nonumber \\ 
f_{02}(t_1, t_2) &\neq 0, \qquad (x_{t_2}, y_{t_2}) \neq (0, 0) \qquad \textnormal{(but small)}.  \label{eval_f1_pe7} 
\end{align}
Let us now define the following quantities: 
\begin{align*}
\mathrm{G}_1 &:= \mathrm{F} 
-\frac{\xt \mathrm{F}_{\xt}}{4} -\frac{\yt \mathrm{F}_{\yt}}{2}, ~~\mathrm{G}_2 := \mathrm{F}-\frac{\yt \mathrm{F}_{\yt}}{2}, 
~~\mathrm{G} := \mathrm{F}+ 4 \mathrm{G}_1 - 2 \mathrm{G}_2. 
\end{align*}
Note that $\mathrm{G}_1$ depends linearly on $f_{12}(t_1, t_2)$
and is independent of $f_{40}(t_1, t_2)$, 
$f_{21}(t_1, t_2)$ and $f_{02}(t_1, t_2)$; 
$\mathrm{G}_2$ depends linearly on $f_{40}(t_1, t_2)$ 
and $f_{21}(t_1, t_2)$
and is independent of $f_{12}(t_1, t_2)$, 
and $f_{02}(t_1, t_2)$; finally 
$\mathrm{G}$ depends linearly on $f_{02}(t_1, t_2)$ 
and $f_{40}(t_1, t_2)$
and is independent of $f_{12}(t_1, t_2)$, 
and $f_{21}(t_1, t_2)$. 
Next, we claim that $\xt \neq 0$ and $\yt \neq 0$; we will 
justify that at the end. Assuming this claim, we observe that 
\eqref{eval_f1_pe7} combined with \eqref{pe7_nv} 
is equivalent to
\begin{align}
\mathrm{G}_1 & = 0, ~~\mathrm{G}_2 = 0, ~~\mathrm{G} = 0, ~~f_{21}(t_1, t_2) = \frac{5f_{12}(t_1, t_2)f_{40}(t_1, t_2) + f_{02}(t_1, t_2) f_{50}(t_1, t_2)}{10 f_{31}(t_1, t_2)},  \nonumber \\
\A^{f(t_1, t_2)}_4 & =0,  ~~f_{02}(t_1, t_2) \neq 0, ~~ \xt \neq 0, ~~\yt \neq 0 \qquad \textnormal{(but small)}.  \label{eval_f1_pe7_modified} 
\end{align}
We will now construct solutions for  
\eqref{eval_f1_pe7_modified}. 
First of all, using $\mathrm{G} =0$ we can solve for $f_{02}(t_1, t_2)$ as a function of 
$f_{40}(t_1, t_2)$, $\xt$ and $\yt$. Next, 
using that $\mathrm{G}_1 =0$, we get $f_{12}(t_1, t_2)$ as a function of $\xt$ and $\yt$. Finally, using 
$\mathrm{G}_2 =0$, the value of $f_{02}(t_1, t_2)$, $f_{12}(t_1, t_2)$  
from the previous two equations and the value of $f_{21}(t_1, t_2)$ from \eqref{eval_f1_pe7_modified}, 
we get $f_{40}(t_1, t_2)$ in terms of $\xt$ and $\yt$. Plugging the expression back in, we get 
$f_{12}(t_1, t_2)$, $f_{21}(t_1, t_2)$, $f_{02}(t_1, t_2)$ and $f_{40}(t_1, t_2)$ in terms of $\xt$ 
and $\yt$.\\
\hf \hf Next, let us define $\mathrm{L} := \xt^3 /\yt^2.$ We note that any expression involving $\xt$ and $\yt$ can be re written in  terms of $\xt$, $\yt$ 
and $\mathrm{L}$ so that the highest power of $\xt$ is $2$; replace $\xt^3$ by $\mathrm{L} \yt^2$. 
Using that fact, 
we conclude    
\begin{align*}
\A^{f(t_1, t_2)}_4 \frac{\xt^4}{\yt^3}&= -\frac{f_{03}(t_1, t_2)^2}{3} + \frac{f_{03}(t_1, t_2) f_{31}(t_1, t_2)}{3} \mathrm{L} + \mathrm{E}_1(\xt, \yt, \mathrm{L}) =0 
\end{align*}
where $\mathrm{E}_1(0,0, \mathrm{L}) =0$. 
Hence, by the Implicit function theorem, we conclude that 
\begin{align*}
\mathrm{L} &= \frac{f_{03}(t_1, t_2)}{f_{31}(t_1, t_2)} + \mathrm{E}_2(\xt, \yt) 
\end{align*}
 where $\mathrm{E}_2(0,0)=0$. 
Hence, $\xt$ and $ \yt$ are parametrized by 
\[ \yt = u^3, ~~\xt = \alpha u^2 + \mathrm{O}(u^3) \qquad \textnormal{where} ~~\alpha := \sqrt[3]{\frac{f_{03}(t_1, t_2)}{f_{31}(t_1, t_2)}}, ~~\textnormal{a branch of the cube root.}  \]
Plugging in all this we get 
\begin{align}
f_{02}(t_1, t_2) &= \frac{f_{03}(t_1, t_2)}{3} u^3 + \mathrm{O}(u^4), 
~~f_{12}(t_1, t_2)= -\frac{f_{03}(t_1, t_2)}{3 \alpha } u + \mathrm{O}(u^2), ~~f_{21}(t_1, t_2) = \mathrm{O}(u^3), \nonumber \\ 
f_{40}(t_1, t_2) &= \mathrm{O}(u^2), ~~f_{02}(t_1, t_2)^3 \A^{f(t_1, t_2)}_6  =  -\frac{10 f_{03}(t_1, t_2)^2 f_{31}(t_1, t_2)^2 }{9 } u^6 + \mathrm{O}(u^7). \label{pe7_multiplicity}
\end{align}
Equation \eqref{pe7_multiplicity} implies that \eqref{psi_pa6_neq_0_a1_pa5_functional} holds. \qed \\

\begin{cor}
\label{a1_pa5_mult_is_5_around_pe7}
Let $\W \lra \D \times \P^2 \times \P T\P^2$ be a vector bundle such that 
the rank of $\W$ is same as the dimension of $ \Delta \PP \E_{7}$ and 
$\Q: \D \times \P^2 \times \P T\P^2 \lra \W$  a \textit{generic} 
smooth section. Suppose $(\ff,\p, \lp) \in \Delta \PP \E_{7} \cap \Q^{-1}(0)$. 
Then the section $$ \pi_2^*\us_{\PP \A_{6}} \oplus \Q: \ov{\ov{\A}_1 \circ \PP \A}_{5} \lra \pi_2^* (\UL_{\PP \A_{6}}) \oplus \W$$
vanishes around $(\ff, \p,  \lp)$ with a multiplicity of $6$.
\end{cor}

\pf Follows from the fact that the sections 
induced by $f_{02}$, $f_{21}$, $f_{12}$ and $f_{40}$ (the corresponding functionals) 
are transverse to the zero set over $\Delta\ov{\PP \A}_3$ \footnote{Take partial derivative with respect to 
$f_{02}$, $f_{21}$, $f_{12}$ and $f_{40}$.}, the fact that  $\Q$ is generic and \eqref{pe7_multiplicity}. \qed \\

\ni This finishes the proof of Lemma \ref{cl_two_pt} (\ref{a1_pa5_cl}).\\



\textbf{Proof of Lemma \ref{cl_two_pt} (\ref{a1_pd4_cl}):}  By definition of $\Delta \mr$ in \eqref{new_defn_delta}, 
it suffices to show  that 
\begin{align}
\{ (\ff, \p, l_{\p}) \in \ov{\ov{\A}_1 \circ \PP \D}_4: \pi_2^* \us_{\PP \D_5}(\ff, \lp) =0 \} &= \Delta \ov{\PP \D}_{6}. \label{a1_pd4_is_pd6}  
\end{align} 
Observe that 
\begin{align}
\ov{\ov{\A}_1 \circ \PP \D}_4 \cap \Delta (\PP \D_4 \cup \PP \D_5 \cup \PP \E_6) & = \varnothing. \label{a1_pd4_intersetc_pd4_is_empty}
\end{align}
This follows from \eqref{a1_pa2_intersect_pd4_empty}, \eqref{a1_pa3_pd5_dual_eqn} and \eqref{a1_pd4_intersect_pe6_is_empty} 
combined with 
the fact that 
\[ (\PP \D_4 \cup \PP \D_5) \cap (\ov{\PP \D_5^{\vee}} \cup \ov{\PP \D}_6) = \varnothing. \] 
Equation \eqref{a1_pd4_intersetc_pd4_is_empty} implies that the lhs of \eqref{a1_pd4_is_pd6} is a subset of its rhs. 
Next, we will simultaneously, prove the following two statements: 
\begin{align}
\ov{\ov{\A}_1 \circ \PP \D}_4 & \supset \Delta \PP \D_6, \label{a1_pa4_supset_pd6}\\ 
\ov{\ov{\A}_1 \circ \PP \D}_5 \cap \Delta \PP \D_6 & = \varnothing, \label{a1_pa5_intersetc_pd6_is_empty} 
\end{align}
Since $\ov{\ov{\A}_1 \circ \PP \D}_4$ is a closed set,  \eqref{a1_pa4_supset_pd6} implies that the rhs of \eqref{a1_pd4_is_pd6} is a subset of its lhs. 

\begin{claim}
\label{claim_pd5_actual_subset_of_a1_pa3}
Let $(\ff,\p, \lp) \in \Delta \PP \D_6$.
Then there exists a solution 
$$ (\ff(t_1, t_2), \p(t_1, t_2),  l_{\p(t_1)} ) \in \ov{ (\D \times \P^2) \circ \PP \D}_4$$ 
\textit{near} $(\ff, \p, \lp)$ to the set of equations
\begin{align}
\pi_1^* \ds_{\A_0}(\ff(t_1, t_2), \p(t_1, t_2),  l_{\p(t_1)}) & = 0, ~\pi_1^* \ds_{\A_1}(\ff(t_1, t_2), \p(t_1, t_2),  l_{\p(t_1)}) = 0, 
~\p(t_1, t_2) \neq \p(t_1). \label{pd5_actual_limit_a1_pa3_functional_eqn}
\end{align}
Moreover, such a solution  
lies in $\ov{\A}_1 \circ \PP \D_4$, i.e., 
\begin{align}
\pi_2^* \us_{\PP \D_5}(\ff(t_1, t_2), \p(t_1, t_2), l_{\p(t_1)}) & \neq 0. \label{psi_pa4_neq_0_a1_pa3_functional_new_pd5_actual}
\end{align}
In particular, $(\ff(t_1, t_2), \p(t_1, t_2),  l_{\p(t_1)} )$ does not lie in $\ov{\ov{\A}_1 \circ \PP \D}_5$. 
\end{claim}
\ni Note that claim \ref{claim_pd5_actual_subset_of_a1_pa3} implies  \eqref{a1_pa4_supset_pd6} and  
\eqref{a1_pa5_intersetc_pd6_is_empty} simultaneously. \\

\pf Choose homogeneous coordinates $[\mathrm{X}: \mathrm{Y}: \mathrm{Z}]$ so that 
$\p = [0:0:1]$ and let ~$\mathcal{U}_{\p}$,   
$\pi_x$, $\pi_y$,  $v_1$, $w$, $v$, $\eta$, $\eta_{t_1}$, $\eta_{t_2}$, 
$x_{t_1}$, $y_{t_1}$, $x_{t_2}$, $y_{t_2}$, 
$f_{ij}(t_1, t_2)$, $\FF$, $\Fx$ and $\Fy$
be exactly the same as defined in the 
proof of claim  
\ref{claim_a4_closure_simultaneous}, except for one difference:
we take $(\ff(t_1, t_2), l_{\p(t_1)})$ to be a point in $\ov{\PP \D}_4$.
Since $(\ff(t_1, t_2), l_{\p(t_1)}) \in \ov{\PP \D}_4$, 
we conclude that 
\[ f_{00}(t_1, t_2) = f_{10}(t_1, t_2) = f_{01}(t_1, t_2) =f_{11}(t_1, t_2)= f_{20}(t_1, t_2)=f_{02}(t_1, t_2) = f_{30}(t_1, t_2)=0.\]
The functional equation \eqref{pd5_dual_limit_a1_pa3_functional_eqn}  
has a solution if and only if 
the following set of equations has a solution (as numbers): 
\begin{align}
\mathrm{F} = 0, \qquad \Fx = 0, \qquad \Fy = 0, \qquad (x_{t_2}, y_{t_2}) \neq (0, 0) \qquad \textnormal{( but small)}. \label{eval_f1_d5_actual} 
\end{align} 
For the convenience of the reader, let us rewrite the expression for $\mathrm{F}$: 
\begin{align*}
\mathrm{F} &: = \frac{f_{21}(t_1, t_2)}{2} x_{t_2}^2 y_{t_2}+ 
\frac{f_{12}(t_1, t_2)}{2} x_{t_2} y_{t_2}^2 + \frac{f_{03}(t_1, t_2)}{6} y_{t_2}^3 + \ldots.
\end{align*}
Since $(\ff, \p, \lp) \in \Delta \PP \D_6$, we conclude that 
\begin{align}
f_{21}, ~f_{40} &= 0 ~~\textnormal{and} ~~f_{12}, ~\D^{f}_7 \neq 0. 
\end{align}
A little bit of thought will reveal 
that the solutions to 
\eqref{eval_f1_d5_actual}
are exactly the same as in \eqref{pd6_hat_y}, 
\eqref{pd6_f21} and \eqref{pd6_f40}, with $f_{02}(t_1, t_2) =0$. 
Since $\alpha \neq 0$, we conclude that $f_{21}(t_1, t_2) \neq 0$ 
for small but non zero $\hat{x}_{t_2}$. 
Hence \eqref{psi_pa4_neq_0_a1_pa3_functional_new_pd5_actual} holds. \qed 

\begin{cor}
\label{psi_pd5_section_vanishes_order_two_around_pd6}
Let $\W \lra \D \times \P^2 \times \P T\P^2$ be a vector bundle such that 
the rank of $\W$ is same as the dimension of $ \Delta \PP \D_6$ and 
$\Q: \D \times \P^2 \times \P T\P^2 \lra \W$  a \textit{generic} 
smooth section. Suppose $(\ff,\p, \lp) \in \Delta \PP \D_6\cap \Q^{-1}(0)$. 
Then the section $$ \pi_2^*\us_{\PP \D_4} \oplus \Q: \ov{\ov{\A}_1 \circ \PP\D}_4 \lra \pi_2^* (\UL_{\PP \D_4}) \oplus \W$$
vanishes around $(\ff, \p,  \lp)$ with a multiplicity of $2$.
\end{cor}
\pf Follows from the fact that the sections induced by $f_{21}$ and $f_{40}$ are transverse to the zero set over 
$\Delta \ov{\PP \D}_4$, \footnote{Take partial derivatives with respect to $f_{21}$ and $f_{40}$} 
the fact that $\Q$ is generic and \eqref{pd6_f21} (combined with $f_{02}(t_1, t_2) =0$). Each branch of 
$\alpha := \sqrt{\frac{\D^{f(t_1, t_2)}_7}{-60 f_{12}(t_1, t_2)}}$  
contributes with a multiplicity of $1$. Hence, the total multiplicity is $2$. \qed \\

\ni Before proceeding further, observe that \eqref{pd5_intersect_a1_pd4_is_empty} implies that 
\begin{align}
\Delta \ov{\mr} & \subset \Delta \ov{\PP \D_6^{\vee}}. \label{pd6_dual_s_is_subset_of_pd6_dual} 
\end{align}


\textbf{Proof of Lemma \ref{cl_two_pt} (\ref{a1_d4_cl}):}  Follows from Lemma \ref{cl_two_pt} (\ref{a1_pd4_cl}), \eqref{pd6_dual_s_is_subset_of_pd6_dual}, 
Lemma \ref{tube_lemma}, \eqref{tube_lemma_X} and 
\eqref{tube_lemma_Y}. \qed \\


\textbf{Proof of Lemma \ref{cl_two_pt} (\ref{a1_pd5_cl}):}  It suffices to prove the following two statements:  
\begin{align}
\{ (\ff, \p, l_{\p}) \in \ov{\ov{\A}_1\circ \PP \D}_5: ~\pi_2^*\us_{\PP \E_6}( \ff, \p, l_{\p} ) \neq 0  \} &=  
\{ (\ff, \p, l_{\p}) \in \Delta \ov{\PP \D}_{7}: ~\pi_2^*\us_{\PP \E_6}( \ff, \p, l_{\p} ) \neq 0  \}
\label{closure_a1_pd5_f12_not_zero}  \\
\{ (\ff, \p, l_{\p}) \in \ov{\ov{\A}_1 \circ \PP \D}_5: ~\pi_2^*\us_{\PP \E_6}( \ff, \p, l_{\p} ) = 0  \} &=  \Delta \ov{\PP \E}_7.   \label{a1_pd5_pe7}
 \end{align}
Let us directly prove a more general version of \eqref{closure_a1_pd5_f12_not_zero}: 
\begin{lmm}
\label{closure_a1_pdk_f12_not_zero}
If $k \geq 5$ then  
\begin{align*}
\{ (\ff, \p, l_{\p}) \in \ov{\ov{\A}_1 \circ  \PP \D}_k: ~\pi_2^*\us_{\PP \E_6}( \ff, \p, l_{\p} ) \neq 0  \} &= 
\{ (\ff, \p, l_{\p}) \in \Delta \ov{\PP \D}_{k+2}: ~\pi_2^*\us_{\PP \E_6}( \ff, \p, l_{\p} ) \neq 0  \}.  
\end{align*}
\end{lmm}
Note that \eqref{closure_a1_pd5_f12_not_zero} is a special case of Lemma \ref{closure_a1_pdk_f12_not_zero}; take $k=5$. 
We will prove the following two facts simultaneously: 
\begin{align}
\{ (\ff, \p, \lp) \in \ov{\ov{\A}_1 \circ \PP \D}_{k} \} \supset \Delta \PP \D_{k+2} \qquad & \forall ~k \geq 5, \label{pdk2_is_subset_of_a1_and_pdk} \\
\{ (\ff, \p, \lp) \in \ov{\ov{\A}_1 \circ \PP \D}_{k+1} \} \cap \Delta \PP \D_{k+2} = \varnothing \qquad  & \forall ~k \geq 4.  \label{pdk2_intersect_a1_and_pd1k+1_is_empty}
\end{align}
It is easy to see that \eqref{pdk2_is_subset_of_a1_and_pdk} and \eqref{pdk2_intersect_a1_and_pd1k+1_is_empty} imply 
Lemma \ref{closure_a1_pdk_f12_not_zero}. We will now prove the following claim:

\begin{claim}
\label{claim_d6_closure_simultaneous}
Let $~(\ff,\p, \lp) \in \Delta \PP \D_{k+2}$ and $ k\geq 5$.
Then there exists a solution 
$$ (\ff(t_1, t_2), \p(t_1, t_2),  l_{\p(t_1)} ) \in \ov{ (\D \times \P^2) \circ \PP \D}_5$$ 
\textit{near} $(\ff, \p, \lp)$ to the set of equations
\begin{align}
\label{closure_a1_dk_f12_not_zero}
\pi_1^* \us_{\A_0}(\ff(t_1, t_2), \p(t_1, t_2),  l_{\p(t_1)}) & = 0, ~\pi_1^* \us_{\A_1}(\ff(t_1, t_2), \p(t_1, t_2),  l_{\p(t_1)}) = 0, \nonumber \\
\pi_2^* \us_{\PP \D_6}(\ff(t_1, t_2), \p(t_1, t_2),  l_{\p(t_1)}) & = 0, \ldots, \pi_2^* \us_{\PP \D_k}(\ff(t_1, t_2), \p(t_1, t_2),  l_{\p(t_1)}) = 0, \nonumber \\
~\pi_2^* \us_{\PP \E_6}(\ff(t_1, t_2), \p(t_1, t_2),  l_{\p(t_1)}) & \neq  0, ~\p(t_1, t_2) \neq \p(t_1). 
\end{align}
Moreover, \textit{any} solution $(\ff(t_1, t_2), \p(t_1, t_2),  l_{\p(t_1)})$ sufficiently close to $(\ff, \p, \lp)$ 
lies in $ \ov{\A}_1 \circ \PP \D_k$, i.e.,
\begin{align}
\pi_2^* \us_{\PP \D_{k+1}}(\ff(t_1, t_2), \p(t_1, t_2),  l_{\p(t_1)}) \neq  0. \label{psi_pd_k_plus_1_does_not_vanish}
\end{align}
In particular $(\ff(t_1, t_2), \p(t_1, t_2),  l_{\p(t_1)})$ \textit{does not} lie in $ \ov{\ov{\A}_1 \circ \PP \D}_{k+1} $.
\end{claim}

\ni It is easy to see that claim \ref{claim_d6_closure_simultaneous} implies \eqref{pdk2_is_subset_of_a1_and_pdk} and \eqref{pdk2_intersect_a1_and_pd1k+1_is_empty} 
simultaneously for all $k\geq 5$.
The fact that \eqref{pdk2_intersect_a1_and_pd1k+1_is_empty} holds for $k=4$ is the content of \eqref{a1_pa5_intersetc_pd6_is_empty}. \\

\pf Choose homogeneous coordinates $[\mathrm{X}: \mathrm{Y}: \mathrm{Z}]$ so that 
$\p = [0:0:1]$ and let ~$\mathcal{U}_{\p}$,   
$\pi_x$, $\pi_y$,  $v_1$, $w$, $v$, $\eta$, $\eta_{t_1}$, $\eta_{t_2}$, 
$x_{t_1}$, $y_{t_1}$, $x_{t_2}$, $y_{t_2}$, 
$f_{ij}(t_1, t_2)$, $\FF$, $\Fx$ and $\Fy$
be exactly the same as defined in the 
proof of claim \ref{claim_pd5_subset_of_a1_pa3}.
Hence 
\[ f_{10}(t_1, t_2) = f_{01}(t_1, t_2) =f_{11}(t_1, t_2)= f_{20}(t_1, t_2)= f_{02}(t_1, t_2) = f_{30}(t_1, t_2)= f_{21}(t_1, t_2)=0.\]
Since $(\ff,\p, \lp) \in \Delta \PP \D_{k+2}$, we conclude that $f_{12}(t_1, t_2) \neq 0$. 
Hence, 
we can make a change of coordinates to write $\mathrm{F}$ as   
\begin{align*} 
\mathrm{F}&= \hat{y}_{t_2}^2 \hat{x}_{t_2} + \frac{\D^{f(t_1, t_2)}_6}{4!} \hat{x}_{t_2}^4 + \frac{\D^{f(t_1, t_2)}_7}{5!} \hat{x}_{t_2}^5+ \ldots 
\end{align*}
The functional equation   \eqref{closure_a1_dk_f12_not_zero}  has a solution if and only if 
the following set of equations has a solution (as numbers): 
\begin{align}
\label{closure_a1_dk_f12_not_zero_numbers}
\hat{y}_{t_2}^2 \hat{x}_{t_2} + \frac{\D^{f(t_1, t_2)}_6}{4!} \hat{x}_{t_2}^4 + \frac{\D^{f(t_1, t_2)}_7}{5!} \hat{x}_{t_2}^5+ \ldots &=0, 
\qquad \hat{y}_{t_2} \hat{x}_{t_2} = 0, \nonumber \\
\hat{y}_{t_2}^2 + \frac{\D^{f(t_1, t_2)}_6}{3!}\hat{x}_{t_2}^3 + \frac{\D^{f(t_1, t_2)}_7}{4!} \hat{x}_{t_2}^4 + \ldots &= 0, \qquad \D^{f(t_1, t_2)}_6, \ldots, \D^{f(t_1, t_2)}_k = 0, \nonumber  \\
(\hat{y}_{t_2}, \hat{x}_{t_2}) & \neq (0,0)   \qquad \textnormal{(but small).}
\end{align}
It is easy to see that the solutions to \eqref{closure_a1_dk_f12_not_zero_numbers} exist given by   
\begin{align}
\D^{f(t_1, t_2)}_6,& \ldots, \D^{f(t_1, t_2)}_k = 0, \nonumber \\
\D^{f(t_1, t_2)}_{k+1} &=  \frac{\D_{k+3}^{f(t_1, t_2)}}{k(k+1)} \hat{x}_{t_2}^2 + \mathrm{O}(\hat{x}_{t_2}^3),   \label{pdk+2_inside_a1_and_pdk_solution}\\ 
\D^{f(t_1, t_2)}_{k+2} &= -\frac{2 \D^{f(t_1, t_2)}_{k+3}}{(k+1)} \hat{x}_{t_2} + \mathrm{O}(\hat{x}_{t_2}^2), \qquad \hat{y}_{t_2} = 0, \qquad \hat{x}_{t_2}  \neq 0 \qquad \textnormal{(but small).} \nonumber 
\end{align}
By \eqref{pdk+2_inside_a1_and_pdk_solution}, it immediately follows that \eqref{psi_pd_k_plus_1_does_not_vanish} holds.  \qed  


\begin{cor}
\label{a1_pdk_mult_is_2_f12_neq_0}
Let $\W \lra \D \times \P^2 \times \P T\P^2$ be a vector bundle such that 
the rank of $\W$ is same as the dimension of $ \Delta \PP \D_{k+2}$ and 
$\Q: \D \times \P^2 \times \P T\P^2 \lra \W$  a \textit{generic} 
smooth section. Suppose $(\ff,\p, \lp) \in \Delta \PP \D_{k+2} \cap \Q^{-1}(0)$. 
Then the section $$ \pi_2^\ast\us_{\PP \D_{k+1}} \oplus \Q: \Delta \ov{\PP \D}_{k} \lra \pi_2^* (\UL_{\PP \D_{k+1}}) \oplus \W$$
vanishes around $(\ff, \p,  \lp)$ with a multiplicity of $2$.
\end{cor}
\pf  This follows from the fact that the sections 
\begin{align*}
\pi_2^\ast\us_{\PP \D_{i}}:& \Delta \ov{\PP \D}_{i-1} - \pi_2^\ast\us_{\PP \E_6}^{-1}(0) \lra  \pi_2^\ast\UL_{\PP \D_{i}}
\end{align*}
are transverse to the zero set for all $6 \leq i \leq k+2$, the fact that 
$\Q$ is generic and \eqref{pdk+2_inside_a1_and_pdk_solution}. \qed \\ 


\ni Next, we will prove \eqref{a1_pd5_pe7}. The lhs of \eqref{a1_pd5_pe7} is a subset of its rhs; this follows from 
\eqref{a1_pd4_intersect_pe6_is_empty} and the fact that $\PP \D_5$ is a subset of $\ov{\PP \D}_4$. To prove the converse, we will prove the following three statements simultaneously:
\begin{align}
\ov{\ov{\A}_1 \circ \PP \D}_5 & \supset \Delta \PP \E_7, \label{a1_pd5_supset_pe7}\\
\ov{\ov{\A}_1 \circ \PP \D}_6 \cap \Delta \PP \E_7 & = \varnothing. \label{a1_pd6_intersect_pe7_empty}\\
\ov{\ov{\A}_1 \circ \PP \E}_6 \cap \Delta \PP \E_7 & = \varnothing. \label{a1_pe6_intersect_pe7_empty}
\end{align}
Since  $\ov{\ov{\A}_1 \circ \PP \D}_5$ is a closed set, \eqref{a1_pd5_supset_pe7} implies that the rhs of 
\eqref{a1_pd5_pe7} is a subset of its lhs.

\begin{claim}
\label{claim_a1_pd6_intersect_pe7_is_empty}
Let $(\ff,\p, \lp) \in \Delta \PP \E_7$.
Then there exist solutions 
$$ (\ff(t_1, t_2), \p(t_1, t_2),  l_{\p(t_1)} ) \in \ov{ (\D \times \P^2) \circ \PP \D}_5$$ 
\textit{near} $(\ff, \p, \lp)$ to the set of equations
\begin{align}
\label{pe7_intersect_a1_pd6_is_empty_functional_eqn}
\pi_1^* \us_{\A_0}(\ff(t_1, t_2), \p(t_1, t_2),  l_{\p(t_1)}) & = 0, ~\pi_1^* \us_{\A_1}(\ff(t_1, t_2), \p(t_1, t_2),  l_{\p(t_1)}) = 0, ~~\p(t_1, t_2) \neq \p(t_1).
\end{align}
Moreover, any such solution sufficiently close to $(\ff, \p, \lp)$ lies in 
$\ov{\A}_1 \circ \PP \D_5$, i.e., 
\begin{align}
\pi_2^* \us_{\PP \D_6}(\ff(t_1, t_2), \p(t_1, t_2), l_{\p(t_1)}) \neq 0. 
\end{align}
In particular, $(\ff(t_1, t_2), \p(t_1, t_2),  l_{\p(t_1)} )$ does not lie in $\ov{\ov{\A}_1 \circ \PP \D}_6$. Since
\begin{align}
\pi_2^* \us_{\PP \E_6}(\ff(t_1, t_2), \p(t_1, t_2), l_{\p(t_1)}) \neq 0 
\end{align}
the solution $(\ff(t_1, t_2), \p(t_1, t_2),  l_{\p(t_1)} )$ does not lie in $\ov{\ov{\A}_1 \circ \PP \E}_6$. 
\end{claim}

\ni Note that claim \ref{claim_a1_pd6_intersect_pe7_is_empty} implies \eqref{a1_pd5_pe7} and \eqref{a1_pd6_intersect_pe7_empty} 
simultaneously. \\ 

\pf Choose homogeneous coordinates $[\mathrm{X}: \mathrm{Y}: \mathrm{Z}]$ so that 
$\p = [0:0:1]$ and let ~$\mathcal{U}_{\p}$,   
$\pi_x$, $\pi_y$,  $v_1$, $w$, $v$, $\eta$, $\eta_{t_1}$, $\eta_{t_2}$, 
$x_{t_1}$, $y_{t_1}$, $x_{t_2}$, $y_{t_2}$, 
$f_{ij}(t_1, t_2)$, $\FF$, $\Fx$ and $\Fy$
be exactly the same as defined in the 
proof of claim \ref{claim_pd5_subset_of_a1_pa3}
except for one difference:
we take $(\ff(t_1, t_2), l_{\p(t_1)})$ to be a point in $\ov{\PP \D}_5$.
Since $(\ff(t_1, t_2), l_{\p(t_1)}) \in \ov{\PP \D}_5$, we conclude that 
\[ f_{10}(t_1, t_2) = f_{01}(t_1, t_2) =f_{11}(t_1, t_2)= f_{20}(t_1, t_2)= f_{02}(t_1, t_2) = f_{30}(t_1, t_2)= f_{21}(t_1, t_2)=0.\]
Since $(\ff,\p, \lp) \in \Delta \PP \E_7$, we conclude that 
\begin{align}
f_{03}, f_{31} &\neq 0,  ~~f_{12}, f_{40} =0. \label{pe7_nv_again} 
\end{align} 
The functional equation \eqref{pe7_intersect_a1_pd6_is_empty_functional_eqn}  
has a solution if and only if 
the following set of equations has a solution (as numbers): 
\begin{align}
\mathrm{F} = 0, \qquad \Fx = 0, \qquad \Fy = 0,  \qquad (x_{t_2}, y_{t_2}) \neq (0, 0) \qquad \textnormal{(but small)}. \label{eval_f1_pe7_again} 
\end{align} 
Let us now define 
\begin{align*}
\mathrm{G} &:= -8\mathrm{F} + 2 \xt \mathrm{F}_{\xt} + 3 \yt \mathrm{F}_{\yt}. 
\end{align*}
We claim that $ \xt \neq 0$ and $ \yt \neq 0$; we will justify that at the end. Assuming that claim we 
conclude that solving \eqref{eval_f1_pe7_again} is equivalent to solving 
\begin{align}
\mathrm{G} = 0, \qquad \Fx = 0, \qquad \Fy = 0,  \qquad (x_{t_2}, y_{t_2}) \neq (0, 0) \qquad \textnormal{(but small)}. \label{eval_f1_pe7_again_modified} 
\end{align}
Note that $\mathrm{G}$ is independent of $f_{12}(t_1, t_2)$ and $f_{40}(t_1, t_2)$. Hence, $\mathrm{G}$ is explicitly given by 
\begin{equation*}
\mathrm{G}= \frac{\T_{03}(\xt,\yt)}{6}\yt^3+ \frac{\T_{31}(\xt,\yt)}{6}\xt^3y+\kq \xt^2\yt^2+\T_{50}(\xt)\xt^5,
\end{equation*}
where   $\T_{03}(0,0) = f_{03}(t_1, t_2)$ and $\T_{31}(0,0) = f_{31}(t_1, t_2)$.  
Using the same argument as in \cite{BM13}, 
there exists a holomorphic function $\Bq(\hat{x}_{t_2})$ and   
constant $\etq$  such that if we make the substitution 
$$\xt = \hat{x}_{t_2} + \etq \hat{y}_{t_2}, \qquad \yt = \hat{y}_{t_2} + \Bq(\hat{x}_{t_2}) \hat{x}_{t_2}^2 $$ 
then $\mathrm{G}$ is given by
\begin{align*}
\mathrm{G} &= \frac{\hat{\T}_{03}(\hat{x}_{t_2}, \hat{y}_{t_2})}{6} \hat{y}_{t_2}^3 + \frac{\hat{\T}_{31}(\hat{x}_{t_2}, \hat{y}_{t_2})}{6} \hat{x}_{t_2}^3 \hat{y}_{t_2},
\end{align*}
where $\hat{\T}_{03} (0,0) = f_{03}(t_1, t_2)$ and $\hat{\T}_{31} (0,0) = f_{31}(t_1, t_2)$. 
We claim that $\hat{y}_{t_2} \neq 0$; we will justify that at the end. Assuming that claim, we conclude from $\mathrm{G} =0$ that 
\begin{align*}
\hat{y}_{t_2} &= u^3, ~~\hat{x}_{t_2} = \alpha u^2 + \mathrm{O}(u^2) \qquad \textnormal{where} \qquad \alpha := \sqrt[3]{-\frac{f_{03}(t_1, t_2)}{f_{31}(t_1, t_2)}} \qquad \textnormal{a branch of the cuberoot.} 
\end{align*}
Using this and the remaining two equations of \eqref{eval_f1_pe7_again_modified} we conclude that
\begin{align}
f_{12}(t_1, t_2) & = - \frac{f_{03}(t_1, t_2)}{3 \alpha} u + \mathrm{O}(u^2), ~~f_{40}(t_1, t_2) = -\frac{4 f_{31}(t_1, t_2)}{\alpha} u + \mathrm{O}(u^2).  \label{f12_and_f40_mult_around_pe7}
\end{align}
It remains to show that $\xt \neq 0$, $\yt \neq 0$ and $\hat{y}_{t_2} \neq 0$. If $\xt =0$, then $\mathrm{F} =0$ implies that 
$f_{03}(t_1, t_2) = \mathrm{O}(\yt)$, 
contradicting \eqref{pe7_nv_again}. Next, if $\yt =0$ then $\mathrm{F}_{\yt} =0$ implies that $f_{31} (t_1, t_2) = \mathrm{O}(\xt)$, contradicting \eqref{pe7_nv_again}.
Finally, if $\hat{y}_{t_2} =0$, then $\mathrm{F}_{\yt} =0$ implies that 
\begin{align*}
f_{31}(t_1, t_2) &= -6\mathrm{B}(0) f_{12}(t_1, t_2) + \mathrm{O}(\xt).  
\end{align*}
As $f_{12}(t_1, t_2)$ and $\xt$ go to zero, $f_{31}(t_1, t_2)$ goes to zero, contradicting \eqref{pe7_nv_again}. \qed 

\begin{cor}
\label{psi_pe6_and_pd6_section_vanishes_order_one_around_pe7}
Let $\W \lra \D \times \P^2 \times \P T\P^2$ be a vector bundle such that 
the rank of $\W$ is same as the dimension of $ \Delta \PP \E_7$ and 
$\Q: \D \times \P^2 \times \P T\P^2 \lra \W$  a \textit{generic} 
smooth section. Suppose $(\ff,\p, \lp) \in \Delta \PP \E_7\cap \Q^{-1}(0)$. 
Then the sections 
\begin{align*}
\pi_2^*\us_{\PP \D_6} \oplus \Q: \ov{\ov{\A}_1 \circ \PP\D}_5 \lra \pi_2^* (\UL_{\PP \D_6}) \oplus \W, 
~~\pi_2^*\us_{\PP \E_6} \oplus \Q: \ov{\ov{\A}_1 \circ \PP\D}_5 \lra \pi_2^* (\UL_{\PP \E_6}) \oplus \W  
\end{align*}
vanish around $(\ff, \p,  \lp)$ with a multiplicity of $1$.
\end{cor}
\pf Follows from the fact 
that the sections induced by $f_{12}$ and $f_{40}$ are transverse to the zero set over 
$\Delta \ov{\PP \D}_5$, \footnote{Take partial derivatives with respect to $f_{12}$ and $f_{40}$.}
$\Q$ is generic and \eqref{f12_and_f40_mult_around_pe7}. \qed

\section{Euler class computation} 
\label{Euler_class_computation}

\hf\hf We are ready to prove the recursive formulas stated in section \ref{algorithm_for_numbers}.\\

\ni \textbf{Proof of \eqref{algoa1a1}:}  
Let $\Q: \D \times \P^2 \times \P^2 \lra \WL$ 
be a generic smooth section to
\bgd
\WL := \bigg({\textstyle \bigoplus}_{i=1}^{\delta_d -(n+2)} \pi_{\D}^*\gD^* 
\bigg)\oplus\Big({\textstyle \bigoplus}_{i=1}^{n}  \pi_2^* \gP^*\Big) \lra \D \times \P^2 \times \P^2.
\edd
Note that  
\begin{align*}
\Num(\A_1 \A_1, n) = \big\langle e(\WL), ~[\ov{\A_1\circ\A}_1] \big\rangle = | \pm (\A_1\circ \A_1) \cap  \Q^{-1}(0)|.
\end{align*}
By Lemma \ref{cl_two_pt} (\ref{a1a_minus_1_cl}) 
\begin{align*} 
\ov{\A}_1 \times \P^2 & = \ov{\ov{\A}_1 \circ (\D \times \P^2)} = \ov{\A}_1 \circ (\D \times \P^2) \du \Delta \ov{\A}_1.
\end{align*} 
The sections 
\begin{align}
 \pi_2^*\ds_{\A_0}:\ov{\A}_1 \times \P^2 - \Delta \ov{\A}_1   \lra \pi_2^* \DL_{\A_0}, \qquad \pi_2^*\ds_{\A_1}: \pi_2^*\ds_{\A_0}^{-1}(0) \lra \pi_2^*\DV_{\A_1}
\end{align}
are transverse to the zero set. 
(cf. Proposition \ref{ift_ml}). 
Hence 
\begin{align}
\big\langle e(\pi_2^*\DL_{\A_0}) e(\pi_2^*\DV_{\A_1}) e(\WL), ~[\ov{\A}_1 \times \P^2] \big\rangle & =  \Num(\A_1 \A_1, n) + 
\mathcal{C}_{\Delta \ov{\A}_1}(\pi_2^*\ds_{\A_0} \oplus \pi_2^*\ds_{\A_1} \oplus \Q), \label{A1A1_Euler_Class}
\end{align}
where $\mathcal{C}_{\Delta \ov{\A}_1}(\pi_2^*\ds_{\A_0} \oplus \pi_2^*\ds_{\A_1} \oplus \Q)$ 
is the contribution of the section 
$\pi_2^*\ds_{\A_0} \oplus \pi_2^*\ds_{\A_1} \oplus \Q$ to the Euler class from 
$\Delta \ov{\A}_1$. 
The lhs of \eqref{A1A1_Euler_Class}, as computed by splitting principle and a case by case check, is 
\begin{align}
\big\langle e(\pi_2^*\DL_{\A_0}) e(\pi_2^*\DV_{\A_1}) e(\WL), ~[\ov{\A}_1 \times \P^2] \big\rangle  &= \Num(\A_1, 0) \times \Num(\A_1, n) \label{a1a1_Euler_class_Main_stratum}.
\end{align}
In fact, the usual answer one arrives at is
\bgd
\Num(\A_1,n)+3(d-1)\Num(\A_1,n+1)+3(d-1)^2\Num(\A_1,n+2).
\edd
One then uses a result from \cite{BM13}:
\begin{align} 
\Num(\A_1,n) = \begin{cases}
3(d-1)^{2},&\textnormal{if}~n=0;\\
3(d-1),&\textnormal{if}~n=1;\\
1,&\textnormal{if}~n=2;\\
0,&\textnormal{otherwise}.
\end{cases}
\end{align}
\ni Next, we  compute 
$\mathcal{C}_{\Delta \ov{\A}_1}(\pi_2^*\ds_{\A_0}\oplus\pi_2^*\ds_{\A_1} \oplus \Q)$.
Note that ~$ \ov{\A}_1 = \A_1 \du \ov{\A}_2.$
By claim \ref{a1_section_contrib_from_a1_and_a2} we get that 
\begin{align}
\mathcal{C}_{\Delta \A_1}(\pi_2^*\ds_{\A_0}\oplus\pi_2^*\ds_{\A_1} \oplus \Q) & = \big\langle e(\pi_2^* \DL_{\A_0})e(\WL) , ~[\Delta \ov{\A}_1] \big\rangle  = \Num(\A_1, n) + d \Num(\A_1, n+1),  \label{a1a1_a1_contribution}\\
\mathcal{C}_{\Delta \ov{\A}_2}(\pi_2^*\ds_{\A_0}\oplus\pi_2^*\ds_{\A_1} \oplus \Q) &= 3 \Num(\A_2, n).  \label{a1a1_a2_contribution}
\end{align}
It is easy to see that \eqref{a1a1_Euler_class_Main_stratum}, \eqref{a1a1_a1_contribution} and \eqref{a1a1_a2_contribution} 
prove \eqref{algoa1a1}. \qed 
\begin{rem}
In \cite{Z1} a different method is used to compute $\mathcal{C}_{\Delta \A_1}(\pi_2^\ast\ds_{\A_0}\oplus\pi_2^*\ds_{\A_1} \oplus \Q)$. The author later on pointed out this simpler method to the second author of this paper. 
\end{rem}

\ni \textbf{Proof of \eqref{algopa20a1} and \eqref{algopa21a1}:} Let $\W_{n,m,2}^1$ and $\Q$ be as 
in \eqref{generic_Q} with $k=2$. 
By definition, $~\Num(\A_1\PP\A_2, n,m)$ is the signed 
cardinality of the intersection of $\A_1\circ \PP \A_2$ with $\Q^{-1}(0)$. 
By Lemma \ref{cl_two_pt}, statement \ref{a1a1_up_cl} we gather that
\begin{align*}
\ov{\ov{\A}_1 \circ  \hat{\A}^{\#}_1}& = \ov{\A}_1\circ \hat{\A}^{\#}_1 \du  \ov{\A}_1 \circ  (\ov{\hat{\A}_1^{\#}}- \hat{\A}_1^{\#}  ) \du \Delta \ov{\hat{\A}}_3 \\ 
                           & =  \ov{\A}_1\circ \hat{\A}^{\#}_1  \du \ov{\A}_1 \circ  \ov{\PP \A}_2 \du \Delta \ov{\hat{\A}}_3 \qquad 
\textnormal{(by \cite{BM13} (cf. Lemma \ref{cl}, statement \ref{A1cl}).}
\end{align*}
By Proposition \ref{A2_Condition_prp}, the section 
$$\pi_2^*\us_{\PP \A_2} : \ov{\ov{\A}_1 \circ  \hat{\A}^{\#}_1} \lra \pi_2^*\UV_{\PP \A_2}$$ 
vanishes on $\A_1 \circ \PP \A_2$ transversely. 
Hence, the zeros of the section 
$$ \pi_2^*\us_{\PP \A_2} \oplus \Q : \ov{\ov{\A}_1 \circ  \hat{\A}^{\#}_1} \lra \pi_2^*\UV_{\PP \A_2} \oplus \W_{n,m,2}^1, $$
restricted to $\A_1 \circ \PP \A_2$ counted with a sign, is our desired number. 
In other words 
\bgd
\Big\langle e(\pi_2^*\UV_{\PP \A_2}) e(\W_{n,m,2}^1), ~[\ov{\ov{\A}_1 \circ  \hat{\A}^{\#}_1}] \Big\rangle = \Num(\A_1 \PP\A_2, n,m) + \mathcal{C}_{\Delta \ov{\hat{\A}}_3}(\pi_2^*\us_{\PP \A_2} \oplus \Q) 
\edd
where $\mathcal{C}_{\Delta \ov{\hat{\A}}_3}(\pi_2^*\us_{\PP \A_2} \oplus \Q)$ is the contribution of the section to the Euler class from 
$\Delta \ov{\hat{\A}}_3$. Note that $ \pi_2^*\us_{\PP \A_2} \oplus \Q$ vanishes only on $\PP \A_3$ and $\hat{\D}_4$ and not on the entire $\ov{\hat{\A}}_3$.   
By Corollary \ref{pa2_section_mult_around_pa3} and  \ref{pa2_section_mult_around_hat_d4}, 
the contribution from $\PP \A_3$ and $\hat{\D}_4$ are  $2$ and $3$ respectively. This proves the claim.   
\begin{rem}
In the above proof we are using Lemma \ref{cohomology_ring_of_pv} with $M:= \ov{\ov{\A}_1 \circ  \hat{\A}^{\#}_1}$. However, in this case 
$M$ is not a smooth manifold; it is only a pseudocycle. Lemma \ref{cohomology_ring_of_pv} is actually true even when $M$ happens to be a 
pseudocycle. 
\end{rem}

\begin{rem}
A completely different method is used in \cite{Z1} to compute $\Num(\A_1 \A_2, n)$; instead of removing the cusp, the node is removed.
In fact, all the numbers $\Num(\A_1 \X_k,n)$ can also be computed by removing the node, instead of removing the $\X_k$ singularity. However, in order to obtain a recursive formula for the number of degree $d$ curves 
through $\delta_{d} - (\delta+k)$ generic points and having $\delta$-nodes and one singularity of type $\X_k$, 
we have to apply the method employed in this paper 
(i.e., we have to remove the $\X_k$-singularity, not the node). 
This observation is again due to Aleksey Zinger. 
\end{rem}


\ni \textbf{Proof of \eqref{algopa3a1}:} Let $\W_{n,m,3}^1$ and $\Q$ be as in \eqref{generic_Q} with $k=3$.
By Lemma \ref{cl_two_pt}, statement \ref{a1_pa2_cl} we have 
\begin{align*}
\ov{\ov{\A}_1 \circ  \PP \A}_2& = \ov{\A}_1 \circ  \PP \A_2 \du \ov{\A}_1 \circ (\ov{\PP \A}_2- \PP \A_2) \du 
\Big( \Delta \ov{\PP \A}_4 \cup \Delta \ov{\hat{\D}^{\#}_5}\Big), \\
                         &= \ov{\A}_1 \circ  \PP \A_2 \du \ov{\A}_1 \circ (\ov{\PP \A}_3 \cup \ov{\hat{\D}^{\#}_4}) \du 
\Big( \Delta \ov{\PP \A}_4 \cup \Delta \ov{\hat{\D}^{\#}_5}\Big) 
\end{align*}
where the last equality follows from \cite{BM13} (cf. Lemma \ref{cl}, statement \ref{A2cl}). By Proposition \ref{A3_Condition_prp}, the section 
$$ \pi_2^*\us_{\PP \A_3} \oplus \Q : \ov{\ov{\A}_1 \circ \PP \A}_2 \lra \pi_2^*\UL_{\PP \A_3} \oplus \W_{n,m,3}^1 $$
vanishes transversely on $\A_1 \circ \PP \A_3$. By definition, it does not vanish on $\A_1 \circ \hat{\D}_4^{\#}$. 
By Corollary \ref{a1_pak_mult_is_2_Hess_neq_0}, the contribution to the Euler class from the points of $\Delta \PP \A_4$ is $2$. 
Furthermore, by definition the section does not vanish on $\Delta\hat{\D}^{\#}_5$.
Hence 
\bgd
\Big\lan e(\pi_2^*\UL_{\PP \A_3} \oplus \W_{n,m,3}^1), ~~[\ov{\ov{\A}_1 \circ \PP \A}_2] \Big\ran = \Num(\A_1\PP\A_3,n,m) + 2 \Num(\PP \A_4, n, m)  
\edd
which proves the equation.  \qed \\


\ni \textbf{Proof of \eqref{algopa4a1}:} Let $\W_{n,m,4}^1$ and $\Q$ be as in \eqref{generic_Q} with $k=4$.
By Lemma \ref{cl_two_pt}, statement \ref{a1_pa2_cl} we have 
\begin{align*}
\ov{\ov{\A}_1 \circ  \PP \A}_3 &=  
\ov{\A}_1 \circ  \PP \A_3 \du \ov{\A}_1 \circ (\ov{\PP \A}_3- \PP \A_3) \du 
\Big( \Delta \ov{\PP \A}_5 \cup \Delta \ov{\PP \D^{\vee}_5} \Big), \\
&= \ov{\A}_1 \circ  \PP \A_3 \du \ov{\A}_1 \circ (\ov{\PP \A}_4 \cup \ov{\PP \D}_4) \du 
\Big( \Delta \ov{\PP \A}_5 \cup \Delta \ov{\PP \D^{\vee}_5} \Big) 
\end{align*}
where the last equality follows from \cite{BM13} (cf. Lemma \ref{cl}, statement \ref{A3cl}). By Proposition \ref{A3_Condition_prp}, the section 
$$ \pi_2^*\us_{\PP \A_4} \oplus \Q : \ov{\ov{\A}_1 \circ \PP \A}_3 \lra \pi_2^*\UL_{\PP \A_4} \oplus \W_{n,m,4}^1 $$
vanishes transversely on $\ov{\A}_1 \circ \PP \A_4$. It is easy to see that it does not vanish on $\A_1 \circ \PP \D_4$. 
By Corollary \ref{a1_pak_mult_is_2_Hess_neq_0}, the contribution to the Euler class from the points of $\Delta\PP \A_5$ is $2$. 
Moreover, the section does not vanish on $\Delta \PP \D^{\vee}_5$. 
Hence 
\bgd
\Big\langle e(\pi_2^*\UL_{\PP \A_4} \oplus \W_{n,m,4}^1), ~~[\ov{\ov{\A}_1 \circ \PP \A}_3] \Big\rangle = \Num(\A_1\PP\A_4,n,m) + 2 \Num(\PP \A_5, n, m)  
\edd
which proves the equation.  \qed \\ 

\ni \textbf{Proof of \eqref{algopa5a1}:} Let $\W_{n,m,5}^1$ and $\Q$ be as in \eqref{generic_Q} with $k=5$.
By Lemma \ref{cl_two_pt}, statement \ref{a1_pa2_cl} we have 
\begin{align*}
\ov{\ov{\A}_1 \circ  \PP \A}_4 &=  
\ov{\A}_1 \circ  \PP \A_4 \du \ov{\A}_1 \circ (\ov{\PP \A}_4- \PP \A_4) \du 
\Big( \Delta \ov{\PP \A}_6 \cup  
\Delta \ov{\mp} \cup \Delta \ov{\PP \E}_6  \Big), \\
&= \ov{\A}_1 \circ  \PP \A_4 \du \ov{\A}_1 \circ (\ov{\PP \A}_5 \cup \ov{\PP \D}_5) \du 
\Big( \Delta \ov{\PP \A}_6 \cup  
\Delta \ov{\mp} \cup \Delta \ov{\PP \E}_6  \Big)  
\end{align*}
where the last equality follows from \cite{BM13} (cf. Lemma \ref{cl}, statement \ref{A4cl}). By Proposition \ref{A3_Condition_prp}, the section 
$$ \pi_2^*\us_{\PP \A_5} \oplus \Q : \ov{\ov{\A}_1 \circ \PP \A}_4 \lra \pi_2^*\UL_{\PP \A_5} \oplus \W_{n,m,5}^1 $$
vanishes transversely on $\A_1 \circ \PP \A_5$. 
By \cite{BM13}
we conclude that this section  
vanishes on $\A_1 \circ \PP \D_5$ with a multiplicity of $2$. 
By Corollary  \ref{a1_pak_mult_is_2_Hess_neq_0} and \ref{a1_pa4_mult_is_5_around_pe6}, 
the contribution to the Euler class from the points of $\Delta\PP \A_6$ and $\Delta\PP \E_6$ are $2$ and $5$ respectively. 
Since the dimension of $\PP \D_7$ is one less than the rank of  $\pi_2^*\UL_{\PP \A_5} \oplus \W_{n,m,5}^1$ and $\Q$ is generic, 
the section 
does not vanish on $\Delta \ov{\PP \D}_7$. Since $\ov{\mp}$ is a subset of $\Delta \ov{\PP \D}_7$ 
(by \eqref{one_a1_one_pa4_f02_zero_f12_not_zero_is_pd7}),  
the section 
does not vanish on $\ov{\mp}$ either. 
Hence
\begin{align*}
\Big\langle e(\pi_2^*\UL_{\PP \A_5} \oplus \W_{n,m,5}^1), ~[\ov{\ov{\A}_1 \circ \PP \A}_4] \Big\rangle & = \Num(\A_1\PP\A_5,n,m) +2\Num(\A_1\PP\D_5,n,m)\\ 
                                                                                    & ~~ + 2 \Num(\PP \A_6, n, m)+ 5\Num(\PP \E_6, n ,m)  
\end{align*}
which proves the equation.  \qed \\ 

\ni \textbf{Proof of \eqref{algopa6a1}:} Let $\W_{n,m,6}^1$ and $\Q$ be as in \eqref{generic_Q} with $k=6$.
By Lemma \ref{cl_two_pt}, statement \ref{a1_pa5_cl} we have 
\begin{align*}
\ov{\ov{\A}_1 \circ  \PP \A}_5 & = \ov{\A}_1 \circ  \PP \A_5 \du \ov{\A}_1 \circ (\ov{\PP \A}_5- \PP \A_5) \du 
\Big( \Delta \ov{\PP \A}_7 \cup \Delta  \ov{\mq} \cup \Delta \ov{\PP \E}_7 \Big) \\
&=  \ov{\A}_1 \circ  \PP \A_5 \du \ov{\A}_1 \circ (\ov{\PP \A}_6 \cup \ov{\PP \D}_6 \cup \ov{\PP \E}_6 ) \du 
\Big( \Delta \ov{\PP \A}_7 \cup \Delta  \ov{\mq} \cup \Delta \ov{\PP \E}_7 \Big) 
\end{align*}
where the last equality follows from \cite{BM13} (cf. Lemma \ref{cl}, statement \ref{A5cl}). By Proposition \ref{A3_Condition_prp}, the section 
$$ \pi_2^*\us_{\PP \A_6} \oplus \Q : \ov{\ov{\A}_1 \circ \PP \A}_5 \lra \pi_2^*\UL_{\PP \A_6} \oplus \W_{n,m,6}^1 $$
vanishes transversely on $\A_1 \circ \PP \A_6$. 
By \cite{BM13}
we conclude that this section  
vanishes on $\A_1 \circ \PP \D_6$ with a multiplicity of $4$. 
By Corollary \ref{a1_pak_mult_is_2_Hess_neq_0} and \ref{a1_pa5_mult_is_5_around_pe7}, 
the contribution to the Euler class from the points of $\Delta\PP \A_7$ and $\Delta\PP \E_7$ are $2$ and $6$ respectively. 
Since the dimension of $\PP \D_8$ is one less than the rank of  $\pi_2^*\UL_{\PP \A_6} \oplus \W_{n,m,6}^1$ and $\Q$ is generic, 
the section 
does not vanish on $\Delta \ov{\PP \D}_8$. Since $\Delta \mq$ is a subset of $\Delta \ov{\PP \D}_8$ 
(by \eqref{one_a1_one_pa5_f02_zero_f12_not_zero_is_pd8}), the section does not vanish 
on $\Delta \ov{\mq}$ either. Hence 
\begin{align*}
\Big\langle e(\pi_2^*\UL_{\PP \A_6} \oplus \W_{n,m,6}^1), ~[\ov{\ov{\A}_1 \circ \PP \A}_5] \Big\rangle & = \Num(\A_1\PP\A_6,n,m) +4\Num(\A_1\PP\D_6,n,m)\\ 
                                                                                    & ~~ + 2 \Num(\PP \A_7, n, m)+ 6\Num(\PP \E_7, n ,m)  
\end{align*}
which proves the equation.  \qed \\ 

\ni \textbf{Proof of \eqref{algopd4a1}:} Let $\W_{n,0,4}^1$ and $\Q$ be as in \eqref{generic_Q} with $k=4$ and $m=0$.
By Lemma \ref{cl_two_pt}, statement \ref{a1_pa2_cl} we have 
\begin{align*}
\ov{\ov{\A}_1 \circ  \PP \A}_3 &= 
\ov{\A}_1 \circ  \PP \A_3 \du \ov{\A}_1 \circ (\ov{\PP \A}_3- \PP \A_3) \du 
\Big( \Delta \ov{\PP \A}_5 \cup \Delta \ov{\PP \D^{\vee}_5} \Big) \\
&= \ov{\A}_1 \circ  \PP \A_3 \du \ov{\A}_1 \circ (\ov{\PP \A}_4 \cup \ov{\PP \D}_4) \du 
\Big( \Delta \ov{\PP \A}_5 \cup \Delta \ov{\PP \D^{\vee}_5} \Big)
\end{align*}
where the last equality follows from \cite{BM13} (cf. Lemma \ref{cl}, statement \ref{A3cl}). By Proposition \ref{A3_Condition_prp}, the section 
$$ \pi_2^*\us_{\PP \D_4} \oplus \Q : \ov{\ov{\A}_1 \circ \PP \A}_3 \lra \pi_2^*\UL_{\PP \D_4} \oplus \W_{n,m,4}^1 $$
vanishes transversely on $\A_1 \circ \PP \D_4$. It is easy to see that the section does not vanish on $\A_1 \circ \PP \A_4$. 
By Corollary \ref{a1_pak_mult_is_2_Hess_neq_0} and \ref{psi_pd4_section_vanishes_order_two_around_dual_d5}, 
the contribution to the Euler class from the points of $\Delta\PP \A_5$ and $\PP \D_5^{\vee}$ are $2$ and $2$ respectively. 
Hence 
\begin{align}
\Big\langle e(\pi_2^*\UL_{\PP \D_4} \oplus \W_{n,0,4}^1), ~[\ov{\ov{\A}_1 \circ \PP \A}_3] \Big\rangle & = \Num(\A_1\PP\A_4,n,0) + 2 \Num(\PP \A_5, n, 0) \nonumber \\
&  +2\big\langle e(\W_{n,0,4}^1), ~[\ov{\PP \D^{\vee}_5}] \big\rangle. \label{pd5_dual_prelim_a1_pa4}
\end{align}
Since the map $\pi: \PP \D_5^{\vee} \lra \D_5$ is one to one, we conclude that 
\begin{align}
\big\langle e(\W_{n,0,4}^1), ~[\ov{\PP \D^{\vee}_5}] \big\rangle &= \Num(\D_5, n). \label{pd5_dual_d5_equality_numbers}
\end{align}
This follows from the same argument as in \cite{BM13}
Equations \eqref{pd5_dual_d5_equality_numbers} and \eqref{pd5_dual_prelim_a1_pa4} imply 
\eqref{algopd4a1}.\\  
\hf \hf Here is an alternative way to compute the lhs of \eqref{pd5_dual_d5_equality_numbers}. 
Recall that 
\[ \ov{\PP \D}_4 = \PP \D_4 \cup \ov{\PP \D}_5 \cup \ov{\PP \D^{\vee}_5}. \]
The section 
$\us_{\PP \D_5^{\vee}}: \ov{\PP \D}_4 \lra \UL_{\PP \D_5^{\vee}}$ vanishes transversely on $\PP \D_5^{\vee}$ and 
does not vanish on $\PP \D_5$. Hence 
\begin{align*}
\big\langle e(\UL_{\PP \D_5^{\vee}} \oplus \W_{n,0,5}^0), ~[\ov{\PP \D}_4] \big\rangle &= \big\langle e(\W_{n,0,5}^0), ~ [\ov{\PP \D^{\vee}_5}]\big\rangle =  \Num(\D_5, n).  
\end{align*}
It is easy to check directly that these two methods give the same answer for $\Num(\D_5, n)$. \qed \\

\ni \textbf{Proof of \eqref{algopd4a1_lambda} :} Let $\W_{n,1,4}^1$ and $\Q$ be as in \eqref{generic_Q} with $k=4$ and $m=1$.
By Lemma \ref{cl_two_pt}, statement \ref{a1_d4_cl} we have 
\begin{align}
\ov{\ov{\A}_1 \circ  \hat{\D}_4} & = \ov{\A}_1 \circ  \hat{\D}_4 \du \ov{\A}_1 \circ (\ov{\hat{\D}}_4- \hat{\D}_4) \du 
\Big( \Delta \ov{\hat{\D}}_6 \Big) \nonumber \\
\implies \ov{\ov{\A}_1 \circ  \hat{\D}^{\#}_4} &= \ov{\A}_1 \circ  \hat{\D}_4^{\#} \du \ov{\A}_1 \circ (\ov{\hat{\D}^{\#}_4}- \hat{\D}_4^{\#}) \du 
\Big( \Delta \ov{\hat{\D}^{\# \flat}_6} \Big) \qquad (\textnormal{since} ~~\ov{\hat{\D}}_4 = \ov{\hat{\D}^{\#}_4} ~~\textnormal{and} 
~~\ov{\hat{\D}}_6 = \ov{\hat{\D}^{\# \flat}_6}) \nonumber \\ 
&= \ov{\A}_1 \circ  \hat{\D}_4^{\#} \du \ov{\A}_1 \circ \ov{\PP \D}_4 \du 
\Big( \Delta \ov{\hat{\D}^{\# \flat}_6} \Big) \qquad \textnormal{(by definition).} \nonumber
\end{align}
By Proposition \ref{PD4_Condition_prp}, the section 
$$ \pi_2^*\us_{\PP \A_3} \oplus \Q : \ov{\ov{\A}_1 \circ  \hat{\D}^{\#}_4} \lra \pi_2^*\UL_{\PP \A_3} \oplus \W_{n,m,4}^1 $$
vanishes transversely on $\A_1 \circ \PP \D_4$. By definition, the section does not vanish on $\hat{\D}^{\# \flat }_6$. 
Hence 
\begin{align*}
\Big\langle e(\pi_2^*\UL_{\PP \A_3} \oplus \W_{n,1,4}^1), ~[\ov{\ov{\A}_1 \circ  \hat{\D}^{\#}_4}] \Big\rangle & = \Num(\A_1\PP\D_4,n,1)
\end{align*}
which proves the equation. \qed \\ 

\ni \textbf{Proof of \eqref{algopd5a1}:} Let $\W_{n,m,5}^1$ and $\Q$ be as in \eqref{generic_Q} with $k=5$. 
By Lemma \ref{cl_two_pt}, statement \ref{a1_pd4_cl} we have 
\begin{align*}
\ov{\ov{\A}_1 \circ  \PP \D}_4 & = \ov{\A}_1 \circ  \PP \D_4 \du \ov{\A}_1 \circ (\ov{\PP \D}_4- \PP \D_4) \du 
\Big( \Delta \ov{\PP \D}_6  \cup \Delta \ov{\mr} \Big)\\
&= \ov{\A}_1 \circ  \PP \D_4 \du \ov{\A}_1 \circ (\ov{\PP \D}_5 \cup \ov{\PP \D^{\vee}_5})  \du 
\Big( \Delta \ov{\PP \D}_6 \cup \Delta \ov{\mr} \Big)
\end{align*}
where the last equality follows from \cite{BM13} (cf. Lemma \ref{cl}, statement \ref{D4cl}). By Proposition \ref{D4_Condition_prp}, the section 
$$ \pi_2^*\us_{\PP \D_5}^{\mathbb{L}} \oplus \Q : \ov{\ov{\A}_1 \circ \PP \D}_4 \lra \pi_2^*\UL_{\PP \D_5} \oplus \W_{n,m,5}^1 $$
vanishes transversely on $\A_1 \circ \PP \D_5$. Moreover, it does not vanish on $\ov{\A}_1 \circ \PP \D_5^{\vee}$ by definition. 
By Corollary \ref{psi_pd5_section_vanishes_order_two_around_pd6}, 
the contribution to the Euler class from the points of $\Delta\PP \D_6$ is $2$.  
The section does not vanish on $\Delta \mr$ by definition. 
Since the dimension of $\PP \D_6^{\vee}$ is same as the dimension of $\PP \D_6$ 
and $\Q$ is generic, by \eqref{pd6_dual_s_is_subset_of_pd6_dual}, the section does not vanish on $\Delta \ov{\mr}$.  
Hence 
\begin{align*}
\Big\langle e(\pi_2^*\UL_{\PP \D_5} \oplus \W_{n,m,5}^1), ~[\ov{\ov{\A}_1 \circ \PP \D}_4] \Big\rangle & = \Num(\A_1\PP\D_5,n,m) + 2 \Num(\PP \D_6, n, m)
\end{align*}
which proves the equation. \qed \\

\ni \textbf{Proof of \eqref{algopd6a1} and \eqref{algope6a1}:} Let $\W_{n,m,6}^1$ and $\Q$ be as in \eqref{generic_Q} with $k=6$. 
By Lemma \ref{cl_two_pt}, statement \ref{a1_pd5_cl} we have 
\begin{align*}
\ov{\ov{\A}_1 \circ  \PP \D}_5 & = \ov{\A}_1 \circ  \PP \D_5 \du \ov{\A}_1 \circ (\ov{\PP \D}_5- \PP \D_5) \du 
\Big( \Delta \ov{\PP \D}_7 \cup \Delta \ov{\PP \E}_7 \Big) \\
&= \ov{\A}_1 \circ  \PP \D_5 \du \ov{\A}_1 \circ (\ov{\PP \D}_6 \cup \ov{\PP \E}_6) \du 
\Big( \Delta \ov{\PP \D}_7 \cup \Delta \ov{\PP \E}_7 \Big) 
\end{align*}
where the last equality follows from \cite{BM13} (cf. Lemma \ref{cl}, statement \ref{D5cl}). By Propositions \ref{D6_Condition_prp} and \ref{E6_Condition_prp}, the sections 
\[ \pi_2^*\us_{\PP \D_6} \oplus \Q : \ov{\ov{\A}_1 \circ \PP \D}_5 \lra \pi_2^*\UL_{\PP \D_6} \oplus \W_{n,m,6}^1, 
~~\pi_2^*\us_{\PP \E_6} \oplus \Q : \ov{\ov{\A}_1 \circ \PP \D}_5 \lra \pi_2^*\UL_{\PP \E_6} \oplus \W_{n,m,6}^1 \]
vanishes transversely on $\A_1 \circ \PP \D_6$ and $\A_1 \circ \PP \E_6$ respectively. 
Moreover, they do not vanish on $\A_1 \circ \PP \E_6$ and $\A_1 \circ \PP \D_6$ respectively.  
By Corollary \ref{a1_pdk_mult_is_2_f12_neq_0} and \ref{psi_pe6_and_pd6_section_vanishes_order_one_around_pe7}
the contribution of the section  $\pi_2^*\us_{\PP \D_6} \oplus \Q$ 
to the Euler class from the points of $\Delta\PP \D_7$ and $\Delta \PP \E_7$ are $2$ and $1$ respectively. 
By Corollary \ref{psi_pe6_and_pd6_section_vanishes_order_one_around_pe7}, the contribution of the section    
$\pi_2^*\us_{\PP \E_6} \oplus \Q$ from the points of $\Delta \PP \E_7$ is $1$; moreover it does not vanish on $\Delta \PP \D_7$. 
Hence 
\begin{align*}
\Big\langle e(\pi_2^*\UL_{\PP \D_6} \oplus \W_{n,m,6}^1), ~[\ov{\ov{\A}_1 \circ \PP \D}_5] \Big\rangle & = \Num(\A_1\PP\D_6,n,m) + 2 \Num(\PP \D_7, n, m)+\Num(\PP \E_7, n, m), \\
\Big\langle e(\pi_2^*\UL_{\PP \E_6} \oplus \W_{n,m,6}^1), ~[\ov{\ov{\A}_1 \circ \PP \D}_5] \Big\rangle & = \Num(\A_1\PP\E_6,n,m) + \Num(\PP \E_7, n, m)
\end{align*}
which prove equations \eqref{algopd6a1} and \eqref{algope6a1}.  \qed

\appendix

\section{Low degree checks} 

\ni \textbf{Verification of the number $\Num(\A_1\A_1,0) = 3(d-2)(d-1)(3d^2-3d-11)$:} \\
\ni $d=2:$ The only way a conic can have $2$ nodes is if it is 
a double line. There are no double lines through $3$ generic points.\\
\ni $d=3:$ The only way a cubic can have $2$ nodes is if it 
breaks into a line and a conic. Hence the number of cubics with $2$
unordered nodes is $ \binom{7}{2}= 21.$ \\

\ni \textbf{Verification of the number $\Num(\A_1\A_2,0) = 12(d-3)(3d^3-6d^2-11d+18)$:} \\
\ni $d=3:$ There are no cubics with one node and one cusp.\\

\ni \textbf{Verification of the number $\Num(\A_1\A_3,0) = 6(d-3)(25d^3-71d^2-122d+280)$:} \\
\ni $d=3:$ There are no cubics with one node and one tacnode.\\ 
\ni $d=4:$ There are two possibilities here. The curve
could break into a line and a cubic. The number of lines through
a given point and tangent to a fixed cubic is $6$. The number of cubics
through through $8$ points, tangent to a given line is $4$. Hence the total 
number of quartics with one node and one tacnode that breaks into a line
and a cubic is $$\binom{10}{1}\times 6+ \binom{10}{2}\times 4 = 240.$$
The number of genus zero quartics with one node 
and one tacnode is $1296$ (cf. \cite{g0pr}, pp. $91$). 
Hence the total number of quartics with
one node and one tacnode is $1536$. \\ 

\ni \textbf{Verification of the number $\Num(\A_1 \D_4,n)$ and $\Num(\A_1\PP \D_4,n,1)$:} \\ 
\ni $d=4:$ These numbers can be verified  by direct geometric means for all values 
of $n$ in the case of quartics (cf. \cite{BM_Detail}). 


\section{Standard facts about Chern classes and projectivized bundle} 

\begin{lmm}{\bf (\cite{MiSt}, Theorem 14.10)}
\label{total_chern_class_of_tpn}
Let $\P^n$ be the $n$-dimensional complex projective space and $\gamma \lra \P^n$ the tautological line bundle. 
Then the total Chern class $c(T\P^n)$ is given by $~c(T\P^n) = (1+ c_1(\gamma^*))^{n+1}.$
\end{lmm}

\begin{lmm}{\bf (\cite{BoTu}, pp. 270)}
\label{cohomology_ring_of_pv}
Let $V\lra M$ be a complex vector bundle, of rank $k$, over a smooth manifold $M$ and $\pi: \P V \lra M$ the \textit{projectivization} of 
$V$. Let $\G\lra \P V$ be the tautological line bundle over $\P V$ and $\lm = c_1(\G^*)$. 
There is a linear isomorphism
\begin{align}
H^*(\P V; \mathbb{Z}) & \cong  H^*(M;\mathbb{Z}) \otimes H^*(\P^{k-1};\mathbb{Z})
\end{align}
and an isomorphism of rings
\begin{align}
H^*(\P V; \mathbb{Z})& \cong  \frac{H^*(M;\mathbb{Z})[\lm]}{\lan \lm^k + \lm^{k-1} \pi^*c_1(V) + \lm^{k-2} \pi^*c_2(V) + \ldots +\pi^*c_{k}(V)\ran}.
\end{align}
In particular, if $\omega \in H^{*}(M; \mathbb{Z}) $ is a top cohomology class then 
\bgd
\lan \pi^*(\omega) \lm^{k-1}, [\P V] \ran = \lan \omega, [M] \ran,
\edd
i.e., $\lm^{k-1}$ is a \textit{cohomology extension} of the fibre.
\end{lmm}

\bibliographystyle{siam}
\bibliography{Myref_bib.bib}

\begin{thebibliography}{1}

\bibitem{BM13}
{\sc S.~Basu and R.~Mukherjee}, {\em Enumeration of curves with one singular
  point}.
\newblock available at \url{http://arxiv.org/abs/1308.2902}.

\bibitem{BM_Detail}
\leavevmode\vrule height 2pt depth -1.6pt width 23pt, {\em Enumeration of
  curves with singularities: Further details}.
\newblock available at \url{https://www.sites.google.com/site/ritwik371/home}.

\bibitem{BoTu}
{\sc R.~Bott and L.~W. Tu}, {\em Differential forms in algebraic topology},
  vol.~82 of Graduate Texts in Mathematics, Springer-Verlag, New York, 1982.

\bibitem{Kaz}
{\sc M.~{\`E}. Kazarian}, {\em Multisingularities, cobordisms, and enumerative
  geometry}, Uspekhi Mat. Nauk, 58 (2003), pp.~29--88.

\bibitem{MiSt}
{\sc J.~W. Milnor and J.~D. Stasheff}, {\em Characteristic classes}, Princeton
  University Press, Princeton, N. J., 1974.
\newblock Annals of Mathematics Studies, No. 76.

\bibitem{Z1}
{\sc A.~Zinger}, {\em Counting plane rational curves: old and new approaches}.
\newblock available at \url{http://arxiv.org/abs/math/0507105}.

\bibitem{g0pr}
{\sc A.~Zinger}, {\em Counting rational curves of arbitrary shape in projective
  spaces}, Geom. Topol., 9 (2005), pp.~571--697 (electronic).

\end{thebibliography}

\vspace*{0.4cm}

\hf {\small D}{\scriptsize EPARTMENT OF }{\small M}{\scriptsize ATHEMATICS, }{\small RKM V}{\scriptsize IVEKANANDA }{\small U}{\scriptsize NIVERSITY, }{\small H}{\scriptsize OWRAH, }{\small WB} {\footnotesize 711202, }{\small INDIA}\\
\hf{\it E-mail address} : \texttt{somnath@maths.rkmvu.ac.in}\\

\hf {\small D}{\scriptsize EPARTMENT OF }{\small M}{\scriptsize ATHEMATICS, }{\small I}{\scriptsize NSTITUTE OF }{\small M}{\scriptsize ATHEMATICAL }{\small S}{\scriptsize CIENCES, }{\small C}{\scriptsize HENNAI }{\footnotesize 600113, }{\small INDIA}\\
\hf{\it E-mail address} : \texttt{ritwikm@imsc.res.in}\\[0.2cm]

\end{document}